%% file: LB_NCVX_affineConstraint.tex
\documentclass[moor,sglanonrev]{informs4}
 
\OneAndAHalfSpacedXI % Current default line spacing
\usepackage{amsmath,amssymb,amsfonts}

\usepackage[hidelinks]{hyperref}
\usepackage[sort&compress]{natbib}
\bibpunct[, ]{[}{]}{,}{n}{}{,}%
\def\bibfont{\small}%

\usepackage{url}
\usepackage{lineno}
\usepackage{tikz}

\RequirePackage{endnotes}
\usepackage[small, margin=1cm]{caption}

\usepackage{appendix}
\usepackage{color}
\definecolor{strcolor}{rgb}{0.6, 0.2, 0.6}
\definecolor{commentcolor}{rgb}{0.3125, 0.5, 0.3125}
\definecolor{keycol}{rgb}{0, 0, 1}

\usepackage{rotating}
\usepackage{fancyvrb}
\usepackage{tcolorbox}
\usepackage{bbm}

\usepackage{listings}

\TheoremsNumberedThrough     % Preferred (Theorem 1, Lemma 1, Theorem 2)
\ECRepeatTheorems
\JOURNAL{MATHEMATICS OF OPERATIONS RESEARCH }

\EquationsNumberedThrough    % Default: (1), (2), ...
\gdef\AQ#1{}
\gdef\CQ#1{}

\MANUSCRIPTNO{MOOR-0001-2024.00}

\input{macros.tex}

\begin{document}
	\def\COPYRIGHTHOLDER{INFORMS}%
	\def\COPYRIGHTYEAR{2017}%
	\def\DOI{\fontsize{7.5}{9.5}\selectfont\sf\bfseries\noindent https://doi.org/10.1287/opre.2017.1714\CQ{Word count = 9740}}
	
	\RUNAUTHOR{Liu, Lin, and Xu} 
	
	\RUNTITLE{Lower Complexity Bounds of FOMs for Affinely Constrained Composite Optimization}
	
	\TITLE{Lower Complexity Bounds of First-order Methods for Affinely Constrained Composite Non-convex Problems}
	
	\ARTICLEAUTHORS{%
		\AUTHOR{Wei Liu}
		%,\textsuperscript{a} Second Author,\textsuperscript{b} Third Author,\textsuperscript{c} Fourth Author,\textsuperscript{c}
		
		\AFF{Department of Mathematical Sciences, Rensselaer Polytechnic Institute, Troy, NY, \EMAIL{liuw16@rpi.edu}}
		%\textsuperscript{b}School of Industrial Engineering, Good College, Collegeville, Maine 01234 \EMAIL{secauth@goodcoll.edu}; 
		%\textsuperscript{c}Their Common Affiliation \EMAIL{thauth@anywhere.edu, fourauth@anywhere.edu}
		
		%mirko.janc@informs.org
		\AUTHOR{Qihang Lin}
		
		\AFF{Department of Business Analytics, University of Iowa, Iowa City, IA, \EMAIL{qihanglin@uiowa.edu}}
		
		\AUTHOR{Yangyang Xu}
		
		\AFF{Department of Mathematical Sciences, Rensselaer Polytechnic Institute, Troy, NY,
		\EMAIL{ xuy21@rpi.edu}}
	}
	
	%\comm{Would it better to change the title to be ``Iteration Complexity of First-order Methods for Affinely Constrained Composite Non-convex Non-smooth Problems: Lower Bound and Upper Bound''.}
	
	%\comm{When we say near stationary point, shall we use ``(near) stationary point''?}
		
\ABSTRACT{
Many recent studies on first-order methods (FOMs) focus on \emph{composite non-convex non-smooth} optimization with linear and/or nonlinear function constraints. Upper (or worst-case) complexity bounds have been established for these methods. However, little can be claimed about their optimality as no lower bound is known, except for a few special \emph{smooth non-convex} cases. In this paper, we make the first attempt to establish lower complexity bounds of FOMs for solving a class of composite non-convex non-smooth optimization with linear constraints. Assuming two different first-order oracles, we establish lower complexity bounds of FOMs to produce a (near) $\epsilon$-stationary point of a problem (and its reformulation) in the considered problem class, for any given tolerance $\epsilon>0$.  
{Our lower bounds indicate that the existence of a non-smooth convex regularizer can evidently increase the difficulty of an affinely constrained regularized problem over its nonregularized counterpart.}
In addition, we show that our lower bound {of FOMs with the second oracle} is tight, with a difference of up to a logarithmic factor from an upper complexity bound established in the extended arXiv version of this paper; see Liu et al.~\cite{liu2023first}. 
}

\SUBJECTCLASS{MSC2000: 90C26, 90C06, 90C60, 49M37, 68Q25, 65Y20; \\  OR/MS subject classification: Computational complexity.}

\AREAOFREVIEW{Continuous optimization.}
	
	
\maketitle

	\section{Introduction}

	First-order methods (FOMs) have attracted increasing attention because of their efficiency in solving large-scale problems arising from machine learning and other areas. 
	The recent studies on FOMs have focused on non-convex problems, and one of the actively studied topics is the oracle complexity for finding a near-stationary point under various assumptions. In this paper, we explore this topic for problems with a composite non-convex non-smooth objective function and linear equality constraints, formulated as %\YX{Let's change $F$ to $F_0$ here and all other places for (P), as $F$ is used in~\eqref{eq:model-spli}.}
	\begin{equation}\label{eq:model}
		\begin{aligned}
			\min_{ \vx\in\RR^d} \,\,& F_0( \vx):= f_0( \vx) + g( \vx), \ \st  \vA \vx +  \vb=\mathbf{0}.
			%\\& \textbf{h}( \vx):=[h_1( \vx),\ldots, h_m( \vx)]\leq 0,
		\end{aligned}
	\tag{P}
	\end{equation}
	Here,  $\vA\in\RR^{n\times d}$, $\vb\in\RR^{n}$, $f_0:\mathbb{R}^d \rightarrow \mathbb{R}$ is $L_f$-smooth and potentially non-convex, {i.e., 	$
		\left\|\nabla f_0\left(\vx\right)- \nabla f_0\left(\vx'\right)\right\| \leq L_f\left\|\vx-\vx'\right\|, \forall\, \vx,\,\vx'\in \mathbb{R}^d
		$}, and $g: \mathbb{R}^d \rightarrow \mathbb{R} \cup\{+\infty\}$ is a proper lower semicontinuous convex function but potentially non-smooth.

	Due to non-convexity, computing or even approximating a global optimal solution for problem~\eqref{eq:model} is intractable in general; see Murty and Kabadi~\cite{murty1985some}. Hence, we focus on using an FOM to find a (near) $\epsilon$-stationary point of problem~\eqref{eq:model} for a given tolerance $\epsilon > 0$; see Definition~\ref{def:eps-pt-P}.  An FOM finds a (near) $\epsilon$-stationary point by querying information from some oracles, which typically dominates the runtime of the method, so its efficiency can be measured by the number of oracles it queries, 
	which is defined as the method's \emph{oracle complexity}. The goal of this paper is to establish lower bounds for the oracle complexity of a class of FOMs to find a (near) $\epsilon$-stationary point of problem~\eqref{eq:model} that satisfies certain properties. %and, additionally, to present an FOM that can nearly achieve this lower complexity bound. 

 \subsection{Problem Class, Oracles, and Algorithm Classes} We will study the lower complexity bounds of FOMs for solving the following problem class.

% \begin{prob}\label{prob:class-1}
% Any problem instance in this class is in the form of \eqref{eq:model} where $\mathrm{relint}(\dom(g)) \neq \emptyset$, $\inf_{\vx} F_0(\vx) > -\infty$ and $f_0$ is $L_f$-smooth,  i.e., 
%			$
%			\left\|\nabla f_0\left(\vx\right)- \nabla f_0\left(\vx'\right)\right\| \leq L_f\left\|\vx-\vx'\right\|, \forall\, \vx,\,\vx'\in \mathbb{R}^d.
%			$ 
%\end{prob}

\begin{prob}\label{prob:class-2} Any problem instance in this class is in the form of \eqref{eq:model} and satisfies the conditions: (i) $\inf_{\vx} F_0(\vx) > -\infty$; (ii) $f_0$ is $L_f$-smooth;
	%$\mathrm{relint}(\dom(g)) \neq \emptyset$,  and  and   
	(iii) $
	g(\vx)= \bar{g}(\bar{\vA} \vx +\bar{\vb}),
	$
	where $\bar{\vA}\in \RR^{\bar{n}\times d}$, $\bar{\vb}\in\RR^{\bar{n}}$, %\blue{$\bar{n}\leq d$} 
	and $\bar{g}: \mathbb{R}^{\bar n} \rightarrow \mathbb{R}\cup\{+\infty\}$ is a proper lower semicontinuous convex function with $\mathrm{relint}(\dom(\bar{g}))\neq\emptyset$ and $\dom(\bar{g})$ not a singleton. %but is not single-valued. %in the form of \eqref{eq:model} and satisfies
	%\begin{itemize}
	%			\item[\textnormal{(a)}]  $\nabla f_0$ is $L_f$-Lipschitz continuous,  i.e., 
	%			$
	%			\left\|\nabla f_0\left(\vx\right)- \nabla f_0\left(\vx'\right)\right\| \leq L_f\left\|\vx-\vx'\right\|, \forall\, \vx,\,\vx'\in \mathbb{R}^d.
	%			$ 
	%			\item[\textnormal{(b)}] $\inf_{\vx} F_0(\vx) > -\infty$.
	%			\item[\textnormal{(c)}]	$
	%			g(\vx)= \bar{g}(\bar{\vA} \vx +\bar{\vb}),
	%			$
	%			where $\bar{\vA}\in \RR^{\bar{n}\times d}$, $\bar{\vb}\in\RR^{\bar{n}}$ and $\bar{g}: \mathbb{R}^{\bar n} \rightarrow \mathbb{R} \cup\{+\infty\}$ is a proper lower semicontinuous convex function. %but potentially non-smooth.
	%\end{itemize}
\end{prob}	
 
Notice that we require   $\mathrm{relint}(\dom(\bar g)) \neq \emptyset$ in Problem Class~\ref{prob:class-2} to ensure the existence of a KKT point. 
%The requirement that $g$ is not single valued 
In addition, with the requirement $\dom(\bar{g})$ not to be a singleton, we also 
exclude the case where $\bar g$ is the indicator function on $\{\vzero\}$, for which the resulting problem looks like an affinely constrained \emph{composite non-smooth} problem but is actually a \emph{smooth} problem. 
%\blue{If we set $\bar{\vA}=\vI$ and $\bar{\vb}=\vzero$, then our problem class degenerates to a general problem. }
 Undoubtedly, an algorithm's oracle complexity depends on what oracle information the algorithm can utilize and what operations it can perform by using the oracle. 
Since we focus on FOMs, we assume that a first-order oracle is accessible and each generated iterate is a certain combination of the oracle information. %Due to the structural difference between the two problem classes, 
We assume different first-order oracles and consider two different classes of FOMs toward solving Problem Class~\ref{prob:class-2}.		

%Respectively for Problem Classes~\ref{prob:class-1} and \ref{prob:class-2}, 
We assume the following oracles accessible at any given input. 
%\comm{Do we have to introduce $\vb$ in the oracles?}
\begin{align}
&\mathrm{ORACLE}_1(\vx, \vz, \eta)\mbox{ returns } \big(\nabla f_0(\vx), \vA\vx, \vA\zz \vz, \prox_{\eta g}(\vx)\big), \forall\, \vx\in\RR^d, \vz\in\RR^n, \eta >0, \label{oracle-1}\\
&\mathrm{ORACLE}_2(\vx, \vy, \vz, \eta)\mbox{ returns } \big(\nabla f_0(\vx), \bar\vA\vx, \vA\vx, \bar\vA\zz \vy, \vA\zz \vz, \prox_{\eta \bar g}(\vy)\big),\label{oracle-2}\\
&\hspace{8cm} \forall\, \vx\in\RR^d, \vy\in\RR^{\bar n}, \vz\in\RR^n, \eta >0,\nonumber
\end{align}
where the proximal mapping is defined as
\begin{eqnarray}
	\label{eq:proxg}
	\prox_{\eta g}(\vx):=\arg\min_{ \vx'} \left\{ \textstyle g(\vx')+\frac{1}{2\eta}\|\vx'-\vx\|^2 \right\}, \forall\, \vx\in \RR^d, \eta > 0.
\end{eqnarray}
When $g(\vx)= \bar{g}(\bar{\vA} \vx +\bar{\vb})$,  it can be significantly more difficult to evaluate $\prox_{\eta g}(\cdot)$ than $\prox_{\eta \bar g}(\cdot)$. In this case, $\mathrm{ORACLE}_1$ is stronger than $\mathrm{ORACLE}_2$. Notice that both oracles can return a bundle of vectors. However, if only one returned vector such as $\nabla f_0(\vx)$ is needed in an update, one can still call the oracle, and we still count the call once. %even if not all returned vectors are used in one update. 
We will study lower complexity bounds of FOMs that belong to the following algorithm classes for solving Problem Class~\ref{prob:class-2}. An FOM in Algorithm Class~\ref{algo:class-1} uses $\mathrm{ORACLE}_1$ in \eqref{oracle-1}, while an FOM in Algorithm Class~\ref{algo:class-2} uses $\mathrm{ORACLE}_2$ in \eqref{oracle-2}.
\begin{algo}\label{algo:class-1}
Given $\vx^{(0)}\in \RR^d$, the iterates are generated such that %sequence $\left\{\vx^{(t)}\right\}_{t=1}^{\infty}$ satisfies that,  
for any $t\geq 1$, $\vx^{(t)} \in  
		\mathbf{span}\left( \big\{\vxi^{(t)}, \prox_{\eta_t g}(\vxi^{(t)}) \big\} \right)$, where $\eta_t > 0$ and $\vxi^{(t)}\in \mathbf{span}\left(\left\{\vx^{(s)}, \nabla f_0(\vx^{(s)}),  \vA\zz\vb, \vA\zz \vA\vx^{(s)}\right\}_{s=0}^{t-1}\right).$
\end{algo}

\begin{algo}\label{algo:class-2}
Given $(\vx^{(0)}, \vy^{(0)})\in \RR^d \times \RR^{\bar n}$, the sequence $\left\{(\vx^{(t)},\vy^{(t)})\right\}_{t=1}^{\infty}$ is generated such that,  for all $t\geq 1$, 
		\begin{equation*}
			\begin{aligned}
				&\vx^{(t)}\in \mathbf{span}\left(\left\{\vx^{(s)}, \nabla f_0(\vx^{(s)}), \vA\zz\vb,\vA\zz \vA\vx^{(s)}, \bar\vA\zz\bar\vb,  \bar\vA\zz \bar\vA\vx^{(s)},\bar\vA\zz\vy^{(s)}\right\}_{s=0}^{t-1}\right),\\
				&\vy^{(t)} \in %\mathbf{span}(\mathcal{M}):=\mathbf{span}\left\{ \left\{\nabla f_i\left(\vx_i^{(0)}\right)\right\}^{m}_{i=1}, \A\zz A\vx^{(0)}, \ldots, \left\{\nabla f_i\left(\vx_i^{(t-1)}\right)\right\}^{m}_{i=1}, A\zz A\vx^{(t-1)}, \{\nabla g(\vx^{(i)})\}_{i=1}^{t-1}\right\}. 
				\mathbf{span}\left(\big\{\vxi^{(t)},\prox_{\eta_t \bar{g}}(\vxi^{(t)}) \big\}\right), \text{ where } \eta_t >0 \text{ and } %\\ &
				\vxi^{(t)}\in \mathbf{span}\left(\left\{\bar\vb, \vy^{(s)}, \bar\vA\bar\vA\zz\vy^{(s)}, \bar\vA\vx^{(s)}\right\}_{s=0}^{t-1}\right).
			\end{aligned}		
		\end{equation*}
\end{algo}

\noindent\textbf{Examples.}~Several existing FOMs while applied to solve \eqref{eq:model} belong to Algorithm Class~\ref{algo:class-1} (e.g., Kong et al.~\cite{kong2019complexity} , Melo et al.~\cite{melo2020iteration1}, Xu~\cite{xu2021first}, Liu et al.~\cite{liu2022linearly}, Liu et al.~\cite{lin2022complexity},  Zhang and Luo~\cite{zhang2022global}). On instances in Problem Class~\ref{prob:class-2}, Algorithm Class~\ref{algo:class-2} includes FOMs that are designed based on quadratic penalty or augmented Lagrangian method (ALM) and linearized alternating direction method of multipliers (ADMM) (see Goncalves et al.~\cite{goncalves2017convergence}, Xu~\cite{xu2017accelerated}) for solving the equivalent Splitting Problem %\footnote{\revise{Suppose that $(\vx^*,\vy^*)$ is an $\epsilon$-stationary point of problem \eqref{eq:model-spli} and $\bar{g}$ is Lipschitz continuous, then $\vx^*$ is $O(\epsilon)$ close to an $\epsilon$-stationary point of problem \eqref{eq:model}.}}
		\begin{equation}
				\label{eq:model-spli}
		\begin{aligned}		
			\min_{ \vx\in\RR^d, \vy\in\RR^{\bar n}} &~F(\vx,\vy):=f_0(\vx) + \bar{g}( \vy),\ 
			 \st~\vA \vx + \textbf{b} =0,\,\, \vy = \bar{\vA} \vx + \bar{\textbf{b}}.
		\end{aligned}		\tag{SP}
		\end{equation}		

%for any $\vx\in\mathbb{R}^d$ and $\eta>0$. We say  $\vx^{(t)}$ is generated by the \emph{$t$-th iteration}\footnote{Without further specification, one iteration of an FOM will call the oracle once. We will specify the type of iteration if this does not hold, e.g., by outer- or inner-iteration.} of the algorithm. 

	{By the Lagrangian function, the problems in \eqref{eq:model} and \eqref{eq:model-spli} can be formulated into bilinear saddle-point problems. Hence, existing FOMs on solving minimax-structured optimization can be applied. Examples of FOMs include the proximal dual implicit accelerated gradient algorithm (Thekumparampil~\cite{thekumparampil2019efficient}), the gradient descent-ascent method (e.g. Lin et al.~\cite{ lin2020gradient}, Xu et al.~\cite{xu2023unified}) and the gradient descent maximization method (e.g. Lu et al.~\cite{lu2020hybrid}, Lin et al.~\cite{lin2020gradient}, Xu~\cite{xu2024d-GDMax}) if an FOM is applied to perform the dual maximization. When applied to the saddle-point reformulation of \eqref{eq:model} (resp. \eqref{eq:model-spli}), these FOMs belong to Algorithm Class~\ref{algo:class-1} (resp. Algorithm Class~\ref{algo:class-2}) if the dual iterate is initialized properly. In particular, let us focus on problem~\eqref{eq:model}. Let $\vz^{(t)}$ for $t=0,1,\dots$ be the Lagrangian multiplier (i.e., dual iterate) to the linear constraint $\vA\vx +\vb=\vzero$ 
	generated at the $t$-th iteration of an FOM for a minimax problem. Then, the aforementioned FOMs belong to Algorithm Class~\ref{algo:class-1} %or Algorithm Class~\ref{algo:class-2} 
	if they are initialized with  $\vz^{(0)}= \vzero$ or $\vz^{(0)}=\vA\vx^{(0)}+\vb$. 
	Although $\vz^{(t)}$ is not defined in Algorithm Class~\ref{algo:class-1}, %or Algorithm Class~\ref{algo:class-2}, 
	it can be proved %by induction 
	that these algorithms can be implemented equivalently by updating $\vx^{(t)}$ and  $\widehat\vx^{(t)}:=\vA^\top\vz^{(t)}$ %(where $\widetilde\vA$ is the coefficient matrix) 
	without calculating $\vz^{(t)}$ explicitly, and $\vx^{(t)}$ and  $\widehat\vx^{(t)}$ are generated by following the condition of Algorithm Class~\ref{algo:class-1}. %(resp. Algorithm Class~\ref{algo:class-2}). 
}

It should be noted that FOMs in Algorithm Classes~\ref{algo:class-1} and \ref{algo:class-2} are allowed to access $\prox_{\eta g}$ or $\prox_{\eta \bar g}$  
but not the projection onto the affine set $\{\vx\in\mathbb{R}^d~|~\vA\vx +  \vb=\mathbf{0}\}$. 
It is also worth noting that our considered Algorithm Classes do not include subgradient methods. In general, computing $\prox_{\eta g}$ is more difficult than computing a subgradient of $g$, and finding an $\epsilon$-stationary point of non-smooth optimization is impossible for subgradient methods with finite iterations; see Zhang et al.~\cite{zhang2020complexity}. %In fact, if we disallow computing the proximal mapping of $g$ but instead give the algorithm access to its subgradients, the algorithm is limited to the class of subgradient methods for non-smooth optimization, in which case finding an $\epsilon$-stationary point is  impossible with finite oracle complexity~\cite{zhang2020complexity}. %\YX{Even if $g$ is convex?} \comm{It would be more than $O(\epsilon^{-4})$ iterations. This part is rewritten by Qihang.}
%at least $O(\epsilon^{-4})$~\cite{arjevani2022lower}. 
%{\color{red}Qihang: I don't think~\cite{arjevani2022lower} has this result. They only consider stochastic case. }
{Moreover, notice that though an FOM in Algorithm Class~\ref{algo:class-2} can utilize the compositional structure of $g$ specified in Condition (iii) of Problem Class~\ref{prob:class-2}, one FOM in Algorithm Class~\ref{algo:class-1} does not assume such a structure and can be applied to the case with a general convex function~$g$. However, we still restrict our consideration to a structured $g$ while establishing lower complexity bounds of FOMs in Algorithm Class~\ref{algo:class-1}, because it enables us to \emph{explicitly} quantify, via the condition number $\kappa([\bar\vA; \vA])$ defined in \eqref{eq:joint_condition_number} below, the difficulty caused by the interaction of $g$ and the affine constraint $\vA\vx+\vb=\vzero$.}

\subsection{Background and Main Question} FOMs have been developed with theoretically proved oracle complexity for finding a (near) $\epsilon$-stationary point of problem~\eqref{eq:model} or a more general problem.  We refer readers to Bian and Chen~\cite{bian2015linearly}, Haeser et al.~\cite{haeser2019optimality}, Kong et al.~\cite{kong2019complexity}, Jiang et al.~\cite{jiang2019structured}, Melo et al.~\cite{melo2020iteration1}, O'Neill and Wright~\cite{o2021log} for problems with affine constraints,  and Cartis et al.~\cite{cartis2014complexity}, Cartis et al.~\cite{cartis2017corrigendum}, Sahin et al.~\cite{sahin2019inexact}, Xie and Wright~\cite{xie2021complexity}, Grapiglia and Yuan~\cite{grapiglia2021complexity}, Lin et al.~\cite{lin2022complexity}, Xiao et al.~\cite{xiao2023dissolving}, Kong et al.~\cite{kong2023iteration}, Liu et al.~\cite{liu2022inexact}, Dahal et al.~\cite{dahal2023damped}, and Huang and Lin~\cite{huang2023single} for problems with nonlinear constraints. These results are stated as  \emph{upper bounds} for the maximum number of oracles those algorithms require to reach a (near) $\epsilon$-stationary point. In special cases of problem~\eqref{eq:model}, some existing algorithms' oracle complexity is known to be non-improvable (also called optimal) because it matches, up to constant factors, a theoretical \emph{lower bound} of oracle complexity (see Nemirovskij and Yudin~\cite{nemirovskij1983problem}), which is the minimum number of oracles an algorithm  %satisfying \revise{a certain linear span assumption} 
needs to find a (near) %\revise{(near)} 
$\epsilon$-stationary point. Below we discuss two special cases of \eqref{eq:model}. %are as follows.

\vspace{0.1cm}
\noindent\textbf{Composite non-convex optimization.}~When the linear constraint $\vA \vx +  \vb=\mathbf{0}$ does not exist,  problem~\eqref{eq:model} becomes a composite non-convex optimization problem 
	\begin{equation}
		\label{eq:composite}
	\min_{\vx}F_0(\vx) = f_0(\vx)+g(\vx). 
	\end{equation}
When $g\equiv 0$, a lower complexity bound of $O(L_f \Delta_{F_0}\epsilon^{-2})$ is given by Carmon et al.~\cite{carmon2020lower} to produce an $\epsilon$-stationary point $\bar\vx$ by an FOM, i.e., $\|\nabla f_0(\bar\vx)\| \le \epsilon$, where $\Delta_{F_0}:= F_0(\vx^{(0)}) - \inf_\vx F_0(\vx)$. On the other side, the steepest gradient descent can reach the same-order complexity bound, see Cartis et al.~\cite{Cartis2012how}. % find an $\epsilon$-stationary point by $O(L_f\epsilon^{-2})$ first-order oracles~\cite{Cartis2012how}. 
Hence, the  lower bound of $O(L_f \Delta_{F_0} \epsilon^{-2})$ is tight. Since smooth non-convex optimization is a subclass of composite non-convex optimization \eqref{eq:composite} if no restriction is imposed on $g$, $O(L_f \Delta_{F_0} \epsilon^{-2})$ is a valid lower bound for FOMs on solving \eqref{eq:composite}. It is well-known that the proximal gradient method can find an $\epsilon$-stationary point of~\eqref{eq:composite} within  $O(L_f \Delta_{F_0} \epsilon^{-2})$ iterations; see Lan~\cite{lan2020first} for example. This implies that $O(L_f \Delta_{F_0} \epsilon^{-2})$ is also a tight lower bound of FOMs for solving composite non-convex optimization, and we can make the following claim. 
\begin{claim}\label{claim-1}
If $\prox_{\eta g}$ is accessible for any $\eta>0$, the convex regularizer $g$ in composite non-convex optimization \eqref{eq:composite} does not increase the lower complexity bound for FOMs or make the problem harder.  
\end{claim}

\vspace{0.1cm}
\noindent\textbf{Affinely constrained smooth non-convex optimization.}~When $g\equiv0$, problem~\eqref{eq:model} becomes a linear equality-constrained smooth non-convex optimization problem %with linear equality constraints, i.e., 
	\begin{align}
		\label{eq:affinenonconvex}
			\min_\vx&~f_0(\vx), \ \text{ s.t. } \vA\vx+\vb=\mathbf{0}.
	\end{align}
	
	It is known that, if exact projection onto the feasible set $\{\vx~|~\vA\vx +  \vb=\mathbf{0}\}$ is allowed, the projected gradient method can find an $\epsilon$-stationary point in $O(L_f \Delta_{f_0} \epsilon^{-2})$ iterations, where $\Delta_{f_0}:= f_0(\vx^{(0)}) - \inf_\vx f_0(\vx)$. When exact projection is prohibited, %by \revise{the linear span assumption}, 
	one can perform inexact projection through the matrix-vector multiplications with $\vA$ and $\vA^{\top}$, which are allowed in $\mathrm{ORACLE}_1$ and $\mathrm{ORACLE}_2$. %by \revise{the linear span assumption}. 
	This way, with $O(\kappa(\vA)\log(\epsilon^{-1}))$ multiplications, one can project any point to the feasible set with a $\mathrm{poly}(\epsilon)$ error. Here, {$\mathrm{poly}(\epsilon)$ denotes a polynomial of $\epsilon$,} %\YX{How about $\mathrm{poly}(\epsilon)$ error?} \revise{It would be better}. 
	$\kappa(\vA)$ is the condition number of $\vA$ defined as
	\begin{align*}
		\kappa(\vA):= \sqrt{\frac{\lambda_{\max}(\vA\vA\zz)}{\lambda^+_{\min}(\vA\vA\zz)}},
	\end{align*}
	where $\lambda^+_{\min}(\vA\vA\zz)$  and $\lambda_{\max}(\vA\vA\zz)$ are the smallest positive and the largest eigenvalues of $\vA\vA\zz$, respectively. Using this inexact projection as a subroutine, one can easily develop an inexact projected gradient method to \eqref{eq:affinenonconvex} with oracle complexity of  $O(\kappa(\vA)\log(\epsilon^{-1})L_f \Delta_{f_0} \epsilon^{-2})$. This complexity matches the lower bound $O(\kappa(\vA)L_f \Delta_{f_0} \epsilon^{-2})$ in Sun and Hong~\cite{sun2019distributed} up to a logarithmic factor and thus is nearly optimal.

%These two examples suggest that a lower bound of oracle complexity is valuable as it informs algorithm designers which algorithms can be potentially improved for higher efficiency and which are non-improvable without additional assumptions. 

By a similar discussion for composite non-convex optimization, we notice that $O(\kappa(\vA)L_f \Delta_{F_0} \epsilon^{-2})$ is also a valid lower complexity bound of FOMs for solving \eqref{eq:model} in Problem Class~\ref{prob:class-2}, as it permits the case of $g\equiv 0$ or $\bar g\equiv 0$. However, no existing FOM can (nearly) reach this complexity. Hence, it is unknown if this lower bound of $O(\kappa(\vA)L_f \Delta_{F_0} \epsilon^{-2})$ is (nearly) tight or not for solving affinely constrained composite non-convex optimization~\eqref{eq:model}. This leads to a natural question as follows. 

\begin{tcolorbox}
\textbf{Question:} For Problem Class~\ref{prob:class-2}, will the existence of a \emph{non-smooth} convex regularizer $g$ make the problem harder when $\prox_{\eta g}$ is accessible for any $\eta>0$? If yes, what will a higher lower complexity bound be for FOMs in Algorithm Class~\ref{algo:class-1}?		
	\end{tcolorbox}
We will provide an affirmative answer to the question above. In short, for FOMs in Algorithm Class~\ref{algo:class-1} on solving Problem Class~\ref{prob:class-2}, \textbf{{our new lower complexity bound can be significantly higher than $O(\kappa(\vA)L_f \Delta_{F_0} \epsilon^{-2})$}}. This is different from Claim~\ref{claim-1}. But we cannot show the tightness of the new lower bound, and we leave it to the future work. Nevertheless, we also establish a lower bound of FOMs in Algorithm Class~\ref{algo:class-2} on solving Problem Class~\ref{prob:class-2}, for which we are able to show its tightness.

	\subsection{Contributions}
		
		Our first major contribution is to establish a lower bound of the oracle complexity of FOMs in Algorithm Class~\ref{algo:class-1} for finding a (near) $\epsilon$-stationary point of problem instances in Problem Class~\ref{prob:class-2}.  
		This is achieved by building a hard instance that belongs to Problem Class~\ref{prob:class-2}. %and thus also to Problem Class~\ref{prob:class-1}. 
		We show that any FOM in Algorithm Class~\ref{algo:class-1} needs at least $O({\kappa([\bar{\vA};\vA]) L_f \Delta_{F_0}} \epsilon^{-2})$ calls to $\mathrm{ORACLE}_1$ to find a (near) $\epsilon$-stationary point of the instance we design; see Theorem~\ref{thm:lower}. Here  
		\begin{align}
			\label{eq:joint_condition_number}
		\Delta_{F_0}:= F_0(\vx^{(0)})-\inf_\vx F_0(\vx), \ [\bar{\vA};\vA]:= \left[
		\begin{array}{c}
			\bar{\vA}\\
			\vA
		\end{array}
		\right]\text{ and }	\kappa([\bar{\vA};\vA]):= \sqrt{\frac{\lambda_{\max}([\bar{\vA};\vA] [\bar{\vA};\vA]^\top )}{\lambda^+_{\min}([\bar{\vA};\vA] [\bar{\vA};\vA]^\top )}}.
		\end{align}
		Our lower complexity bound can be viewed as a generalization of the lower bound $O(\kappa(\vA)L_f \Delta_{f_0} \epsilon^{-2})$ in Sun and Hong~\cite{sun2019distributed} for affinely constrained smooth problem~\eqref{eq:affinenonconvex}. Our result provides a new insight that the difficulty of finding a (near) $\epsilon$-stationary point of problem~\eqref{eq:model} depends on the interaction between the affine constraints and the regularizer $g$. %characterized by $\kappa([\bar{\vA};\vA])$. 
		In contrast to Claim \ref{claim-1} for an affinely constrained \emph{smooth} non-convex optimization, we show that for \emph{non-smooth} composite non-convex optimization \eqref{eq:model}, the non-smooth regularization term $g$ can {make the problem significantly harder} %\blue{increase the lower complexity bound of FOMs} 
		even if $\prox_{\eta g}$ is accessible %\blue{at each iteration} 
		for any $\eta>0$; see Claim \ref{claim2} below.

Our second contribution is to show that %under Assumption~\ref{ass:linearspan3}, 
the minimum oracle complexity is  $O({\kappa([\bar{\vA};\vA]) L_f \Delta_{F_0}} \epsilon^{-2})$ 
for any FOM in Algorithm Class~\ref{algo:class-2} to find a near $\epsilon$-stationary point of problem instances in Problem Class~\ref{prob:class-2}. %that satisfies Assumption~\ref{assume:problemsetup}(a,c) and the condition $$\Delta_F :=F(\vx^{(0)},\vy^{(0)})-\inf_{\vx,\vy} F(\vx,\vy)<\infty.$$ 
The lower bound is established by using the same hard instance as that we use to prove the lower bound for FOMs in  Algorithm Class~\ref{algo:class-1}. %Moreover, under Assumption~\ref{ass:linearspan3}, the minimum oracle complexity 
%for an algorithm to find a near $\epsilon$-stationary point of problem~\eqref{eq:model} is $O({\kappa([\bar{\vA};\vA]) L_f \Delta_{F_0}} \epsilon^{-2})$, {the same as that we establish under Assumption~\ref{ass:linearspan} which allows the more expensive operator $\prox_{\eta g}$.}
In Section~4 of the extended arXiv version (see Liu et al.~\cite{liu2023first}) of this paper, we give an inexact proximal gradient (IPG) method that belongs to Algorithm Class~\ref{algo:class-2} for solving  Problem Class~\ref{prob:class-2}. It can find a near $\epsilon$-stationary point of any instance in Problem Class~\ref{prob:class-2} by calling $\mathrm{ORACLE}_2$ at most $\tilde O({\kappa([\bar{\vA};\vA]) L_f \Delta_{F_0}} \epsilon^{-2})$ times, which has a difference of up to a logarithmic factor from our lower bound. This shows that our lower oracle complexity bound is nearly tight. Similar results are shown in a concurrent paper (see Zhu et al.~\cite{zhu2023optimal}), which was released just a few weeks after the first version of our extended arXiv paper (see Liu et al.~\cite{liu2023first}). However, for non-convex problems, Zhu et al.~\cite{zhu2023optimal} only studied lower complexity bounds of FOMs in our Algorithm Class~\ref{algo:class-2}, and in addition, they do not require $\mathrm{relint}(\dom(\bar g)) \neq \emptyset$ and thus could directly use the result from Sun and Hong~\cite{sun2019distributed}. In contrast, our problem class is sharper due to the requirement of $\mathrm{relint}(\dom(\bar g)) \neq \emptyset$, and a new proof technique must be employed to handle a \emph{real} composite non-convex optimization with affine constraints.

%\comm{Can we claim that IPG is a contribution in this paper?}

\subsection{Related Work}\label{sec:relatedwork}
%The lower bound of iteration complexity is first introduced by Nemirovski and Yudin and tells algorithm designers which is the most efficient algorithm for a class of problems. The information-based bounds are established for a general class of problems and are significant in determining the potential for improvement of a numerical method for general use. %We then give some examples.

The \emph{proximal gradient} (PG) method can find an $\epsilon$-stationary point of the composite non-smooth non-convex problem~\eqref{eq:composite} 
within $O(L_f\epsilon^{-2})$ iterations~\cite{nesterov2012make}, which matches the lower bound in Carmon et al.~\cite{carmon2020lower} and Carmon et al.~\cite{carmon2021lower}. 
When there is no affine constraints and $g\equiv0$ in problem~\eqref{eq:model}, our lower-bound complexity result is reduced to the lower bound in Carmon et al.~\cite{carmon2020lower}.		
For the affinely constrained non-convex smooth problem~\eqref{eq:affinenonconvex}, Sun and Hong~\cite{sun2019distributed} show a lower-bound complexity of  $O(\kappa(\vA)L_f\epsilon^{-2})$ for finding an $\epsilon$-stationary point. 
In the same paper, they give an FOM that achieves this lower bound. 
When $g\equiv 0$, our complexity lower bounds are reduced to their lower bound. Lower complexity bounds of FOMs have also been established for \emph{smooth} non-convex strongly-concave min-max problems in Li et al.~\cite{li2021complexity}, and Zhang et al.~\cite{zhang2021complexity}.  
	
	Before our work, there only exist upper bounds of oracle complexity for finding a (near) $\epsilon$-stationary point of \eqref{eq:model} and~\eqref{eq:model-spli} with a non-smooth $g$ or $\bar g$. For instance, Kong et al.~\cite{kong2019complexity} develop a quadratic-penalty accelerated inexact proximal point method that finds an $\epsilon$-stationary point of problem~\eqref{eq:model} with oracle complexity $O(\epsilon^{-3})$. Lin et al.~\cite{lin2022complexity} study a method similar to Kong et al.~\cite{kong2019complexity} and show that oracle complexity of $O(\epsilon^{-5/2})$ is sufficient. The  ALM is another effective approach for problem~\eqref{eq:model}. The oracle complexity of ALM for problem~\eqref{eq:model} has been studied by Hong~\cite{hong2016decomposing}, Hajinezhad and Hong~\cite{hajinezhad2019perturbed}, Melo et al.~\cite{melo2020iteration1}, and Zhang and Luo~\cite{zhang2022global}. For example, the inexact proximal accelerated ALM by Melo et al.~\cite{melo2020iteration1} achieves oracle complexity of $O(\epsilon^{-5/2}$) and, in a special case where $g(\vx)$ is the indicator function of a polyhedron, the smoothed proximal ALM by Zhang and Luo~\cite{zhang2022global} improves the complexity to $O(\epsilon^{-2}$). ADMM is an effective algorithm for optimization with a separable structure  like that in problem~\eqref{eq:model-spli}. ADMM and its variants have been studied by Melo and Monterio~\cite{melo2017iteration}, Gonçalves et al.~\cite{goncalves2017convergence}, Jiang et al.~\cite{jiang2019structured}, Zhang and Luo~\cite{zhang2020proximal}, Yashtini~\cite{yashtini2022convergence}, and Hong et al.~\cite{hong2016convergence} for constrained non-convex  optimization problems including problem~\eqref{eq:model-spli}. For example, it is shown (see Gonçalves et al.~\cite{goncalves2017convergence}, Jiang et al.~\cite{jiang2019structured}, and Yashtini~\cite{yashtini2022convergence}) that ADMM finds an $\epsilon$-stationary point of  problem~\eqref{eq:model-spli} with oracle complexity of $O(\epsilon^{-2})$. 
	%However,~\eqref{zhang2022global} has to assume that $g(\vx)$ is the indicator function of a convex polyhedron

	The aforementioned methods for problem~\eqref{eq:model} all belong to Algorithm Class~\ref{algo:class-1}, and the aforementioned methods for problem~\eqref{eq:model-spli} all belong to Algorithm Class~\ref{algo:class-2}.  However, the oracle complexity of those methods for finding an  $\epsilon$-stationary point either does not match or is not comparable with our lower-bound complexity for the corresponding problems.  
	Specifically, the oracle complexity of ADMM in Yashtini~\cite{yashtini2022convergence} for solving problem~\eqref{eq:model-spli} depends on the Kurdyka-Łojasiewicz (KŁ)
	coefficient, which is not directly comparable with our lower bound. The oracle complexity $O(\kappa^{2}([\vA;\bar{\vA}])L_f^2\Delta_F \epsilon^{-2})$ of ADMM for solving problem~\eqref{eq:model-spli} is presented in Goncalves et al.~\cite{goncalves2017convergence}.
	The results in Jiang et al.~\cite{jiang2019structured} of ADMM are only applicable to problems~\eqref{eq:model} and~\eqref{eq:model-spli} with a separable structure and $\bar{\vA}=\vI_d$ or $g\equiv 0$, for which case their oracle complexity is $O(\kappa^{2}(\vA)L_f^2\Delta \epsilon^{-2})$ with $\Delta = \Delta_{F_0}$ or $\Delta = F(\vx^{(0)}, \vy^{(0)}) - \inf_{\vx,\vy}F(\vx, \vy)$.  Zhang and Luo~\cite{zhang2022global} study the complexity of ALM for problem~\eqref{eq:model} when $g$ is the indicator function of a polyhedral set and the exact projection onto the polyhedral set can be computed. They obtained complexity of $O(\widehat\kappa^{2}L_f^3\Delta_{F_0} \epsilon^{-2})$, where $\widehat\kappa$ is a joint condition number of the equality $\vA \vx + \textbf{b} =\vzero$  and the inequality defining the polyhedral set. When $g\equiv 0$, their complexity is reduced to $O(\kappa^{2}(\vA)L_f^3\Delta_{F_0} \epsilon^{-2})$. Zhang et al.~\cite{zhang2022iteration} further extend the results  to problems with nonlinear convex inequality constraints. However, their algorithm requires exact projection to the set defined by the inequality constraints, which is impractical for many applications and does not belong to Algorithm Class~\ref{algo:class-1}. 
	
	\subsection{Notations and Definitions}
	
%	\comm{Can we introduce the hoffman bound $\theta$ here?}
%	\comm{Can we introduce the complexity operators $\Theta$ and $O$ here?}
%	\YX{no need}

	For any $a\in\RR$, we use $\lceil a\rceil$ to denote the smallest integer that is no less than $a$. % and $\lceil a\rceil_+$ the smallest nonnegative integer that is no less than $a$. %$\cB_\delta(\vz)$ denotes a ball with radius $\delta$ and center $\vz$. %If $z=\vzero$, we simply use $\cB_\delta$.
	% We define $\cB_\delta^+$ as the intersection of $\cB_\delta$ with the nonnegative orthant, so in the $n$-dimensional space, $\cB_\delta^+=\cB_\delta\cap\RRR_+^n$. We use $V_m(\delta)$ for the volume of $\cB_\delta$ in the $m$-dimensional space. $[n]$ denotes the set $\{1,\ldots,n\}$. 
	$\mathbf{Null}(\vH)$ represents the null space of a matrix $\vH$. For any vector $\vz$, $[\vz]_j$ denotes its $j$-th coordinate. 	%$\mathbf{diag}(\vz)$ denotes a diagonal matrix with $\vz$ on the diagonal, and for any square matrix $\vA$, $\Diag(\vA)$ is a vector that takes the diagonal of $\vA$
	%Given two integers $a_1$ and $a_2$, we use  $a_1  \mid a_2$ ($a_1 \nmid a_2$) to represent that $a_2$ is (isn't) divisible by $a_1$. 
	We denote $\mathbf{1}_{p}$ for an all-one vector in $\mathbb{R}^p$, and $\mathbf{0}$ to represent an all-zero vector when its dimension is clear from the context. $\vI_p$ denotes a $p\times p$ identity matrix and  %$\vJ_p$ a $(p-1)\times p$ matrix with the following structure:
	\begin{equation}
		\label{eq:matrixJ}
		%			\vI_{\bar{d}}:=\left[\begin{array}{cccc}
			%			 1 &&&\\
			%		 &1 &&\\
			%			&& \ddots &  \\
			%			&&  & 1 \\
			%		\end{array}\right]
		%	\quad\text{and}	\quad
		\vJ_{p}:=\left[\begin{array}{ccccc}
			-1 & 1 &&&\\
			& -1 &1 &&\\
			&& \ddots & \ddots &  \\
			&&  &-1 & 1 \\
		\end{array}\right] \in \RR^{(p-1)\times p}.
	\end{equation}
A vector $\vx$ is said $\omega$-close to another vector $\widehat\vx$ if $\|\vx -\widehat\vx\| \le \omega$. For any set $\cX$, we denote $\iota_\cX$ as its indicator function, i.e., $\iota_\cX(\vx) = 0$ if $\vx\in\cX$ and $+\infty$ otherwise, and $\mathrm{relint}(\cX)$ denotes the relative interior of $\cX$. 	
	We use $\otimes$ for the Kronecker product, $\mathbf{co}(\mathcal{S})$ for the convex hull of a set $\mathcal{S}$ and $\overline{\mathbf{co}}(\cS)$ for the closure of $\mathbf{co}(\mathcal{S})$. %We use $O$, $\Theta$, and $o$ with standard meanings, while in the complexity result statement, $\widetilde O$ has a similar meaning as $O$ but suppresses a polynomial of $|\log\epsilon|$ for a given error tolerance $\epsilon>0$. 
	For a function $f:\mathbb{R}^d\rightarrow\mathbb{R}\cup \{+\infty\}$, its directional derivative at $\vx$ along a direction $\vv\in\mathbb{R}^d$  is defined as 
	\begin{equation}\label{eq:def-dir-der}
	f^{\prime}(\vx ; \vv)=\lim _{s \downarrow 0} \frac{f(\vx+s \vv)-f(\vx)}{s},
	\end{equation}
where $s\downarrow 0$ means $s\rightarrow 0$ and $s>0$.
$\partial g$ denotes the subdifferental of a closed convex function $g$.	%Note that the directional derivative exists for $f_0$, $g$ and $F$ under our Assumption~\ref{assume:problemsetup}. 
	
\begin{definition}\label{def:eps-pt-P}
Given $\epsilon \ge0$, a point $\vx^*$ is called an $\epsilon$-stationary point of~\eqref{eq:model} if for some $\vgamma\in\RR^n$,
\begin{equation}
		\label{eq:epsta}
		\max\left\{\dist\left(\mathbf{0},\nabla f_0(\vx^*)+\vA\zz \vgamma %+\nabla \textbf{h}(\vx^*)\zz \nu 
		+\partial g(\vx^*) \right), \left\| \vA\vx^*+\vb\right\| %,\left\| \nu\zz \textbf{h}(\vx^*)\right\|,\left\| (-\nu)_+ \right\|
		\right\}\leq \epsilon,
	\end{equation}
and a point $(\vx^*, \vy^*)$ is called an $\epsilon$-stationary point of~\eqref{eq:model-spli} if for some $\vz_1\in \RR^{\bar{n}}$ and $\vz_2\in\RR^n$,
\begin{equation}
	\label{eq:kktviofgesub}
	\begin{aligned}
		&	\max\left\{ \dist (\vzero, \partial \bar{g}(\vy^*) - \vz_1), \|\nabla f_0(\vx^*) + \bar{\vA}\zz \vz_1 + \vA\zz \vz_2\|,  \|\vy^*-\bar{\vA}\vx^*-\bar{\vb}\|,  \|\vA\vx^*+\vb\| \right\}
		 \leq\epsilon.
	\end{aligned}
\end{equation}	
When $\epsilon = 0$, we simply call $\vx^*$ and $(\vx^*, \vy^*)$ stationary points of problems~\eqref{eq:model} and~\eqref{eq:model-spli}, respectively.	We say that $\bar{\vx}$ is a near $\epsilon$-stationary point of~\eqref{eq:model} if $\|\bar\vx - \vx^*\| = O(\epsilon)$ where $\vx^*$ is an $\epsilon$-stationary point of~\eqref{eq:model}.
%	with $\omega= O(\epsilon)$.
\end{definition}	
	
	\subsection{Organization}
	
	The rest of this paper is organized as follows. In Section~\ref{sec:lb1}, we present lower-bound complexity results of FOMs in Algorithm Class~\ref{algo:class-1} for solving Problem Class~\ref{prob:class-2}. %problem~\eqref{eq:model}   under Assumption \ref{ass:linearspan}. 
	Then in Section~\ref{sec:extension}, we show the lower-bound complexity results of FOMs in Algorithm Class~\ref{algo:class-2} for solving Problem Class~\ref{prob:class-2}.
	Also we show the tightness of our lower complexity bounds for Algorithm Class~\ref{algo:class-2}. Concluding remarks and a few open questions are given in Section~\ref{sec:conclusion}.
	
	\section{Lower Bound of Oracle Complexity for Algorithm Class~\ref{algo:class-1}}\label{sec:lb1}
	
	In this section, for a given $\epsilon>0$, we derive a lower bound of the oracle complexity of an FOM in Algorithm Class~\ref{algo:class-1} to find a (near) $\epsilon$-stationary point of problem instances in Problem Class~\ref{prob:class-2}.   
	
	\subsection{A Challenging Instance  $\cP$ in Problem Class~\ref{prob:class-2}}
Our hard instance is built based on the instance in Sun and Hong~\cite{sun2019distributed} by carefully separating its affine constraints. Some of its linear constraints are kept as hard constraints in our instance, and the others are used to design a non-smooth convex term in our objective function.  The smooth term in the objective of our instance is almost the same as that in the instance of Sun and Hong~\cite{sun2019distributed}, and a slight difference is in some coefficients that we choose to satisfy the smoothness condition with a given (but arbitrary) parameter $L_f>0$.		
 More precisely, we let $m_1$ and $m_2$ be positive integers such that $m_1m_2$ is even and $m_1\geq2$. Let $m=3m_1m_2$ and $\bar{d}$ be an odd positive integer such that $\bar{d}\geq 5$. Here, $m$ will be proportional to $\kappa([\bar{\vA}; \vA])$ for our instance; see lemma \ref{lem:condH}. Also, set $d=m{\bar{d}}$, and  let  
\begin{eqnarray}
	\label{eq:xblock}
	\vx=\left(\vx_1\zz, \ldots,\vx_m\zz\right)\zz\in\mathbb{R}^{d}, \text{ with }\vx_i\in\mathbb{R}^{\bar d}, \, i=1,\dots,m.
\end{eqnarray}
	%Let $\vI_{\bar{d}}\in\mathbb{R}^{\bar{d}\times \bar{d}}$ be an identity matrix and  $\vJ_{m}\in\mathbb{R}^{(m-1)\times m}$ defined as
Moreover, we define a matrix  $\vH\in\mathbb{R}^{(m-1)\bar{d}\times m\bar{d}}$ by
\begin{equation}
	\label{eq:matrixAstar}
	\newcommand{\zm}{
			\left[\begin{array}{ccccc}
			-\vI_{\bar{d}} & \vI_{\bar{d}} &&&\\
			& -\vI_{\bar{d}} &\vI_{\bar{d}} &&\\
			&& \ddots & \ddots &  \\
			&&  &-\vI_{\bar{d}} & \vI_{\bar{d}} \\
		\end{array}\right]
	}		
		\vH:=mL_f\cdot \vJ_{m}\otimes \vI_{\bar{d}}=
		mL_f\cdot
		\left.
		\,\smash[b]{\underbrace{\!\zm\!}_{\textstyle\text{$m$ blocks}}}\,
		\right\}\text{$m-1$ blocks.}
		\vphantom{\underbrace{\zm}_{\textstyle\text{$m$ blocks}}}		
\end{equation}
%
%with $\vx_i\in\mathbb{R}^{\bar d}$ for $i=1,\dots,m$ .  
Define 
\begin{equation}
\label{eq:indexsetM}
\begin{aligned}
\mathcal{M}:=\{im_1| i=1,2,\dots,3m_2-1\},& \quad \mathcal{M}^C:=\{1,2,\dots,m-1\}\backslash  \mathcal{M}, \\ n=(m-3m_2){\bar{d}},& \quad \bar n=(3m_2-1){\bar{d}},
\end{aligned}
\end{equation}
and let  
\begin{eqnarray}
	\label{eq:AandAbar}
\bar{\vA}:=mL_f\cdot \vJ_{\mathcal{M}}\otimes \vI_{\bar{d}},\quad 
\vA:=mL_f\cdot \vJ_{\mathcal{M}^C}\otimes \vI_{\bar{d}}, \quad \bar\vb = \vzero \in \RR^{\bar n}, \quad \vb = \vzero \in \RR^n,
\end{eqnarray}
where $\vJ_{\mathcal{M}}$ and $\vJ_{\mathcal{M}^C}$ are the rows of $\vJ_{m}$ indexed by $\mathcal{M}$ and $\mathcal{M}^C$, respectively.

{We partition the matrix $\vH$ into two submatrices $\bar\vA$ and $\vA$ as described in \eqref{eq:AandAbar} for two primary goals. First, as we will argue later above Claim~\ref{claim2}, %we aims to maintain relatively low condition numbers for matrices $\bar{\vA}, \vA$, and a high condition number for $[\bar{\vA}; \vA]$. This ensures that 
it holds $\frac{\kappa([\bar{\vA};\vA])}{\kappa(\vA)} \ge \frac{3 m_2}4$. This fact will be used to show that the existence of a nonsmooth $g$ can make the problem significantly harder. 
%can be large enough. 
Second, it is crucial that with the constructed $\bar\vA$ and $\vA$, we can still obtain zero-respecting sequences, illustrated in Figure~\ref{fig:1} below. %can still be effectively estimated under our partitioning scheme (see Figure \ref{fig:1}).
}

%Let $\vy=(\vy_1\zz, \ldots,\vy_{3m_2-1}\zz)\zz\in\mathbb{R}^{(3m_2-1)\bar{d}}$ with $\vy_i\in\mathbb{R}^{\bar d}$ for $i=1,\dots,3m_2-1$. 
Furthermore, given $\epsilon \in(0,1)$, we define $\bar g:\mathbb{R}^{\bar{n}}\rightarrow \mathbb{R}$ and $g:\RR^d \rightarrow \RR$ as
\begin{equation}
	\label{eq:gbar}
\bar g(\vy):=\frac{\beta}{mL_f}\|\vy\|_1=\max\left\{\vu^\top\vy~\Big|~\|\vu\|_\infty\leq \frac{\beta}{mL_f}\right\},
\end{equation}
and \begin{equation}
	\label{eq:g}
	g(\vx):=\bar{g}(\bar\vA\vx)=\beta\sum_{i\in\mathcal{M}}  \|\vx_{i }-\vx_{i+1}\|_1,
\end{equation}
where $\vx$ has the block structure in~\eqref{eq:xblock}, $\mathcal{M}$ is defined in \eqref{eq:indexsetM}, and $\beta$ is a constant satisfying %\YX{Do we need $\beta\ge 1$?}	 \revise{we don't}
\begin{equation}
	\label{eq:betachoice}
	\beta>(50\pi+1+\|\vA\|)\sqrt{m}\epsilon.
\end{equation} 
%Using $\bar\vA$ in~\eqref{eq:AandAbar}, we set $g$ in problem~\eqref{eq:model} to be

%Since $m_1\geq 2$, each summand in~\eqref{eq:g} depends on different variables $\vx_{i }$'s. It is easy to see that $\vx_1=\vx_2=\cdots=\vx_m$ if and only if $\vA\vx+\vb=\mathbf{0}$ and $g(\vx)=0$ according to our choices above. 

Finally, we design the smooth function $f_0$ to complete the instance of~\eqref{eq:model} in Problem Class~\ref{prob:class-2}. Let $\Psi:\RR\mapsto\RR$ and $\Phi:\RR\mapsto\RR$ be defined as
	\begin{equation}
		\label{eq:PsiPhi}
	\Psi(u):=\left\{\begin{array}{ll}
		0, & \text{ if }u \leq 0, \\
		1-e^{-u^2}, & \text{ if } u>0,
	\end{array}\right.  \quad \text{and}\quad\Phi(v):=4 \arctan v+2 \pi.
	\end{equation}
In addition, for each $j=1,\ldots,\bar d$, we define function $\varphi\left(\cdot, j\right):\RR^{{\bar{d}}} \rightarrow \RR$ as
   \begin{equation}
   	\label{eq:Theta}
   	\varphi\left(\vz, j\right):=\left\{
   	\begin{array}{ll}
	-\Psi(1) \Phi\left([\vz]_1\right), &\text{ if } j=1,  \\[0.1cm]
		\Psi\left(-[\vz]_{j-1}\right) \Phi\left(-[\vz]_j\right)-\Psi\left([\vz]_{j-1}\right) \Phi\left([\vz]_j\right), &\text{ if } j = 2, \dots,{\bar{d}},
		\end{array}\right.
	\end{equation}	
and	for $i=1,\dots,m$, we define $h_i:\RR^{{\bar{d}}} \rightarrow \RR$ as
	\begin{equation}
		\label{eq:hi}
h_i\left(\vz\right):= \left\{\begin{array}{ll}
	\varphi\left(\vz, 1\right)+3 \sum_{j=1}^{\lfloor {\bar{d}} / 2\rfloor} \varphi\left(\vz, 2 j\right), & \text{ if }i \in\left[1, \frac{m}{3}\right], \\ [0.1cm]
\varphi\left(\vz, 1\right), &\text{ if } i \in\left[\frac{m}{3}+1, \frac{2 m}{3}\right], \\ [0.1cm]
\varphi\left(\vz, 1\right)+3 \sum_{j=1}^{\lfloor {\bar{d}} / 2\rfloor} \varphi\left(\vz, 2 j+1\right), &\text{ if } i \in\left[\frac{2 m}{3}+1, m\right].
\end{array}\right.
	\end{equation}	
Now for the given $\epsilon\in(0,1)$ and $L_f > 0$,  we define $f_i:\RR^{{\bar{d}}} \rightarrow \RR$ as
\begin{equation}
	\label{eq:fi}
	f_i\left(\vz\right):=\frac{300  \pi \epsilon^2}{mL_f} h_i\left(\frac{\sqrt{m}L_f \vz }{150\pi\epsilon}\right), \forall\, \vz\in \RR^{\bar d}, \forall\, i=1,\dots,m,
\end{equation}	
%Let $d=m{\bar{d}}$ and $\vx=(\vx_1\zz, \ldots,\vx_m\zz)\zz\in\mathbb{R}^{d}$ with $\vx_i\in\mathbb{R}^{\bar d}$ for $i=1,\dots,m$. 
and let $f_0: \RR^d \rightarrow \RR$ be
\begin{equation}
	\label{eq:f0}
f_0(\vx):= \sum^{m}_{i=1}  f_i(\vx_i), \forall\, \vx\in \RR^d \text{ with the structure in }\eqref{eq:xblock}.
\end{equation}

Putting all the components given above, we obtain a specific instance of~\eqref{eq:model}. We formalize it in the following definition.
%\revise{$(\bar{\vA}, \bar{\vb})$ are neglected in the following definition.}
\begin{definition}[instance $\cP$]
	\label{def:hardinstance}
Given $\epsilon\in(0,1)$ and $L_f > 0$, let $m_1, m_2$ and $\bar d$ be integers such that $m_1\ge2$ is even and $\bar d\ge 5$ is odd. We refer to as \emph{instance $\cP$} the instance of problem~\eqref{eq:model} where $f_0$ is given in~\eqref{eq:f0} with each $f_i$ defined in~\eqref{eq:fi}, $g$ is given in~\eqref{eq:g} with $\beta$ satisfying~\eqref{eq:betachoice}, and $(\vA, \bar\vA, \vb, \bar\vb)$ is given in~\eqref{eq:AandAbar}. %called instance $\mathcal{P}(\bar{d}, m_1, m_2, \epsilon, L_f)$, or simply instance~$\mathcal{P}$, if it satisfies $d=m{\bar{d}}$, 
%$n=(m-3m_2){\bar{d}}$, $\bar{n}=(3m_2-1){\bar{d}}$, $m=3m_1m_2$, $\vb=\mathbf{0}$, $\bar{\vb}=\mathbf{0}$, definitions~\eqref{eq:xblock},~\eqref{eq:AandAbar},
%\eqref{eq:f0},~\eqref{eq:g}, and~\eqref{eq:betachoice}, where $m_1$ and $m_2$ be positive integers such that $m_1m_2$ is even and $m_1\geq2$, and $\bar{d}$ is an odd positive integer such that $\bar{d}\geq 5$. 
\end{definition}

\subsection{Properties of Instance $\cP$}

In order to show the challenge of instance $\cP$ for an FOM in Algorithm Class~\ref{algo:class-1},  we give a few facts and properties about $\cP$. First, notice that $\bar{\vA}$ and $\vA$ in~\eqref{eq:AandAbar} are two block submatrices of $\vH$ in rows. 
 It is easy to obtain the following proposition.
 \begin{proposition}\label{prop:consensus}
 By the definitions in~\eqref{eq:xblock} through~\eqref{eq:g}, it holds
 \begin{enumerate}
 \item[\textnormal{(a)}] $\vx_1=\vx_2=\cdots=\vx_m$ if and only if $\vH\vx=\mathbf{0}$; 
 \item[\textnormal{(b)}] $\vx_i=\vx_{i+1}$ for $i\in \mathcal{M}$ if and only if $\bar\vA\vx=\mathbf{0}$ or equivalently $g(\vx)=0$; 
 \item[\textnormal{(c)}] $\vx_i=\vx_{i+1}$ for $i\in \mathcal{M}^C$ if and only if $\vA\vx=\mathbf{0}$; 
 \item[\textnormal{(d)}] $\vx_1=\vx_2=\cdots=\vx_m$ if and only if $\vA\vx=\mathbf{0}$ and $g(\vx)=0$.
 \end{enumerate} 
 \end{proposition}	
		
Second, it is straightforward to have	$\left\{\frac{\partial h_i(\vz)}{\partial [\vz]_j}\right\}$ based on three cases of $j$, also given in Sun and Hong~\cite{sun2019distributed}. %for $i=1,\dots,m$ and $j=1,\dots,\bar{d}$. They are derived under three cases of $j$.
	
	\noindent Case i) When $j=1$,
	\begin{equation}
		\label{eq:deriveh3}
		\frac{\partial h_i(\vz)}{\partial [\vz]_j} =\left\{\begin{array}{ll}
			-\Psi(1) \Phi^{\prime}([\vz]_{1})+3\left[-\Psi^{\prime}(-[\vz]_{1}) \Phi(-[\vz]_{2})-\Psi^{\prime}([\vz]_{1}) \Phi([\vz]_{2})\right],& \text{ if }i \in\left[1, \frac{m}{3}\right], \\[0.1cm]
			-\Psi(1) \Phi^{\prime}([\vz]_{1}), &\text{ if }i \in\left[\frac{m}{3}+1, m\right].
		\end{array}\right.
	\end{equation}

\noindent	Case ii) When $j$ is even, 
	\begin{equation}
		\label{eq:deriveh1}
		\frac{\partial h_i(\vz)}{\partial [\vz]_j} =\left\{\begin{array}{ll}
			3\left[-\Psi(-[\vz]_{j-1}) \Phi^{\prime}(-[\vz]_j)-\Psi([\vz]_{j-1}) \Phi^{\prime}([\vz]_j)\right], & \text{ if }i \in\left[1, \frac{m}{3}\right], \\[0.1cm]
			0,			& \text{ if }i \in\left[\frac{m}{3}+1, \frac{2m}{3}\right], \\[0.1cm]
			3\left[-\Psi^{\prime}(-[\vz]_j) \Phi(-[\vz]_{j+1})-\Psi^{\prime}([\vz]_j) \Phi([\vz]_{j+1})\right], & \text{ if }i \in\left[\frac{2m}{3}+1, m\right].
		\end{array}\right.
	\end{equation}
	
\noindent	Case iii) When $j$ is odd and $j\neq 1$, 
	\begin{equation}
		\label{eq:deriveh2}
	\frac{\partial h_i(\vz)}{\partial [\vz]_j} = \left\{\begin{array}{ll}
		3\left[-\Psi^{\prime}(-[\vz]_j) \Phi(-[\vz]_{j+1})-\Psi^{\prime}([\vz]_j) \Phi([\vz]_{j+1})\right], & \text{ if }i \in\left[1, \frac{m}{3}\right], \\[0.1cm] 
		0,	& \text{ if }i \in\left[\frac{m}{3}+1, \frac{2m}{3}\right], \\[0.1cm] 
		3\left[-\Psi(-[\vz]_{j-1}) \Phi^{\prime}(-[\vz]_j)-\Psi([\vz]_{j-1}) \Phi^{\prime}([\vz]_j)\right], & \text{ if }i \in\left[\frac{2m}{3}+1, m\right].
		\end{array}\right.
	\end{equation}

%Next we construct an instance of~\eqref{eq:model}.  
%with $f_i$ defined in~\eqref{eq:fi} and choose $g$ to be

%Then we choose $\vb$ and $\vA$ in~\eqref{eq:model} as 

%Since $\vb\in\mathbb{R}^{(m-3m_2){\bar{d}}}$ and $\vA\in\RR^{(m-3m_2){\bar{d}}\times m\bar{d}}$, we have $d=m{\bar{d}}$ and $n=(m-3m_2){\bar{d}}$ in this instance of~\eqref{eq:model}. 
%We then summarize the instance of problem~\eqref{eq:model} constructed above as follows.

%We will show next that 
Thirdly, instance~$\mathcal{P}$ belongs to Problem Class~\ref{prob:class-2}. By~\eqref{eq:g}, it is clear that $\bar g$ is convex and $\mathrm{relint}(\dom(\bar g)) \neq \emptyset$. %Assumption~\ref{assume:problemsetup}(c) clearly holds. %with $\bar g$ in~\eqref{eq:gbar}, $\bar{\vA}$ in~\eqref{eq:AandAbar} and $\bar{\vb}=\mathbf{0}$.  
To verify the other requirements, %Assumptions~\ref{assume:problemsetup}(a) and~\ref{assume:problemsetup}(b), 
we need the next lemma.
\begin{lemma}
	\label{lem:functions}
Let $\Psi$, $\Phi$ and $h_i$ be given in~\eqref{eq:PsiPhi} and~\eqref{eq:hi}.	The following statements hold:
\begin{enumerate}	
	\item[\textnormal{(a)}]  %$\Psi,  \Psi^{\prime}, \Phi$ and $\Phi^{\prime}$ are nonnegative and bounded, and  $\Psi^{\prime \prime}$ and $\Phi^{\prime \prime}$ are bounded. Specifically, it holds that 
	$\Psi(u)=0$ and $\Psi'(u)=0$ for all $u\leq 0$, and $0 \leq \Psi(u)<1$, $0 \leq \Psi^{\prime}(u) \leq \sqrt{2/e}$,
	$0<\Phi(v)<4 \pi$, and $0<\Phi^{\prime}(v) \leq 4$ for all $v\in\RR$.
\item	[\textnormal{(b)}] 
	$
	\Psi(u) \Phi^{\prime}(v)>1
	$
	for all $u$ and  $v$ satisfying $u \geq 1$ and $|v|<1$.
	
\item[\textnormal{(c)}] $h_i(\mathbf{0})-\inf_{\vz} h_i\left(\vz\right)\leq 10\pi {\bar{d}}$ for $i=1,2,\ldots,m$.
	
\item[\textnormal{(d)}] $\nabla {h}_i$ is $75 \pi$-Lipschitz continuous  for $i=1,2,\ldots,m$.
	
\item[\textnormal{(e)}]  ${h}_i$ is $25\pi\sqrt{{\bar{d}}}$-Lipschitz continuous for $i=1,2,\ldots,m$.
\end{enumerate}
\end{lemma}
\proof{Proof of Lemma \ref{lem:functions}.}
	(a)-(d) are directly from Sun and Hong~\cite[Lemma 3.1]{sun2019distributed}. By \eqref{eq:deriveh3}--\eqref{eq:deriveh2}, we obtain that, for any $\vz\in\RR^{\bar{d}}$, 
	\begin{equation}
		\begin{aligned}
			\label{eq:partialhmax}
			\left|\frac{\partial h_i(\vz)}{\partial [\vz]_j}\right| \leq \max \left\{\sup _u\left|\Psi(1) \Phi^{\prime}(u)\right|+6 \sup _u\left|\Psi^{\prime}(u)\right| \sup _v|\Phi(v)|, 6 \sup _u\left|\Psi(u)\right| \sup _v\left|\Phi^{\prime}(v)\right|\right\}< 25\pi, \forall\, i, j,
		\end{aligned}
	\end{equation}
	where the second inequality is from Lemma~\ref{lem:functions}(a) and the fact that $\sup _v\left|\Psi(1) \Phi^{\prime}(v)\right|\leq 4(1-e^{-1})<\pi$ by the definition of $\Phi$. Hence, $\|\nabla h_i(\vz) \|\leq 25\pi\sqrt{{\bar{d}}}$. Thus (e) holds, and we complete the proof.
	\null\hfill$\square$
\endproof
From the definition of $f_i$ in \eqref{eq:fi} and the chain rule, we have 
\begin{equation}\label{eq:grad-f-i-chain}
 	[\nabla f_i(\vz)]_j  =  \frac{2\epsilon}{\sqrt{m}} \left[\nabla {h}_i\left(\frac{\sqrt{m} L_f \vz}{150\pi\epsilon}\right)\right]_j,\quad \forall j=1,\dots,m,
 \end{equation}
 which together with Lemma~\ref{lem:functions} clearly indicates the properties below about $\{f_i\}_{i=0}^m$.  
	
	\begin{lemma}
	\label{cor:boundf}
Let $\{f_i\}_{i=0}^m$	be defined in~\eqref{eq:fi} and~\eqref{eq:f0}. Then	%The following statements hold:
\begin{enumerate}	
	\item[\textnormal{(a)}] $ f_i(\mathbf{0})-\inf _{\vz} f_i\left(\vz\right)\leq 3000\pi^2 {\bar{d}}\epsilon^2/mL_f$ for $i=1,\dots,m$, and $ f_0(\mathbf{0})-\inf _{\vx} f_0\left(\vx\right)\leq 3000\pi^2 {\bar{d}}\epsilon^2/L_f$.
	
	\item[\textnormal{(b)}] $\nabla f_i$ is $L_f$-Lipschitz continuous for $i=0,1,\dots,m$.
	
	\item[\textnormal{(c)}] $f_i$ is $\frac{50\pi\epsilon\sqrt{{\bar{d}}}}{\sqrt{m}}$-Lipschitz continuous for $i=1,\dots,m$, and $f_0$ is ${50\pi\epsilon\sqrt{m{\bar{d}}}}$-Lipschitz continuous.
\end{enumerate}	
	\end{lemma}
%\begin{proof}
%	(a), (b) and (c) is  from  Lemma~\ref{lem:functions} (c), (d) and (e), respectively, given with the definitions of $f_i$ and $f_0$. 	
%	%(b) is  from  Lemma~\ref{lem:functions}(d) and the definitions of $f_i$ and $f_0$.
%\end{proof} 

%The  following proposition follows  Lemma~\ref{cor:boundf} (a) and (c) and the earlier discussions.
By Lemma~\ref{cor:boundf} and $\bar g \ge0$, we immediately have the   following proposition. %show that instance $\cP$ satisfies Assumptions~\ref{assume:problemsetup}(a) and~\ref{assume:problemsetup}(b).
\begin{proposition}
Instance~$\mathcal{P}$ given in Definition~\ref{def:hardinstance} belongs to Problem Class~\ref{prob:class-2}.
\end{proposition}

The following lemma characterizes the joint condition number of $\bar{\vA}$ and $\vA$ defined in~\eqref{eq:joint_condition_number}, which will appear in our lower bound of oracle complexity.  

\begin{lemma}
\label{lem:condH}
Let $\vA$ and $\bar\vA$ be given in~\eqref{eq:AandAbar}. Then $\frac{m}{4} \le \kappa([\bar{\vA}; \vA])=\kappa(\vH)={\frac{\sin(\frac{(3m_1m_2-1)\pi}{6m_1m_2})}{\sin(\frac{\pi}{6m_1m_2})}}<m$.
\end{lemma}
\proof{Proof of Lemma \ref{lem:condH}.}
%{\color{red} Qihang: please give a formal proof to get rid of ``$\approx$''.}
The first equality holds because $\vH$ are split into $\bar{\vA}$ and $\vA$ in rows. %and  $\| [\bar{\vA}; \vA]\|= \|\vH\|$. 
	Let 
	\begin{equation*}
		\begin{aligned}
			\bar{\vH}=mL_f\left[\begin{array}{cccc}
				-1 & 1 &\\
				& \ddots & \ddots &  \\
				&  &-1 & 1 \\
			\end{array}\right]\in\RR^{(m-1)\times m}.
		\end{aligned}
	\end{equation*}
	We then have $\vH \vH\zz=(\bar{\vH}\otimes \vI_{\bar{d}})(\bar{\vH}\otimes \vI_{\bar{d}})\zz=\left(\bar{\vH}\bar{\vH}\zz\right) \otimes\left( \vI_{\bar{d}} \vI_{\bar{d}}^{\top}\right)=\left(\bar{\vH}\bar{\vH}\zz\right) \otimes\vI_{\bar{d}} $. Let $\lambda_i(\bar{\vH}\bar{\vH}\zz)$ be the $i$-th largest eigenvalue of $\bar{\vH}\bar{\vH}\zz$. 
	Since $\bar{\vH}\bar{\vH}\zz$ is tridiagonal and Toeplitz, its eigenvalues have closed forms (see Gray~\cite{gray2006toeplitz}):  
	$$\textstyle \lambda_i(\bar{\vH}\bar{\vH}\zz)= 4m^2L_f^2\sin^2\left(\frac{i\pi}{6m_1m_2}\right), \forall\, i=1,2,\ldots, m-1.$$ It then yields that $\kappa([\bar{\vA}; \vA])=\kappa(\vH )=\kappa(\bar{\vH})= {\frac{\sin(\frac{(3m_1m_2-1)\pi}{6m_1m_2})}{\sin(\frac{\pi}{6m_1m_2})}}$. Because  $\sin(z) \le 1, \forall\, z$, $\sin(z)\geq \frac{2z}{3}$ for $z\in[0,\pi/12]$, and $m_1m_2\geq 2$, we have 
		%\begin{equation*}
$	\frac{\sin(\frac{(3m_1m_2-1)\pi}{6m_1m_2})}{\sin(\frac{\pi}{6m_1m_2})}\leq\frac{1}{\frac{\pi}{9m_1m_2}}=\frac{3m}{\pi}<m.
$	%\end{equation*}
Also,	because  $z\geq\sin(z)\geq \frac{z}{2}$ for $z\in[0,\pi/2]$ and $m_1m_2\geq 2$, we have 
	%\begin{equation*}
${\frac{\sin(\frac{(3m_1m_2-1)\pi}{6m_1m_2})}{\sin(\frac{\pi}{6m_1m_2})}}\geq	\frac{\frac{(m-1)\pi}{4m}}{\frac{\pi}{2m}}=\frac{m-1}{2}\geq \frac{m}{4}.
$	%\end{equation*}
 Hence, we have obtained all desired results and complete the proof.	
	\null\hfill$\square$
\endproof

		\subsection{An Auxiliary Problem and Its Properties}

	To establish a lower bound of the oracle complexity of FOMs in Algorithm Class~\ref{algo:class-1} for solving Problem Class~\ref{prob:class-2}, %under Assumptions~\ref{assume:problemsetup} and \ref{ass:linearspan}, 
	we consider an auxiliary problem of instance~$\mathcal{P}$ in this subsection and analyze its properties. The Auxiliary Problem is given as follows:
	\begin{equation}
		\label{eq:model3}
		\begin{aligned}
			\min_{ \vx\in\RR^d} \,\,&   f_0(\vx), 
			\,\,\st  \vH \vx =  \mathbf{0},
		\end{aligned}
		\tag{AP}
	\end{equation}
	where $\vH$ and $f_0$  are defined in~\eqref{eq:f0} and~\eqref{eq:matrixAstar}. An $\epsilon$-stationary point $\vx^*$ of problem~\eqref{eq:model3} satisfies
	\begin{equation}
		\label{eq:KKTviomodel2}
		\max\left\{\left\|\vH\vx^*\right\|, \min _{\vgamma'\in\RR^{(m-1){\bar{d}}}}\left\|\nabla  f_0(\vx^*)+\vH^{\top} \vgamma'\right\|\right\}\leq \epsilon.
	\end{equation}
	The next lemma characterizes the relationship between the (near-)stationary points of  instance $\cP$ and the auxiliary problem~\eqref{eq:model3}. %The proof techniques  have been utilized in~\cite{liu2022linearly,liu2022inexact}.
		\begin{lemma}\label{cor:kktvio2}
	%\label{lem:kktequiv}
Let $\epsilon >0$ be given in Definition~\ref{def:hardinstance} for instance $\cP$. Then for any $\widehat\epsilon \in [0, \epsilon]$, an $\widehat\epsilon$-stationary point of instance~$\mathcal{P}$ is also an $\widehat\epsilon$-stationary point of  the auxiliary problem~\eqref{eq:model3}. 
\end{lemma}

\proof{Proof of Lemma \ref{cor:kktvio2}.}
			Suppose $\vx^*$ is an $\widehat\epsilon$-stationary point of instance~$\mathcal{P}$, i.e., %According to~\eqref{eq:epsta}, 
			 for some $\vgamma\in\RR^n$ and $\vxi\in\partial g(\vx^*)$ such that 
			\begin{equation}
				\label{eq:KKTe}
				\|\nabla f_0(\vx^*)+ \vxi+\vA\zz \vgamma\|\leq  \widehat\epsilon \quad\text{ and }\quad\|\vA\vx^*\|\leq \widehat\epsilon.
			\end{equation}
			By the definitions of $g$ in \eqref{eq:gbar} and~\eqref{eq:g}, there exists $\vu\in\mathbb{R}^{\bar{n}}$ such that $\vxi=\bar\vA^{\top}\vu$ and thus the first condition in~\eqref{eq:KKTe} becomes
			\begin{equation}
				\label{eq:KKTHe}
				\|\nabla f_0(\vx^*)+ \bar\vA^{\top}\vu+\vA\zz \vgamma \|\leq  \widehat\epsilon.
			\end{equation}

Hence, in order to show that $\vx^*$ is an $\widehat\epsilon$-stationary point of  problem~\eqref{eq:model3}, we only need to prove~$\|\vH\vx^*\|\leq~\widehat\epsilon$. To show this, we first notice from the definition in~\eqref{eq:def-dir-der} that, for any $\vv\in\mathbf{Null}(\vA)$, 
\begin{equation}
	\label{eq:dderivegeq0epsilon}
	\begin{aligned}
		F_0^{\prime}(\vx^*;\vv)= \vv\zz\nabla  f_0(\vx^*)+ g^{\prime}(\vx^*;\vv) \geq \vv\zz\nabla  f_0(\vx^*) +\vv\zz \vxi =\vv\zz\left(\nabla f_0(\vx^*)+ \vxi+\vA\zz \vgamma\right)\geq -\widehat\epsilon\|\vv\|,
	\end{aligned}
\end{equation}
where the first inequality follows from Clarke~\cite{clarke1990optimization} 
\begin{equation}\label{eq:result-dir-der}
g^{\prime}(\vx;\vv)=\sup_{\vxi'\in\partial g(\vx)}\vv\zz \vxi', \forall\, \vx\in \dom(g), \forall\, \vv,
\end{equation} 
and the second inequality is by~\eqref{eq:KKTe} and the Cauchy-Schwarz inequality.

%Next, we prove that $\|\vH\vx^*\|\leq\epsilon$. Since $\|\vA\vx^*\|\leq \epsilon$, it suffices to 

Second, we claim $\bar\vA\vx^*=\mathbf{0}$, namely,  $\vx^*_{i}=\vx^*_{i+1}$ for all $i\in\mathcal{M}$, where $\mathcal{M}$ is defined in~\eqref{eq:indexsetM}. Suppose this claim is not true. Then for some $\bar{i}\in \mathcal{M}$, it holds $\vx^*_{\bar i}\neq\vx^*_{\bar i+1}$. 
We let %There exists $\bar{j}\in\{1,2,\ldots,\bar{d}\}$ such that  
\begin{equation}\label{eq:set-delta}
\bar{\vxi}:=\vx^*_{\bar{i}+1}-\vx^*_{\bar{i}}\neq \mathbf{0},
%\bar{\vxi}:=\left[\vx^*_{\bar{i}m_1+1}\right]_{\bar{j}}-\left[\vx^*_{\bar{i} m_1}\right]_{\bar{j}}\neq \mathbf{0}_{\bar d}.
\end{equation}
and $\vv^*=(\vv_1^{*\top},\vv_2^{*\top},\dots,\vv_m^{*\top})^\top$ where each $\vv_i^*\in\mathbb{R}^{\bar{d}}$ is defined as
$$
\vv_i^*= \bar{\vxi}, \text{ if }i\leq \bar{i}, \text{ and } \vv_i^*= \mathbf{0},	\text{ if }i>\bar{i}. 
$$
It is easy to see that $\vv^*_i=\vv^*_{i+1}$ for any $i\neq \bar{i}$. Thus by {Proposition~\ref{prop:consensus}(c)}, $\vA\vv^*=\mathbf{0}$ and, for any $s\in(0,1)$,
%\revise{$s$ may be conflict with that in Proposition 2.3.} \YX{This is Okay.}
\begin{align}\label{eq:gtvchange}
	\nonumber
	 g\left(\vx^*+s \vv^*\right) %\\	\nonumber
	= & \beta  \sum_{i<\bar{i}}\left\|\vx^*_i+s\bar{\vxi}-\vx^*_{i+1}-s\bar{\vxi}\right\|_1+\beta \|\vx^*_{\bar{i}}+s\bar{\vxi}-\vx^*_{\bar{i}+1}\|_1
	+\beta \sum_{i>\bar{i}}\left\|\vx^*_i-\vx^*_{i+1}\right\|_1 \\	%\nonumber
	= & \beta  \sum_{i<\bar{i}}\left\|\vx^*_i-\vx^*_{i+1}\right\|_1+\beta (1-s)\|\vx^*_{\bar{i}}-\vx^*_{\bar{i}+1}\|_1
	+\beta \sum_{i>\bar{i}}\left\|\vx^*_i-\vx^*_{i+1}\right\|_1 
	= g\left(\vx^*\right)-s\beta\|\bar{\vxi}\|_1, 
\end{align}
where the first and last equalities follow from \eqref{eq:g} and the definition of $\vv^*$, and the second one is by~\eqref{eq:set-delta} and $s\in(0,1)$.
In addition, from \eqref{eq:fi} and \eqref{eq:partialhmax}, we have that, for any $\vz\in\mathbb{R}^{\bar d}$,  
\begin{equation}
	\label{eq:fi_infinitynorm}
\|\nabla f_i(\vz)\|_\infty=\frac{2 \epsilon}{\sqrt{m}} \left\|\nabla h_i\left(\frac{\sqrt{m}L_f \vz }{150\pi\epsilon}\right)\right\|_\infty
%\leq \frac{2 \epsilon}{\sqrt{m}}\max_{\vz}\max_{j}\left|\frac{\partial h_i(\vz)}{\partial z_j}\right|
\leq \frac{50\pi\epsilon}{\sqrt{m}}.
\end{equation}
Moreover, by the definition of $\vv^*_i$ for $i=1,\dots,m$, we have 
$$
f_i(\vx^*_i+s \vv^*_i)- f_i(\vx^*_i)=s\nabla f_i(\vx^*_i+s' \vv^*_i)^\top\vv^*_i\leq s\|\nabla f_i(\vx^*_i+s' \vv^*_i)\|_\infty\|\vv^*_i\|_1 \overset{\eqref{eq:fi_infinitynorm}}\leq \frac{50s\pi\epsilon}{\sqrt{m}}\|\bar{\vxi}\|_1,
$$
where the equality holds from the mean value theorem for some $s'\in(0,s)$. The inequality above, together with \eqref{eq:f0}  and \eqref{eq:gtvchange}, implies
\begin{equation*}
	\begin{aligned}
		&\frac{1}{s}\bigg(F_0(\vx^*+s \vv^*)-F_0(\vx^*)\bigg)=\frac{1}{s}\bigg( f_0(\vx^*+s \vv^*)- f_0(\vx^*)+g(\vx^*+s \vv^*)-g(\vx^*)\bigg)\\
		=&\frac{1}{s}\bigg( \sum_{i=1}^m\big(f_i(\vx_i^*+s \vv_i^*)- f_i(\vx_i^*)\big)+g(\vx^*+s \vv^*)-g(\vx^*)\bigg)\\
		\leq& \frac{1}{s}\left(50\pi s\epsilon\sqrt{m}\|\bar{\vxi}\|_1-\beta s\|\bar{\vxi}\|_1\right)= \left(50\pi \epsilon\sqrt{m}-\beta \right)\|\bar{\vxi}\|_1.
	\end{aligned}
\end{equation*}
Taking the limit of the left-hand side of the inequality above as $s$ approaching zero, we have 
$$\left(50\pi \epsilon\sqrt{m}-\beta \right)\|\bar{\vxi}\|_1\geq F_0^{\prime}(\vx^*;\vv).$$ Thus by \eqref{eq:dderivegeq0epsilon} and the choice of $\vv^*$, we have 
$$
\left(50\pi \epsilon\sqrt{m}-\beta \right)\|\bar{\vxi}\|_1\geq -\widehat\epsilon\|\vv^*\|\geq -\widehat\epsilon\sqrt{\bar{i}}\|\bar{\vxi}\|
\geq -\widehat\epsilon\sqrt{m}\|\bar{\vxi}\|_1.
$$
This leads to a contradiction as $\beta > (50\pi+1)\epsilon\sqrt{m}$ from \eqref{eq:betachoice} and $\widehat\epsilon \le \epsilon$.  Therefore, the claim $\bar\vA\vx^*=\mathbf{0}$ is true. Thus the second condition in~\eqref{eq:KKTe} indicates $\|\vH\vx^*\|\leq\widehat\epsilon$, and we complete the proof. %it holds that $\vx^*_{i}=\vx^*_{i +1}$ for $i\in\mathcal{M}$, and thus $\bar\vA\vx^*=\mathbf{0}$ and $\|\vH\vx^*\|\leq\epsilon$. Recall inequality~\eqref{eq:KKTHe}. We have proved that $\vx^*$ an $\epsilon$-stationary point of problem~\eqref{eq:model3}.
	\null\hfill$\square$
\endproof

By Lemma~\ref{cor:kktvio2}, if $\vx^*$ is not an $\epsilon$-stationary point of the auxiliary instance~\eqref{eq:model3}, it cannot be an $\epsilon$-stationary point of instance~$\mathcal{P}$. In other words, the number of oracles needed to find an $\epsilon$-stationary point of~$\mathcal{P}$ is at least the number of oracles needed to find an $\epsilon$-stationary point of \eqref{eq:model3}. %\comm{The following sentence may not be proper: (AP) is not an instance.} 
Note that the auxiliary instance~\eqref{eq:model3} of instance~$\mathcal{P}$ is the worst-case instance used in Sun and Hong~\cite{sun2019distributed} to establish the lower-bound complexity for affinely constrained smooth optimization. In fact, according to Sun and Hong~\cite{sun2019distributed}, any algorithm that can access $\nabla f_0$ and matrix-vector multiplication with $\vH$ and $\vH^\top$ at any historical solutions needs at least  $\Theta({\kappa(\vH) L_f \Delta_{f_0}} \epsilon^{-2})$ oracles to find an $\epsilon$-stationary point of \eqref{eq:model3}, where $\Delta_{f_0}:= f_0(\vx^{(0)}) - \inf_\vx f_0(\vx)$. However, we cannot directly apply the lower bound here because the instance~$\mathcal{P}$ that we consider has both affine constraints and a non-smooth term $g$, and FOMs in Algorithm Class~\ref{algo:class-1} cannot apply $\vH$ and~$\vH^\top$ for matrix-vector multiplications but instead can use $\vA$ and $\vA^{\top}$ as well as $\prox_{\eta g}$ for any $\eta>0$. Notice that as $\eta\to \infty$, for instance~$\mathcal{P}$, $\prox_{\eta g}(\hat\vx)$ will reduce to the projection of $\hat\vx$ onto the set $\{\vx: \bar\vA\vx =~\vzero\}$. Next, we show that even with the stronger oracle, any FOM in Algorithm Class~\ref{algo:class-1} still needs at least $\Theta({\kappa(\vH) L_f \Delta_{f_0}} \epsilon^{-2})$ oracles to find an $\epsilon$-stationary point of the auxiliary instance~\eqref{eq:model3} of instance~$\mathcal{P}$.

To do so, we need the following lemma, which is essentially the same as that in Sun and Hong~\cite[Lemma 3.3]{sun2019distributed}. Due to the difference of coefficients, we provide a proof in Appendix~\ref{sec:proof-lem} for the sake of completeness. 

\begin{lemma}
	\label{lem:nablaf}
	Let $\{f_i\}_{i=1}^m$ be defined in~\eqref{eq:fi} with $\epsilon>0$. For any $\vz\in\RR^{\bar{d}}$, if $|[\vz]_{\bar{j}}|<\frac{150\pi\epsilon}{\sqrt{m} L_f}$ for some $\bar{j} \in \{1,2,\ldots,{\bar{d}}\}$, then %it holds that
	$
	\left\|\frac{1}{m}\sum_{i=1}^m\nabla f_i(\vz)\right\| > \frac{2 \epsilon}{\sqrt{m}}.
	$
\end{lemma}

The following lemma provides a lower bound to the stationarity measure of a point $\vx$ as a solution to problem~\eqref{eq:model3}. Its proof is given in  Appendix~\ref{sec:proof-lem}. %by with $\left\|\frac{1}{m}\sum_{i=1}^m\nabla f_i(\bar\vx)\right\|$ where $\bar{\vx}= \frac{1}{m} \sum_{i=1}^m \vx_i$.

\begin{lemma}
	\label{lem:kktvio}
	Let $\vx\in\mathbb{R}^d$ be given in~\eqref{eq:xblock}, $\vH$ in~\eqref{eq:matrixAstar}, and $\{f_i\}_{i=0}^m$ in~\eqref{eq:fi} and~\eqref{eq:f0}. %Set $\bar{\vx}:= \frac{1}{m} \sum_{i=1}^m \vx_i$. 
	Then
	$$
	\max\left\{\left\|\vH\vx\right\|, \min _{\vgamma}\left\|\nabla  f_0(\vx)+\vH^{\top} \vgamma\right\|\right\} \geq \frac{\sqrt{m}}{2}\left\|\frac{1}{m}\sum_{i=1}^m\nabla f_i(\bar{\vx})\right\|, \text{ with } \bar{\vx}:= \frac{1}{m} \sum_{i=1}^m \vx_i.
	%\max\left\{\left\|\vv^* -\mathbf{Proj}_{\vH\vx=0}(\vx^*)\right\|, \min _{\gamma}\left\|\nabla  f(\vx)+\vH^{\top} \gamma\right\|\right\} \geq \frac{\sqrt{m}}{2}\left\|\frac{1}{m}\sum_{i=1}^m\nabla f_i(\bar{\vx})\right\|.
	$$
\end{lemma}

The previous two lemmas imply that if there exists $\bar{j}\in\{1,2,\ldots,{\bar{d}}\}$ such that $[\bar{\vx}]_{\bar{j}}=0$, where $\bar{\vx}~=~\frac{1}{m} \sum_{i=1}^m \vx_i$, then $\vx$ cannot be an $\epsilon$-stationary point of the auxiliary problem~\eqref{eq:model3} of instance~$\mathcal{P}$. 

\subsection{Lower Bound of Oracle Complexity}

In this subsection, we provide a lower bound of the oracle complexity of FOMs in Algorithm Class~\ref{algo:class-1} for solving Problem Class~\ref{prob:class-2}, by showing that %any algorithm satisfying Assumption~\ref{ass:linearspan} will need 
a large number of calls to $\mathrm{ORACLE}_1$ will be needed to find a (near) $\epsilon$-stationary point of instance~$\mathcal{P}$. 	

For any integer $t\ge0$, let 
\begin{equation}\label{eq:t-th-iter}
	\vx^{(t)}=(\vx_1^{(t)\top}, \ldots,\vx_m^{(t)\top})^\top \text{ with each }\vx_i^{(t)}\in\RR^{\bar d}, \text{ and } \bar{\vx}^{(t)}= \frac{1}{m} \sum_{i=1}^m \vx_i^{(t)}
\end{equation}
%$\vx^{(t)}=(\vx_1^{(t)\top}, \ldots,\vx_m^{(t)\top})^\top$ 
be the $t$-th iterate of an algorithm and the block average. %$\bar{\vx}^{(t)}= \frac{1}{m} \sum_{i=1}^m \vx_i^{(t)}$.  
Without loss of generality, we assume $\vx^{(0)}=\mathbf{0}$. Otherwise, we can %change variable as follows: we 
change $f_0(\vx)$  to $f_0(\vx-\vx^{(0)})$, $g(\vx)$ to  $g(\vx-\vx^{(0)})$, and $\vA\vx=\vb$ to  $\vA(\vx-\vx^{(0)})=\vb$ in instance $\cP$. Then, the resulting instance becomes $\min_{ \vx\in\RR^d}  F_0( \vx):= f_0( \vx-\vx^{(0)}) + g( \vx-\vx^{(0)}), \ \st  \vA (\vx-\vx^{(0)}) +  \vb=\mathbf{0}$, %$\vx$ in instance~$\cP$ to $\vx-\vx^{(0)}$, 
and we can process the rest of the proof with this new instance. %\YX{Add an explanation why shifting the solution will not change our claimed results.} 
Our lower bound will be established based on the fact that, if $t$ is not large enough, $\textnormal{supp}(\bar{\vx}^{(t)})\subset\{1,\dots,\bar{j}-1\}$ for some $\bar{j}\in\{1,2,\ldots,{\bar{d}}\}$. Hence, $[\bar{\vx}^{(t)}]_{\bar{j}}=0$ and $\vx^{(t)}$ cannot be a (near) $\epsilon$-stationary point due to Lemmas \ref{lem:nablaf} and \ref{lem:kktvio}.
From the way that the iterates are generated in Algorithm Class~\ref{algo:class-1},  $\textnormal{supp}(\bar{\vx}^{(t)})$ is influenced by $\textnormal{supp}(\nabla f_i(\vx_i^{(t)}))$, which is characterized by the following lemma.
					
\begin{lemma}
	\label{lem:iterateguess}
Let $\{f_i\}_{i=1}^m$ be defined in~\eqref{eq:fi}. {Given any $\bar{j}\in\{1,\dots,\bar{d}\}$ and $\vz\in\RR^{\bar{d}}$ with $\textnormal{supp}(\vz)\subset\{1,\dots,\bar{j}-1\}$, it holds that  } %\YX{Do not use italic supp. Instead use $\mathrm{supp}$. Denote $\bar \vz = \frac{\sqrt{m} L_f \vz}{150\pi\epsilon}$ to shorten a few equations.}
%\YX{Change all $\vz_j$ to $[\vz]_j$ for notational consistency.}
\begin{enumerate}
	\item   When $\bar{j}=1$, $\textnormal{supp}(\vz)=\emptyset$, $\vz=\mathbf{0}$, and
	$
	\textnormal{supp}(\nabla f_i(\vz))\subset\{1\}$, for any $i\in[1,m];
	$ %for $i=1,\dots,m$.
	\item  When $\bar{j}$ is even, 
	$$
	\textnormal{supp}(\nabla f_i(\vz))\subset\left\{
	\begin{array}{ll}
	\{1,\dots,\bar{j}\}, & \text{ if }\ i\in\left[1, \frac{m}{3}\right],\\[0.1cm]
	\{1,\dots,\bar{j}-1\}, & \text{ if }\ i\in\left[\frac{m}{3}+1, m\right];
	\end{array}
	\right.
	$$
	\item  When $\bar{j}$ is odd and $\bar{j}\neq 1$, 
	$$
	\textnormal{supp}(\nabla f_i(\vz))\subset\left\{
	\begin{array}{ll}
		\{1,\dots,\bar{j}-1\}, & \text{ if }\ i\in\left[1, \frac{2m}{3}\right],\\[0.1cm]
		\{1,\dots,\bar{j}\}, & \text{ if }\ i\in\left[\frac{2m}{3}+1, m\right].
	\end{array}
	\right.
	$$
\end{enumerate}
%When $\bar{j}=1$ (so $\vz=\mathbf{0}$), the structure above implies $\textnormal{supp}(\nabla f_i(\vz))\subset\{1\}$ for $i=1,\dots,m$.
%When $\bar{j}$ is even, the structure above implies $\textnormal{supp}(\nabla f_i(\vz))\subset\{1,\dots,\bar{j}+1\}$ for $\left[1, \frac{m}{3}\right]$.
%$\textnormal{supp}(\nabla f_i(\vz))\subset\{1,\dots,\bar{j}\}$ for $i=1,\dots,m$.
\end{lemma}

\proof{Proof of Lemma \ref{lem:iterateguess}.}
	Denote $\bar \vz = \frac{\sqrt{m} L_f \vz}{150\pi\epsilon}$ and recall definition~\eqref{eq:grad-f-i-chain}.
 	%Consider $\bar{j}\in\{1,\dots,\bar{d}\}$ and $\vz\in\RR^{\bar{d}}$ with $\textnormal{supp}(\vz)\subset\{1,\dots,\bar{j}-1\}$. 
% Recall the relation
 %		
 	By Lemma~\ref{lem:functions}(a), we have %that,  for $j\geq \bar{j}$, 
 	\begin{eqnarray}
 		\label{eq:zeroPsi}
 		\Psi\left(-[\bar \vz]_j\right)=\Psi\left([\bar \vz]_j\right)=
 		\Psi'\left(-[\bar \vz]_j\right)=\Psi'\left([\bar \vz]_j\right)=0, \forall\, j\geq \bar{j}.
 	\end{eqnarray}
Therefore, by definitions~\eqref{eq:deriveh3}--\eqref{eq:deriveh2}, the support of $\vz$ leads to the following structure of 
	$[\nabla f_i(\vz)]_j$ for $j\geq \bar{j}$:
	
	\begin{itemize}
			\item 
	If $j=1$,
	$
	[\nabla f_i(\vz)]_j=	-\Psi(1) \Phi^{\prime}\left([\bar \vz]_j\right)$, for any $i \in\left[1, m\right];
	$
		\item 
	If $j$ is even, 
	$$
[\nabla f_i(\vz)]_j=\left\{
\begin{array}{ll}
-\frac{6\epsilon}{\sqrt{m}}\left[\Psi\left(-[\bar \vz]_{j-1}\right) \Phi^{\prime}\left(-[\bar \vz]_{j}\right)+\Psi\left([\bar \vz]_{j-1}\right) \Phi^{\prime}\left([\bar \vz]_{j}\right)\right],& \text{ for }  i \in\left[1, \frac{m}{3}\right],\\[0.1cm]
0, & \text{ for } i \in\left[\frac{m}{3}+1, m\right];
\end{array}
\right.
	$$
	\normalsize
	\item 
	If $j$ is odd and $j\neq 1$, 
	%\small
	$$
	[\nabla f_i(\vz)]_j=\left\{
	\begin{array}{ll}
		0, & \text{ for } i   \in\left[1, \frac{2m}{3}\right],\\[0.1cm]
-\frac{6\epsilon}{\sqrt{m}}\left[\Psi\left(-[\bar \vz]_{j-1}\right) \Phi^{\prime}\left(-[\bar \vz]_{j}\right)+\Psi\left([\bar \vz]_{j-1}\right) \Phi^{\prime}\left([\bar \vz]_{j}\right)\right],& \text{ for } i  \in\left[\frac{2m}{3}+1, m\right].
	\end{array}
	\right.
	$$
	\normalsize
\end{itemize} 	
Since $\Psi\left(-[\bar \vz]_{j-1}\right)=\Psi\left([\bar \vz]_{j-1}\right)=0$ for any $j > \bar{j}$, the structures above imply $[\nabla f_i(\vz)]_j = 0,\forall\, j > \bar j$ and thus give the desired claims.
%that $\textnormal{supp}(\nabla f_i(\vz))\subset\{1,\dots,\bar{j}\}$ when $\textnormal{supp}(\vz)\subset\{1,\dots,\bar{j}-1\}$.
	\null\hfill$\square$
\endproof

According to the structure of $\vA$ given in~\eqref{eq:AandAbar},  $\textnormal{supp}((\vA^{\top}\vA\vx)_i)$ is determined by $\textnormal{supp}(\vx_{i-1})$, $\textnormal{supp}(\vx_{i})$ and $\textnormal{supp}(\vx_{i+1})$. Also, $\textnormal{supp}(\prox_{\eta g}(\vx))$ has a similar property according to the definition of $g$  in~\eqref{eq:g}. These properties are formally stated in the following lemma.

\begin{lemma}
	\label{lem:supp}
Let $\vx$ be the structured vector given in~\eqref{eq:xblock}, $\vA$ in~\eqref{eq:AandAbar}, and $g$ be given in~\eqref{eq:g}. Define $\vx_0=\vx_{m+1}=\mathbf{0}\in\RR^{{\bar{d}}}$. The following statements hold:
	\begin{enumerate}
		\item[\textnormal{(a)}] Let  $\widehat\vx=\vA^{\top}\vA\vx=(\widehat\vx_1\zz, \ldots,\widehat\vx_m\zz)\zz$ with each $\widehat\vx_i\in\RR^{{\bar{d}}}$.  Then %It then yields that 
		\begin{equation}\label{eq:supp-relation}
		\textnormal{supp}(\widehat\vx_i)\subset \textnormal{supp}(\vx_{i-1})\cup\textnormal{supp}(\vx_{i})\cup\textnormal{supp}(\vx_{i+1}),\, \forall\, i\in [1, m].
		\end{equation}		
		
		\item[\textnormal{(b)}] For any $\eta>0$, let $\widetilde\vx=\prox_{\eta g}(\vx)=(\widetilde\vx_1\zz, \ldots,\widetilde\vx_m\zz)\zz$ with each $\widetilde\vx_i\in\RR^{{\bar{d}}}$. Then %It then yields that 
		$$\textnormal{supp}(\widetilde\vx_i)\subset \textnormal{supp}(\vx_{i-1})\cup\textnormal{supp}(\vx_{i})\cup\textnormal{supp}(\vx_{i+1}),\, \forall\, i\in [1, m].$$
	\end{enumerate}
\end{lemma}

\proof{Proof of Lemma \ref{lem:supp}.}
	(a) %Note that $\vA^{\top}\vb=\mathbf{0}$ for instance~$\mathcal{P}$. Let $\vx=(\vx_1\zz, \ldots,\vx_m\zz)\zz$ with $\vx_i\in\RR^{{\bar{d}}}$. Let  $\widehat\vx=\vA^{\top}\vA\vx=(\widehat\vx_1\zz, \ldots,\widehat\vx_m\zz)\zz$ with $\widehat\vx_i\in\RR^{{\bar{d}}}$. We define $\vx_0=\vx_{m+1}=\mathbf{0}$. By the definition of $\vA$ in~\eqref{eq:AandAbar}, we have 
The relation in~\eqref{eq:supp-relation} immediately follows from the observation
\begin{equation}\label{eq:ata-struc}
\newcommand{\AtA}{\left[ 
\begin{array}{cccc}
\vB & & & \\
& \vB & & \\
& & \ddots &\\
& & & \vB
\end{array}
\right]
}
\newcommand{\Bmat}{\left[ 
\begin{array}{rrrrr} 
\vI_{\bar d} & - \vI_{\bar d} & & & \\
- \vI_{\bar d} & 2 \vI_{\bar d} & - \vI_{\bar d} & \\
& \ddots & \ddots & \ddots & \\
&  &	- \vI_{\bar d} & 2 \vI_{\bar d} & - \vI_{\bar d}  \\
& & & - \vI_{\bar d} & \vI_{\bar d}
		\end{array} \right]
}
\vA^\top\vA =  \left. \AtA 
		\right\}\text{$3m_2$ blocks}, \text{ with } 
\vB = m^2 L_f^2 
\left.
		\,\smash[b]{\underbrace{\!\Bmat\!}_{\textstyle\text{$m_1$ blocks}}}\,
		\right\}\text{$m_1$ blocks.}
		\vphantom{\underbrace{\Bmat}_{\textstyle\text{$m$ blocks}}}	
\end{equation}
%Hence,~\eqref{eq:supp-relation} immediately follows.	
%	\begin{eqnarray}
%		\nonumber
%		\textnormal{supp}(\widehat\vx_i)&\subset&
%		\left\{
%		\begin{array}{ll}
%			\textnormal{supp}(\vx_{i-1})\cup\textnormal{supp}(\vx_{i})\cup\textnormal{supp}(\vx_{i+1})&\text{ if }i-1,i\notin\mathcal{M},\\
%			\textnormal{supp}(\vx_{i})\cup\textnormal{supp}(\vx_{i+1})&\text{ if }i-1\in\mathcal{M},\\
%			\textnormal{supp}(\vx_{i-1})\cup\textnormal{supp}(\vx_{i})&\text{ if }i\in\mathcal{M}.
%		\end{array}
%		\right.%\\\label{eq:suppA}
%		%&\subset&
%	\end{eqnarray}
	
	(b) Given any $x_1$ and $x_2$ in $\mathbb{R}$ and any $c>0$, consider the following optimization problem in $\mathbb{R}^2$:
	\begin{eqnarray*}
		(\widetilde x_1, \widetilde x_2)=\argmin_{z_1,z_2\in\mathbb{R}}\frac{1}{2}(z_1-x_1)^2+\frac{1}{2}(z_2-x_2)^2+c|z_1-z_2|.
	\end{eqnarray*}
	The optimal solution of this problem is 
	\begin{eqnarray}
		\label{eq:2dproxg}
		(\widetilde x_1, \widetilde x_2)=\left\{
		\begin{array}{ll}
			( (x_1+x_2)/2,  (x_1+x_2)/2),&\text{ if } |x_1-x_2|\leq 2c\\[0.1cm]
			( x_1-c\cdot\text{sign}(x_1-x_2),  x_2+c\cdot\text{sign}(x_1-x_2)),&\text{ if } |x_1-x_2|> 2c.
		\end{array}
		\right.
	\end{eqnarray}
	Recall the definition of $g$ in \eqref{eq:g} and $\textbf{prox}_{\eta g}$, we obtain that
	$$
	\prox_{\eta g}(\vx) = \argmin_{\vy} \eta\beta\sum_{i\in\mathcal{M}}  \|\vx_{i }-\vx_{i+1}\|_1 + \frac{1}{2}\|\vx-\vy\|^2.
	$$
	It then holds that 
	\begin{eqnarray*}
		\widetilde\vx_i = 
		\left\{
		\begin{array}{ll}
			 \argmin_{\vy_i} \frac{1}{2}\|\vx_i-\vy_i\|^2, & \text{ if }i-1,i\notin\mathcal{M},\\[1mm]
			\argmin_{\vy_i} \eta\beta \|\vx_{i }-\vx_{i-1}\|_1 + \frac{1}{2}\|\vx_i-\vy_i\|^2, &\text{ if }i-1\in\mathcal{M},\\[1mm]
			\argmin_{\vy_i} \eta\beta \|\vx_{i }-\vx_{i+1}\|_1 + \frac{1}{2}\|\vx_i-\vy_i\|^2 ,&\text{ if }i\in\mathcal{M}.
		\end{array}
		\right.
	\end{eqnarray*}
	%For any given $\eta>0$, let $\widetilde\vx=\prox_{\eta g}(\vx)=(\widetilde\vx_1\zz, \ldots,\widetilde\vx_m\zz)\zz$ with $\widetilde\vx_i\in\RR^{{\bar{d}}}$. 
%\YX{To make the following claim easier to understand, add one more step by formulating a two-dimensional problem about $([\vx_i]_j, [\vx_{i+1}]_{j})$ for each $i$ and $j$. } 	
	By~\eqref{eq:2dproxg} and the separability property of $\|\cdot\|_1$ and $\|\cdot\|^2$, %the fact that $m_1\geq2$, and the definition of $g$  in~\eqref{eq:g}, 
	we have that for any $j\in\{1,\dots,\bar{d}\}$, 
	\begin{eqnarray*}
		[\widetilde\vx_i]_j=
		\left\{
		\begin{array}{ll}
			[\vx_i]_j, &\text{ if }i-1,i\notin\mathcal{M},\\\null
			([\vx_{i-1}]_j+[\vx_i]_j)/2, &\text{ if }i-1\in\mathcal{M}\text{ and }\left|[\vx_{i-1}]_j-	[\vx_{i}]_j\right|\leq \eta\beta,\\[1mm]
			([\vx_{i}]_j+[\vx_{i+1}]_j)/2, &\text{ if }i\in\mathcal{M}\text{ and }\left|[\vx_{i}]_j-	[\vx_{i+1}]_j\right|\leq \eta\beta,\\[1mm]
			[\vx_i]_j+\eta\beta\cdot\text{sign}([\vx_{i-1}]_j-[\vx_i]_j), &\text{ if }i-1\in\mathcal{M}\text{ and }\left|[\vx_{i-1}]_j-	[\vx_{i}]_j\right|> \eta\beta,\\[1mm]
			[\vx_i]_j-\eta\beta\cdot\text{sign}([\vx_{i-1}]_j-[\vx_i]_j), &\text{ if }i\in\mathcal{M}\text{ and }\left|[\vx_{i-1}]_j-	[\vx_{i}]_j\right|> \eta\beta,
		\end{array}
		\right.
	\end{eqnarray*}
	which implies 
	\begin{eqnarray*}
		\textnormal{supp}(\widetilde\vx_i)&\subset&
		\left\{
		\begin{array}{ll}
			\textnormal{supp}(\vx_{i}),&\text{ if }i-1,i\notin\mathcal{M},\\
			\textnormal{supp}(\vx_{i-1})\cup\textnormal{supp}(\vx_{i}),&\text{ if }i-1\in\mathcal{M},\\
			\textnormal{supp}(\vx_{i})\cup\textnormal{supp}(\vx_{i+1}),&\text{ if }i\in\mathcal{M}.
		\end{array}
		\right.%\%\\label{eq:suppg}
		%&\subset&\textnormal{supp}(\vx_{i-1})\cup\textnormal{supp}(\vx_{i})\cup\textnormal{supp}(\vx_{i+1}).
	\end{eqnarray*}
	The proof is then completed.
	\null\hfill$\square$
\endproof

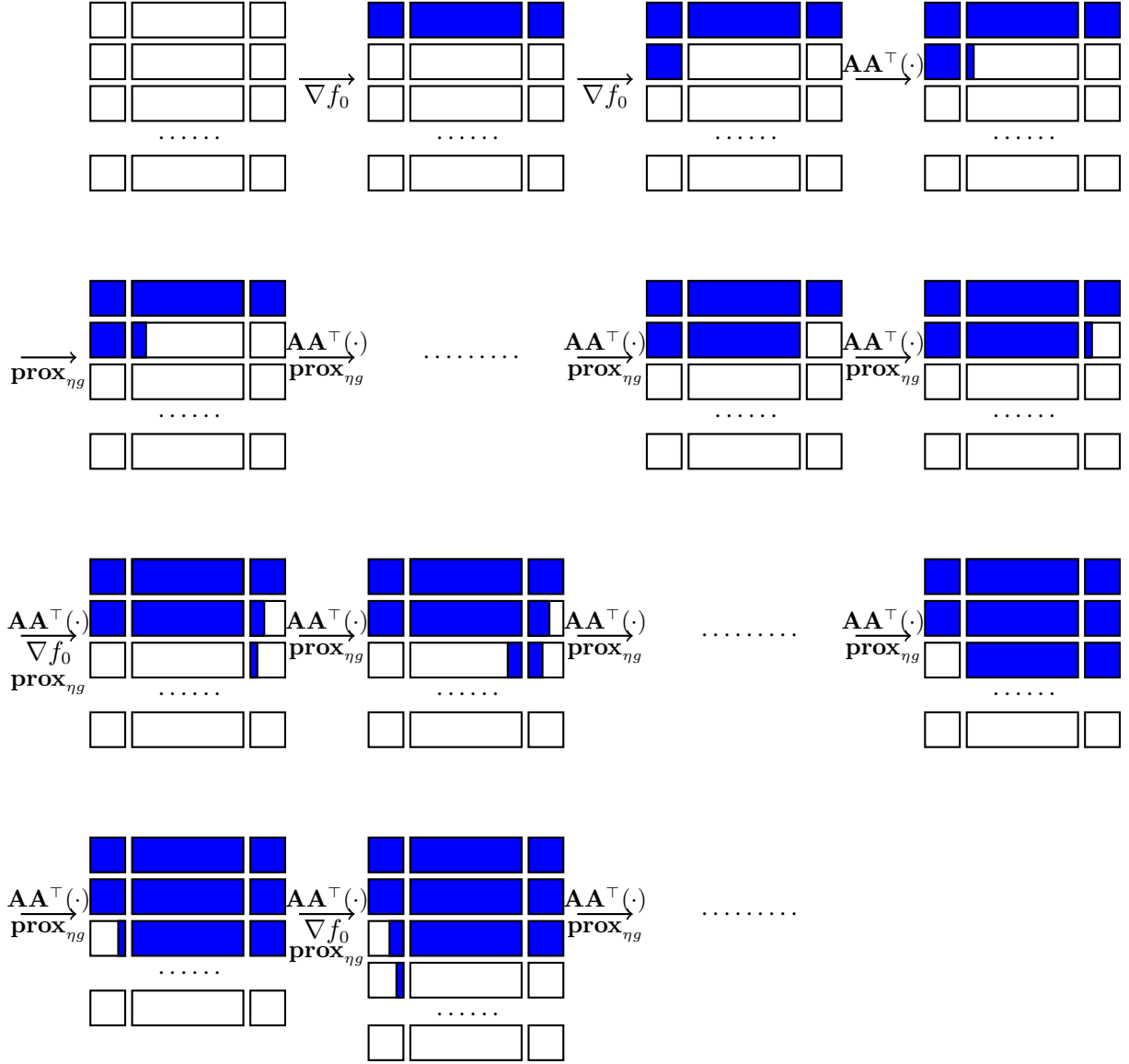
\begin{figure}
	\centering
	\begin{tikzpicture}
		% First row of boxes
		\draw[thick ] (0, 1) rectangle (0.5, 1.5);  % filled black box
		\draw[thick ] (0.6, 1) rectangle (2.2, 1.5);  % large rectangle
		%\draw[thick, fill=blue] (0.6, 1) rectangle (0.7, 1.5);  % part of the large rectangle filled black
		\draw[thick ] (2.3, 1) rectangle (2.8, 1.5);  % small rectangle on right
		% Second row of boxes
		\draw[thick] (0, 0.4) rectangle (0.5, 0.9);  % small rectangle on left
		\draw[thick] (0.6, 0.4) rectangle (2.2, 0.9);  % large rectangle
		\draw[thick] (2.3, 0.4) rectangle (2.8, 0.9);  % small rectangle on right
		
		\draw[thick] (0, -0.2) rectangle (0.5, 0.3);  % small rectangle on left
		\draw[thick] (0.6, -0.2) rectangle (2.2, 0.3);  % large rectangle
		\draw[thick] (2.3, -0.2) rectangle (2.8, 0.3);  % small rectangle on right
		
		\node at (1.45, -0.45) {$\cdots\cdots$};

		\draw[thick] (0, -1.2) rectangle (0.5, -0.7);  % small rectangle on left
		\draw[thick] (0.6, -1.2) rectangle (2.2, -0.7);  % large rectangle
		\draw[thick] (2.3, -1.2) rectangle (2.8, -0.7);  % small rectangle on right
		
		\draw[->, thick] (3, 0.4) -- (3.8, 0.4);
		\node at (3.4, 0.15) {$\nabla f_0$};
		
		\draw[thick, fill=blue] (4, 1) rectangle (4.5, 1.5);  % filled blue box
		\draw[thick, fill=blue] (4.6, 1) rectangle (6.2, 1.5);  % large rectangle
		%\draw[thick, fill=blue] (4.6, 1) rectangle (4.7, 1.5);  % part of the large rectangle filled blue
		\draw[thick, fill=blue] (6.3, 1) rectangle (6.8, 1.5);  % small rectangle on right
		
		% Second row of boxes
		\draw[thick ] (4, 0.4) rectangle (4.5, 0.9);  % small rectangle on left
		\draw[thick] (4.6, 0.4) rectangle (6.2, 0.9);  % large rectangle
		\draw[thick] (6.3, 0.4) rectangle (6.8, 0.9);  % small rectangle on right
		
		\draw[thick] (4, -0.2) rectangle (4.5, 0.3);  % small rectangle on left
		\draw[thick] (4.6, -0.2) rectangle (6.2, 0.3);  % large rectangle
		\draw[thick] (6.3, -0.2) rectangle (6.8, 0.3);  % small rectangle on right
		
		\node at (5.45, -0.45) {$\cdots\cdots$};
		
		\draw[thick] (4, -1.2) rectangle (4.5, -0.7);  % small rectangle on left
		\draw[thick] (4.6, -1.2) rectangle (6.2, -0.7);  % large rectangle
		\draw[thick] (6.3, -1.2) rectangle (6.8, -0.7);  % small rectangle on right
		
		\draw[->, thick] (7, 0.4) -- (7.8, 0.4);
		\node at (7.4, 0.15) {$\nabla f_0$}; 
		
		\draw[thick, fill=blue] (8, 1) rectangle (8.5, 1.5);  % filled blue box
		\draw[thick, fill=blue] (8.6, 1) rectangle (10.2, 1.5);  % large rectangle
		\draw[thick, fill=blue] (10.3, 1) rectangle (10.8, 1.5);  % small rectangle on right
		
		% Second row of boxes
		\draw[thick, fill=blue] (8, 0.4) rectangle (8.5, 0.9);  % small rectangle on left
		\draw[thick] (8.6, 0.4) rectangle (10.2, 0.9);  % large rectangle 
		\draw[thick] (10.3, 0.4) rectangle (10.8, 0.9);  % small rectangle on right
		
		\draw[thick] (8, -0.2) rectangle (8.5, 0.3);  % small rectangle on left
		\draw[thick] (8.6, -0.2) rectangle (10.2, 0.3);  % large rectangle
		\draw[thick] (10.3, -0.2) rectangle (10.8, 0.3);  % small rectangle on right
		
		\node at (9.45, -0.45) {$\cdots\cdots$};
		
		\draw[thick] (8, -1.2) rectangle (8.5, -0.7);  % small rectangle on left
		\draw[thick] (8.6, -1.2) rectangle (10.2, -0.7);  % large rectangle
		\draw[thick] (10.3, -1.2) rectangle (10.8, -0.7);  % small rectangle on right
		
		\draw[->, thick] (11, 0.4) -- (11.8, 0.4);
		
		\node at (11.4, 0.6) {\small $\vA\vA\zz(\cdot)$};
		
		\draw[thick, fill=blue] (12, 1) rectangle (12.5, 1.5);  % filled blue box
		\draw[thick, fill=blue] (12.6, 1) rectangle (14.2, 1.5);  % large rectangle
		%\draw[thick, fill=blue] (12.6, 1) rectangle (12.8, 1.5);  % part of the large rectangle filled blue
		\draw[thick, fill=blue] (14.3, 1) rectangle (14.8, 1.5);  % small rectangle on right
		
		% Second row of boxes
		\draw[thick, fill=blue] (12, 0.4) rectangle (12.5, 0.9);  % small rectangle on left
		\draw[thick] (12.6, 0.4) rectangle (14.2, 0.9);  % large rectangle
		\draw[thick, fill=blue] (12.6, 0.4) rectangle (12.7, 0.9);  % part of the large rectangle filled blue
		\draw[thick] (14.3, 0.4) rectangle (14.8, 0.9);  % small rectangle on right
		
		\draw[thick] (12, -0.2) rectangle (12.5, 0.3);  % small rectangle on left
		\draw[thick] (12.6, -0.2) rectangle (14.2, 0.3);  % large rectangle
		\draw[thick] (14.3, -0.2) rectangle (14.8, 0.3);  % small rectangle on right
		
		\node at (13.45, -0.45) {$\cdots\cdots$};
		
		\draw[thick] (12, -1.2) rectangle (12.5, -0.7);  % small rectangle on left
		\draw[thick] (12.6, -1.2) rectangle (14.2, -0.7);  % large rectangle
		\draw[thick] (14.3, -1.2) rectangle (14.8, -0.7);  % small rectangle on right

		\draw[->, thick] (-1, -3.6) -- (-0.2, -3.6); 
		\node at (-0.6, -3.9) {\small $\prox_{\eta {g}}$};
		
		\draw[thick, fill=blue] (0, -3) rectangle (0.5, -2.5);  % filled black box
		\draw[thick, fill=blue] (0.6, -3) rectangle (2.2, -2.5);  % large rectangle  rectangle filled black
		\draw[thick, fill=blue] (2.3, -3) rectangle (2.8, -2.5);  % small rectangle on right
		
		% Second row of boxes
		\draw[thick, fill=blue] (0, -3.6) rectangle (0.5, -3.1);  % small rectangle on left
		\draw[thick] (0.6, -3.6) rectangle (2.2, -3.1);  % large rectangle
		\draw[thick , fill=blue] (0.6, -3.6) rectangle (0.8, -3.1);  % large rectangle
		\draw[thick ] (2.3, -3.6) rectangle (2.8, -3.1);  % small rectangle on right
		
		\draw[thick] (0, -4.2) rectangle (0.5, -3.7);  % small rectangle on left
		\draw[thick] (0.6, -4.2) rectangle (2.2,-3.7);  % large rectangle
		\draw[thick] (2.3, -4.2) rectangle (2.8,-3.7);  % small rectangle on right
		
		\node at (1.45, -4.45) {$\cdots\cdots$};
		
		\draw[thick] (0, -5.2) rectangle (0.5, -4.7);  % small rectangle on left
		\draw[thick] (0.6, -5.2) rectangle (2.2, -4.7);  % large rectangle
		\draw[thick] (2.3, -5.2) rectangle (2.8, -4.7);  % small rectangle on right
		
		\draw[->, thick] (3, -3.6) -- (3.8, -3.6);
		\node at (3.4, -3.4) {\small $\vA\vA\zz(\cdot)$};
		\node at (3.4, -3.85) {\small $\prox_{\eta {g}}$};
		
		\node at (5.5, -3.6) {$\cdots\cdots\cdots$};
		
		\draw[->, thick] (7, -3.6) -- (7.8, -3.6);
		\node at (7.4, -3.4) {\small $\vA\vA\zz(\cdot)$};
		\node at (7.4, -3.85) {\small $\prox_{\eta {g}}$};
		
		\draw[thick, fill=blue] (8, -3) rectangle (8.5, -2.5);  % filled blue box
		\draw[thick] (8.6, -3) rectangle (10.2, -2.5);  % large rectangle
		\draw[thick, fill=blue] (8.6, -3) rectangle (10.2, -2.5);  % part of the large rectangle filled blue
		\draw[thick, fill=blue] (10.3, -3) rectangle (10.8, -2.5);  % small rectangle on right
		
		% Second row of boxes
		\draw[thick, fill=blue] (8, -3.6) rectangle (8.5, -3.1);  % small rectangle on left
		\draw[thick, fill=blue] (8.6, -3.6) rectangle (10.2, -3.1);  % large rectangle
		\draw[thick] (10.3, -3.6) rectangle (10.8, -3.1);  % small rectangle on right
		
		\draw[thick] (8, -4.2) rectangle (8.5, -3.7);  % small rectangle on left
		\draw[thick] (8.6, -4.2) rectangle (10.2, -3.7);  % large rectangle
		\draw[thick] (10.3, -4.2) rectangle (10.8, -3.7);  % small rectangle on right
		
		\node at (9.45, -4.45) {$\cdots\cdots$};
		
		\draw[thick] (8, -5.2) rectangle (8.5, -4.7);  % small rectangle on left
		\draw[thick] (8.6, -5.2) rectangle (10.2, -4.7);  % large rectangle
		\draw[thick] (10.3, -5.2) rectangle (10.8, -4.7);  % small rectangle on right
		
		\draw[->, thick] (11, -3.6) -- (11.8, -3.6);
		\node at (11.4, -3.85) {\small $\prox_{\eta {g}}$};
		\node at (11.4, -3.4) {\small $\vA\vA\zz(\cdot)$};
		
		\draw[thick, fill=blue] (12, -3) rectangle (12.5, -2.5);  % filled blue box
		\draw[thick, fill=blue] (12.6, -3) rectangle (14.2, -2.5);  % large rectangle
		
		\draw[thick, fill=blue] (14.3, -3) rectangle (14.8, -2.5);  % small rectangle on right
		%\draw[thick, fill=blue] (14.3, -3) rectangle (14.4, -2.5);  % small rectangle on right
		
		% Second row of boxes
		\draw[thick, fill=blue] (12, -3.6) rectangle (12.5, -3.1);  % small rectangle on left
		%\draw[thick, fill=blue] (12.4, -3.6) rectangle (12.5, -3.1);  % small rectangle on left
		\draw[thick, fill=blue] (12.6, -3.6) rectangle (14.2, -3.1);  % large rectangle
		\draw[thick ] (14.3, -3.6) rectangle (14.8, -3.1);  % small rectangle on right
		\draw[thick , fill=blue] (14.3, -3.6) rectangle (14.4, -3.1);  % small rectangle on right
		
		\draw[thick] (12, -4.2) rectangle (12.5, -3.7);  % small rectangle on left
		\draw[thick ] (12.6, -4.2) rectangle (14.2, -3.7);  % large rectangle
		\draw[thick ] (14.3, -4.2) rectangle (14.8, -3.7);  % small rectangle on right
		
		\node at (13.45, -4.45) {$\cdots\cdots$};
		
		\draw[thick] (12, -5.2) rectangle (12.5, -4.7);  % small rectangle on left
		\draw[thick ] (12.6, -5.2) rectangle (14.2, -4.7);  % large rectangle
		\draw[thick ] (14.3, -5.2) rectangle (14.8, -4.7);  % small rectangle on right

		\draw[->, thick] (-1, -7.6) -- (-0.2, -7.6);
		\node at (-0.6, -7.85) {$\nabla f_0$};
		\node at (-0.6, -7.4) {\small $\vA\vA\zz(\cdot)$};
		\node at (-0.6, -8.3) {\small $\prox_{\eta {g}}$};
		
		\draw[thick, fill=blue] (0, -7) rectangle (0.5, -6.5);  % filled black box
		\draw[thick] (0.6, -3) rectangle (2.2, -2.5);  % large rectangle
		\draw[thick, fill=blue] (0.6, -7) rectangle (2.2, -6.5);  % part of the large rectangle filled black
		\draw[thick, fill=blue] (2.3, -7) rectangle (2.8, -6.5);  % small rectangle on right
		
		% Second row of boxes
		\draw[thick, fill=blue] (0, -7.6) rectangle (0.5, -7.1);  % small rectangle on left
		%\draw[thick, fill=blue] (0.3, -7.6) rectangle (0.5, -7.1);  % small rectangle on left
		\draw[thick, fill=blue] (0.6, -7.6) rectangle (2.2, -7.1);  % large rectangle
		\draw[thick] (2.3, -7.6) rectangle (2.8, -7.1);  % small rectangle on right
		\draw[thick, fill=blue] (2.3, -7.6) rectangle (2.5, -7.1);  % small rectangle on right
		
		\draw[thick] (0, -8.2) rectangle (0.5, -7.7);  % small rectangle on left
		\draw[thick] (0.6, -8.2) rectangle (2.2, -7.7);  % large rectangle
		\draw[thick] (2.3, -8.2) rectangle (2.8, -7.7);  % small rectangle on right
		\draw[thick, fill=blue] (2.3, -8.2) rectangle (2.4, -7.7);  % small rectangle on right
		
		\node at (1.45, -8.45) {$\cdots\cdots$};
		
		\draw[thick] (0, -9.2) rectangle (0.5, -8.7);  % small rectangle on left
		\draw[thick] (0.6, -9.2) rectangle (2.2, -8.7);  % large rectangle
		\draw[thick] (2.3, -9.2) rectangle (2.8, -8.7);  % small rectangle on right
		
		\draw[->, thick] (3, -7.6) -- (3.8, -7.6);
		\node at (3.4, -7.4) {\small $\vA\vA\zz(\cdot)$};
		\node at (3.4, -7.85) {\small $\prox_{\eta {g}}$};
		
		\draw[thick, fill=blue] (4, -7) rectangle (4.5, -6.5);  % filled black box
		\draw[thick] (4.6, -7) rectangle (6.2, -6.5);  % large rectangle
		\draw[thick, fill=blue] (4.6, -7) rectangle (6.2, -6.5);  % part of the large rectangle filled black
		\draw[thick, fill=blue] (6.3, -7) rectangle (6.8, -6.5);  % small rectangle on right
		
		% Second row of boxes
		\draw[thick, fill=blue] (4, -7.6) rectangle (4.5, -7.1);  % small rectangle on left
		\draw[thick, fill=blue] (4.6,-7.6) rectangle (6.2, -7.1);  % large rectangle
		\draw[thick ] (6.3,-7.6) rectangle (6.8, -7.1);  % small rectangle on right
		\draw[thick, fill=blue ] (6.3,-7.6) rectangle (6.6, -7.1);  % small rectangle on right
		
		\draw[thick] (4, -8.2) rectangle (4.5, -7.7);  % small rectangle on left
		\draw[thick] (4.6, -8.2) rectangle (6.2, -7.7);  % large rectangle
		\draw[thick, fill=blue] (6,-8.2) rectangle (6.2, -7.7);  % large rectangle
		\draw[thick] (6.3, -8.2) rectangle (6.8, -7.7);  % small rectangle on right
		\draw[thick, fill=blue] (6.3, -8.2) rectangle (6.5, -7.7);  % small rectangle on right
		
		\node at (5.45, -8.45) {$\cdots\cdots$};
		
		\draw[thick] (4, -9.2) rectangle (4.5, -8.7);  % small rectangle on left
		\draw[thick] (4.6, -9.2) rectangle (6.2,-8.7);  % large rectangle
		\draw[thick] (6.3, -9.2) rectangle (6.8, -8.7);  % small rectangle on right
		
		\draw[->, thick] (7, -7.6) -- (7.8, -7.6);
		\node at (7.4, -7.4) {\small $\vA\vA\zz(\cdot)$};
		\node at (7.4, -7.85) {\small $\prox_{\eta {g}}$};
		
		\node at (9.5, -7.6) {$\cdots\cdots\cdots$};
		
		\draw[->, thick] (11, -7.6) -- (11.8, -7.6);
		\node at (11.4, -7.4) {\small $\vA\vA\zz(\cdot)$};
		\node at (11.4, -7.85) {\small $\prox_{\eta {g}}$};
		
		\draw[thick, fill=blue] (12, -7) rectangle (12.5, -6.5);  % filled black box
		\draw[thick] (12.6, -7) rectangle (14.2, -6.5);  % large rectangle
		\draw[thick, fill=blue] (12.6, -7) rectangle (14.2, -6.5);  % part of the large rectangle filled black
		\draw[thick, fill=blue] (14.3, -7) rectangle (14.8, -6.5);  % small rectangle on right
		
		% Second row of boxes
		\draw[thick, fill=blue] (12, -7.6) rectangle (12.5, -7.1);  % small rectangle on left
		\draw[thick, fill=blue] (12.6, -7.6) rectangle (14.2, -7.1);  % large rectangle
		%\draw[thick, fill=blue] (14.0, -7.6) rectangle (14.2, -7.1);  % part of the large rectangle filled black
		\draw[thick, fill=blue] (14.3, -7.6) rectangle (14.8, -7.1);  % small rectangle on right
		
		\draw[thick] (12, -8.2) rectangle (12.5, -7.7);  % small rectangle on left
		\draw[thick, fill=blue] (12.6, -8.2) rectangle (14.2, -7.7);  % large rectangle
		\draw[thick, fill=blue] (14.3, -8.2) rectangle (14.8, -7.7);  % small rectangle on right rectangle on right
		
		\node at (13.45, -8.45) {$\cdots\cdots$};
		
		\draw[thick] (12, -9.2) rectangle (12.5, -8.7);  % small rectangle on left
		\draw[thick] (12.6, -9.2) rectangle (14.2, -8.7);  % large rectangle
		\draw[thick] (14.3, -9.2) rectangle (14.8, -8.7);  % small rectangle on right
		
		\draw[->, thick] (-1, -11.6) -- (-0.2, -11.6);
		\node at (-0.6, -11.4) {\small $\vA\vA\zz(\cdot)$};
		\node at (-0.6, -11.85) {\small $\prox_{\eta {g}}$};
		
		\draw[thick, fill=blue] (0, -11) rectangle (0.5, -10.5);  % filled black box
		\draw[thick] (0.6, -7) rectangle (2.2, -6.5);  % large rectangle
		\draw[thick, fill=blue] (0.6, -11) rectangle (2.2, -10.5);  % part of the large rectangle filled black
		\draw[thick, fill=blue] (2.3, -11) rectangle (2.8, -10.5);  % small rectangle on right
		
		% Second row of boxes
		\draw[thick, fill=blue] (0, -11.6) rectangle (0.5, -11.1);  % small rectangle on left
		%\draw[thick, fill=blue] (0.3, -11.6) rectangle (0.5, -11.1);  % small rectangle on left
		\draw[thick, fill=blue] (0.6, -11.6) rectangle (2.2, -11.1);  % large rectangle
		\draw[thick, fill=blue] (2.3, -11.6) rectangle (2.8, -11.1);  % small rectangle on right
		
		\draw[thick] (0, -12.2) rectangle (0.5, -11.7);  % small rectangle on left
		\draw[thick, fill=blue] (0.4, -12.2) rectangle (0.5, -11.7);  % small rectangle on left
		\draw[thick, fill=blue] (0.6, -12.2) rectangle (2.2, -11.7);  % large rectangle
		\draw[thick, fill=blue] (2.3, -12.2) rectangle (2.8, -11.7);  % small rectangle on right

		\node at (1.45, -12.45) {$\cdots\cdots$};
		
		\draw[thick] (0, -13.2) rectangle (0.5, -12.7);  % small rectangle on left
		\draw[thick] (0.6, -13.2) rectangle (2.2, -12.7);  % large rectangle
		\draw[thick] (2.3, -13.2) rectangle (2.8, -12.7);  % small rectangle on right

		%\draw[->, thick] (15, -7.6) -- (15.8, -7.6);
		%\node at (15.4, -7.85) {$\nabla f_0$};
		%\node at (15.4, -7.4) {\small $\vA\vA\zz(\cdot)$};
		%\node at (15.4, -8.3) {\small $\prox_{\eta {g}}$};
		
		\draw[->, thick] (3, -11.6) -- (3.8, -11.6);
		\node at (3.4, -11.85) {$\nabla f_0$};
		\node at (3.4, -11.4) {\small $\vA\vA\zz(\cdot)$};
		\node at (3.4, -12.2) {\small $\prox_{\eta {g}}$};
		
		\draw[thick, fill=blue] (4, -11) rectangle (4.5, -10.5);  % filled black box
		\draw[thick] (4.6, -11) rectangle (6.2, -10.5);  % large rectangle
		\draw[thick, fill=blue] (4.6, -11) rectangle (6.2, -10.5);  % part of the large rectangle filled black
		\draw[thick, fill=blue] (6.3, -11) rectangle (6.8, -10.5);  % small rectangle on right
		
		% Second row of boxes
		\draw[thick, fill=blue] (4, -11.6) rectangle (4.5, -11.1);  % small rectangle on left
		%\draw[thick, fill=blue] (4.4, -11.6) rectangle (4.5, -11.1);  % part of the large rectangle filled black
		\draw[thick, fill=blue] (4.6, -11.6) rectangle (6.2, -11.1);  % large rectangle
		\draw[thick, fill=blue] (6.3, -11.6) rectangle (6.8, -11.1);  % small rectangle on right
		
		\draw[thick] (4, -12.2) rectangle (4.5, -11.7);  % small rectangle on left
		\draw[thick, fill=blue] (4.3, -12.2) rectangle (4.5, -11.7);  % small rectangle on left
		\draw[thick, fill=blue] (4.6, -12.2) rectangle (6.2, -11.7);  % large rectangle
		\draw[thick, fill=blue] (6.3, -12.2) rectangle (6.8, -11.7);  % small rectangle on right
		
		\draw[thick ] (4, -12.8) rectangle (4.5, -12.3);  % small rectangle on left
		\draw[thick , fill=blue] (4.4, -12.8) rectangle (4.5, -12.3);  % small rectangle on left
		\draw[thick ] (4.6, -12.8) rectangle (6.2, -12.3);  % large rectangle
		\draw[thick] (6.3, -12.8) rectangle (6.8, -12.3);  % small rectangle on right

		\node at (5.45, -13.05) {$\cdots\cdots$};
		
		\draw[thick] (4, -13.7) rectangle (4.5, -13.2);  % small rectangle on left
		\draw[thick] (4.6, -13.7) rectangle (6.2, -13.2);  % large rectangle
		\draw[thick] (6.3, -13.7) rectangle (6.8, -13.2);  % small rectangle on right
		
		\draw[->, thick] (7, -11.6) -- (7.8, -11.6);
		\node at (7.4, -11.4) {\small $\vA\vA\zz(\cdot)$};
		\node at (7.4, -11.85) {\small $\prox_{\eta {g}}$};
		
		\node at (9.5, -11.6) {$\cdots\cdots\cdots$};
		
	\end{tikzpicture}
	\caption{\textbf{Illustration of the zero-respecting sequences.} Each subfigure represents one whole vector $\vx$ in a matrix format, with the first column corresponding to matrix \([\vx_1, \vx_2, \ldots, \vx_{m/3}]\), the second column to \([\vx_{m/3+1}, \ldots, \vx_{2m/3}]\), the last column to \([\vx_{2m/3+1}, \ldots, \vx_m]\), and the $i$-th row representing the row vector $[\vx_1]_i, [\vx_2]_i,\ldots,[\vx_m]_i$. A cell is plotted white if all its elements are zero and otherwise in blue. %if any element is non-zero.
		By Lemmas~\ref{lem:nablaf} and \ref{lem:kktvio}, if any row is zero, %does not include a blue cell, 
		then $\vx$ cannot be an $\epsilon$-stationary point of instance $\mathcal{P}$. 
		Starting from \( \vx = \vzero \), the figure shows how the zero elements are changed to non-zero by using the oracle information. %i.e., a white figure, we then show how the elements become blue. 
		After the first iteration, all elements in the first row can be made non-zero according to Lemma~\ref{lem:iterateguess}(1). 
		Next, 
		all elements in the  {first}   column of the second row can be made non-zero. As the iteration proceeds, in the second column of the second row, the element {\([\vx_{m/3+1}]_2\)}  is the first to become non-zero by the operator \(\vA\vA\zz(\cdot)\).   Then, the operator \(\text{prox}_{\eta g}\) makes the next element \([\vx_{m/3+2}]_2\) non-zero, followed by the operator \(\vA\vA\zz(\cdot)\), which changes \([\vx_{m/3+3}]_2\) to non-zero, and so on. It needs at least \( m/6 \) iterations (i.e., oracles) to make the entire second column of the second row non-zero.
	Then under the action of the operators \(\text{prox}_{\eta g}\) and \(\vA\vA\zz(\cdot)\), the element \([\vx_{2m/3+1}]_2\) can turn to non-zero; by Lemma \ref{lem:iterateguess}(3), \([\vx_{2m/3+1}]_3\) turns to non-zero by using $\nabla f_0$; this process continues. 
	}\label{fig:1}
\end{figure}

Now we are ready to show an important result, in Proposition~\ref{thm:iterateguess}, on how fast $\textnormal{supp}(\bar{\vx}^{(t)})$ can expand with~$t$.
%Using the lemma above, we can 
%show that $\textnormal{supp}(\bar{\vx}^{(t)})\subset\{1,\dots,\bar{j}-1\}$ for some $\bar{j}\in\{1,2,\ldots,{\bar{d}}\}$ if $t$ is not large enough. 
{By Lemmas~\ref{lem:iterateguess} and \ref{lem:supp}, if $\vx^{(0)} = \vzero$, then $\textnormal{supp}(\bar{\vx}^{(t)})$ will grow slowly. Roughly speaking, Proposition~\ref{thm:iterateguess} shows that turning a zero in $\bar{\vx}$ to non-zero will take at least $\frac{m}{6}$ oracles after the second iteration. Figure~\ref{fig:1} illustrates the process of how zero elements turn to non-zero along the iterations by using the first-order oracle information. %Next, we show the zero-respecting sequences~\cite{carmon2020lower} in , which indicates 
It explains the core idea of the proof of Proposition~\ref{thm:iterateguess}.}

\begin{proposition}
	\label{thm:iterateguess}
Suppose an FOM in Algorithm Class~\ref{algo:class-1} is applied to instance~$\mathcal{P}$ from $\vx^{(0)}=\mathbf{0}$ and generates a sequence $\{\vx^{(t)}\}_{t\geq0}$. By notations in~\eqref{eq:t-th-iter}, %Let $\vx^{(t)}=(\vx^{(t)\top}_1, \ldots,\vx^{(t)\top}_m)\zz$ with $\vx^{(t)}_i\in\RR^{{\bar{d}}}$. , 
it holds for any $\bar{j}\in\{2,3,\dots,\bar{d}\}$ that 
\begin{eqnarray}
	\label{eq:outerinduction}
	\textnormal{supp}(\vx_i^{(t)})\subset\{1,\dots,\bar{j}-1\} \text{ for }i=1,\dots,m\text{ and }t\leq 1+ m(\bar{j}-2)/6.
\end{eqnarray}
\end{proposition}
 
\proof{Proof of Proposition \ref{thm:iterateguess}.}
We prove the claim by induction on $\bar{j}$. %Before  that, we derived some properties of the operations allowed by Assumption~\ref{ass:linearspan}. 
%\YX{This proof is too long. Let's state these properties separately in a lemma or a paragraph. Then use it to prove this proposition.}
%
Let $\vxi^{(t)}=(\vxi^{(t)\top}_1, \ldots,\vxi^{(t)\top}_m)\zz$ with $\vxi^{(t)}_i\in~\RR^{{\bar{d}}}$ and $\vzeta^{(t)}=\prox_{\eta_t g}(\vxi^{(t)})=(\vzeta^{(t)\top}_1, \ldots,\vzeta^{(t)\top}_m)\zz$  with $\vzeta^{(t)}_i\in\RR^{{\bar{d}}}, \forall\, i$ for $t\geq1$. 
Since $\vx^{(0)}=\mathbf{0}$, we have $\textnormal{supp}(\nabla f_i(\vx^{(0)}_i))\subset\{1\}, \forall\, i$ from Lemma~\ref{lem:iterateguess}. Notice $\vb = \vzero$. Hence, %By Assumption~\ref{ass:linearspan},  
$\textnormal{supp}(\vxi_i^{(1)})\subset\{1\}, \forall\, i$, which further indicates $\textnormal{supp}(\vzeta_i^{(1)})\subset\{1\}, \forall\, i$ by Lemma~\ref{lem:supp}(b), and thus $\textnormal{supp}(\vx_i^{(1)})\subset\{1\}, \forall\, i$. This proves the claim in~\eqref{eq:outerinduction} for $\bar{j}=2$.
%
%Since it is shown above that  $\vx^{(0)}=\mathbf{0}$ and $\textnormal{supp}(\vx_i^{(1)})\subset\{1\}$ for any $i$, claim~\eqref{eq:outerinduction} holds for $\bar{j}=2$. 
Now suppose that the claim~\eqref{eq:outerinduction} holds for some $\bar{j}\geq 2$. We go to prove it for $\bar{j}+1$.
	
%Consider iteration $t$ with $t\leq 1+ m(\bar{j}-2)/6$. 
According to the hypothesis of the induction, we have 
\begin{equation}\label{eq:induction-result}
\textnormal{supp}(\vx_i^{(r)})\subset\{1,\dots,\bar{j}-1\},  \forall\, i\in [1, m] \text{ and } \forall\, r\leq \bar t:=1+ m(\bar{j}-2)/6.
\end{equation} 
Below we let $\widehat\vx^{(r)}=\vA\zz \vA\vx^{(r)}$ for any $r\geq0$ and consider two cases: $\bar{j}$ is even and $\bar{j}$ is odd. 
 	
\textbf{Case 1}: Suppose $\bar{j}$ is even.  We claim that, for $s=0,1,\dots,\frac{m}{6}$,
\begin{eqnarray}
	\label{eq:innerinduction1}
\textnormal{supp}(\vx_i^{(r)})\subset\left\{
\begin{array}{ll}
	\{1,\dots,\bar{j}\}, &\text{ if } i\in\left[1, \frac{m}{3}+2s\right],\\[0.1cm]
	\{1,\dots,\bar{j}-1\}, & \text{ if } i\in\left[\frac{m}{3}+2s+1, m\right],
\end{array}
\right.\ \forall\, r \le \bar t+s.
\end{eqnarray}
%for any $r\leq t+s$. 
%Since $\textnormal{supp}(\vx_i^{(r)})\subset\{1,\dots,\bar{j}-1\}$ for any $i$ and $r\leq t$, 
{We prove this claim by induction. First,} notice~\eqref{eq:induction-result} implies~\eqref{eq:innerinduction1} for $s=0$. Second, we suppose~\eqref{eq:innerinduction1} holds for some integer $s \in [0, \frac{m}{6})$ { and go to show that \eqref{eq:innerinduction1} also holds for $s+1$}. By Lemma~\ref{lem:iterateguess} and $\frac{m}{3} + 2s \le \frac{2m}{3}$, it holds 
$$
\textnormal{supp}(\nabla f_i(\vx_i^{(r)}))\subset\left\{
\begin{array}{ll}
	\{1,\dots,\bar{j}\}, & \text{ if } i\in\left[1, \frac{m}{3}+2s\right],\\[0.1cm]
	\{1,\dots,\bar{j}-1\}, &\text{ if }  i\in\left[\frac{m}{3}+2s+1, m\right],
\end{array}
\right.\  \forall\, r \le \bar t+s.
$$
%for any $r\leq t+s$. 
{By our induction hypothesis, the claim in \eqref{eq:innerinduction1} holds for some $s \in [0, \frac{m}{6})$, thus we obtain from   Lemma~\ref{lem:supp}(a) that}
$$
\textnormal{supp}(\widehat\vx_i^{(r)})\subset\left\{
\begin{array}{ll}
	\{1,\dots,\bar{j}\}, &\text{ if } i\in\left[1, \frac{m}{3}+2s+1\right],\\[0.1cm]
	\{1,\dots,\bar{j}-1\}, &\text{ if } i\in\left[\frac{m}{3}+2s+2, m\right],
\end{array}
\right.\ \forall\, r \le \bar t+s.
$$
 %for any $r\leq t+s$.  
 Hence, by Algorithm Class~\ref{algo:class-1}, we have
$$
\textnormal{supp}(\vxi_i^{(\bar t+s+1)})\subset\left\{
\begin{array}{ll}
	\{1,\dots,\bar{j}\}, & \text{ if } i\in\left[1, \frac{m}{3}+2s+1\right],\\[0.1cm]
	\{1,\dots,\bar{j}-1\},& \text{ if } i\in\left[\frac{m}{3}+2s+2, m\right],
\end{array}
\right.
$$
and thus it follows from Lemma~\ref{lem:supp}(b) and $\vzeta^{(t)}=\prox_{\eta_t g}(\vxi^{(t)}), \forall\, t\ge1$ that
$$
\textnormal{supp}(\vzeta_i^{(\bar t+s+1)})\subset\left\{
\begin{array}{ll}
	\{1,\dots,\bar{j}\},& \text{ if } i\in\left[1, \frac{m}{3}+2s+2\right],\\[0.1cm]
	\{1,\dots,\bar{j}-1\}, &\text{ if } i\in\left[\frac{m}{3}+2s+3, m\right].
\end{array}
\right.
$$

Now since $\vx_i^{(\bar t+s+1)} \in \mathbf{span}\left(\{\vxi_i^{(\bar t+s+1)}, \vzeta_i^{(\bar t+s+1)}\}\right)$ by Algorithm Class~\ref{algo:class-1}, we have
$$
\textnormal{supp}(\vx_i^{(\bar t+s+1)})\subset\left\{
\begin{array}{ll}
	\{1,\dots,\bar{j}\}, & \text{ if } i\in\left[1, \frac{m}{3}+2s+2\right],\\[0.1cm]
	\{1,\dots,\bar{j}-1\}, & \text{ if } i\in\left[\frac{m}{3}+2s+3, m\right],
\end{array}
\right.
$$
which means~\eqref{eq:innerinduction1} holds for $s+1$ as well. {Therefore, by induction,} \eqref{eq:innerinduction1} holds for all $s=0,1,\dots,\frac{m}{6}$. Let $s=\frac{m}{6}$ in~\eqref{eq:innerinduction1}. We have $\textnormal{supp}(\vx_i^{(r)})\subset\{1,\dots,\bar{j}\}$ for any $i$ and $r\leq \bar t+\frac{m}{6} = 1+ m(\bar{j}-2)/6+\frac{m}{6}=1+m(\bar{j}-1)/6$. 

\textbf{Case 2}: Suppose $\bar{j}$ is odd.  We claim that, for $s=0,1,\dots,\frac{m}{6}$,
\begin{eqnarray}
	\label{eq:innerinduction2}
	\textnormal{supp}(\vx_i^{(r)})\subset\left\{
	\begin{array}{ll}
		\{1,\dots,\bar{j}-1\},& \text{ if } i\in\left[1, \frac{2m}{3}-2s\right],\\[0.1cm]
		\{1,\dots,\bar{j}\},& \text{ if } i\in\left[\frac{2m}{3}-2s+1, m\right],
	\end{array}
	\right. \ \forall\, r \le \bar t+s.
\end{eqnarray}
%for any $r\leq t+s$. 
%Since $\textnormal{supp}(\vx_i^{(r)})\subset\{1,\dots,\bar{j}-1\}$ for any $i$ and $r\leq t$, this claim holds for $s=0$. 
%{\color{red}(YX: revise these arguments similar to what I did to Case 2).} 
Again~\eqref{eq:induction-result} implies~\eqref{eq:innerinduction2} for $s=0$. Suppose it holds for an integer $s\in [0, \frac{m}{6})$. %satisfying $0\leq s\leq\frac{m}{6}$. According to 
Then by Lemma~\ref{lem:iterateguess}, 
$$
\textnormal{supp}(\nabla f_i(\vx_i^{(r)}))\subset\left\{
\begin{array}{ll}
	\{1,\dots,\bar{j}-1\}, & \text{ if }  i\in\left[1, \frac{2m}{3}-2s\right],\\[0.1cm]
	\{1,\dots,\bar{j}\}, & \text{ if }  i\in\left[\frac{2m}{3}-2s+1, m\right],
\end{array}
\right. \ \forall\, r \le \bar t+s.
$$
%for any $r\leq t+s$. 
{By the induction hypothesis, the claim in \eqref{eq:innerinduction2} holds for some $s \in [0, \frac{m}{6})$, %all $r \le \bar t+s$, 
and thus we obtain from   Lemma~\ref{lem:supp}(a) that}
$$
\textnormal{supp}(\widehat\vx_i^{(r)})\subset\left\{
\begin{array}{ll}
	\{1,\dots,\bar{j}-1\},& \text{ if } i\in\left[1, \frac{2m}{3}-2s-1\right],\\[0.1cm]
	\{1,\dots,\bar{j}\},& \text{ if } i\in\left[\frac{2m}{3}-2s, m\right],
\end{array}
\right. \ \forall\, r \le \bar t+s.
$$
%for any $r\leq t+s$.  
Hence, by Algorithm Class~\ref{algo:class-1}, we have
$$
\textnormal{supp}(\vxi_i^{(\bar t+s+1)})\subset\left\{
\begin{array}{ll}
	\{1,\dots,\bar{j}-1\},& \text{ if } i\in\left[1, \frac{2m}{3}-2s-1\right],\\[0.1cm]
	\{1,\dots,\bar{j}\},& \text{ if } i\in\left[\frac{2m}{3}-2s, m\right],
\end{array}
\right.
$$
%In addition,  it is 
and then it follows from Lemma~\ref{lem:supp}(b) $\vzeta^{(t)}=\prox_{\eta_t g}(\vxi^{(t)}), \forall\, t\ge1$ that
$$
\textnormal{supp}(\vzeta_i^{(\bar t+s+1)})\subset\left\{
\begin{array}{ll}
	\{1,\dots,\bar{j}-1\},&\text{ if } i\in\left[1, \frac{2m}{3}-2s-2\right],\\[0.1cm]
	\{1,\dots,\bar{j}\},&\text{ if } i\in\left[\frac{2m}{3}-2s-1, m\right].
\end{array}
\right.
$$
Again since $\vx_i^{(\bar t+s+1)} \in \mathbf{span}\left(\{\vxi_i^{(\bar t+s+1)}, \vzeta_i^{(\bar t+s+1)}\}\right)$ by Algorithm Class~\ref{algo:class-1}, we have
$$
\textnormal{supp}(\vx_i^{(\bar t+s+1)})\subset\left\{
\begin{array}{ll}
	\{1,\dots,\bar{j}-1\},&\text{ if } i\in\left[1, \frac{2m}{3}-2s-2\right],\\[0.1cm]
	\{1,\dots,\bar{j}\},&\text{ if } i\in\left[\frac{2m}{3}-2s-1, m\right],
\end{array}
\right.
$$
which means~\eqref{eq:innerinduction2} holds for $s+1$ as well. By induction,~\eqref{eq:innerinduction2} holds for $s=0,1,\dots,\frac{m}{6}$. Let $s=\frac{m}{6}$ in~\eqref{eq:innerinduction2}. We have $\textnormal{supp}(\vx_i^{(r)})\subset\{1,\dots,\bar{j}\}$ for any $i$ and $r\leq \bar t+\frac{m}{6} = 1+ m(\bar{j}-2)/6+\frac{m}{6}=1+m(\bar{j}-1)/6$. 

Therefore, we have proved that~\eqref{eq:outerinduction} holds for $\bar{j}+1$, when $\bar j$ is either even or odd. By induction,~\eqref{eq:outerinduction} holds for any integer $\bar j \in [2, \bar d]$, and we complete the proof.
	\null\hfill$\square$
\endproof

					Finally, we are ready to give our main result about  the lower bound  of oracle complexity.
					
					\begin{theorem}
						\label{thm:lower}
Let $\epsilon>0$ and $L_f > 0$ be given.				 Then for any $\omega\in [0, \frac{150\pi \epsilon}{L_f})$, there exists an instance of problem~\eqref{eq:model} in Problem Class~\ref{prob:class-2}, i.e., instance~$\mathcal{P}$ in Definition~\ref{def:hardinstance}, such that any FOM in Algorithm Class~\ref{algo:class-1} needs at least 
					$\left \lceil \frac{\kappa([\bar{\vA}; \vA]) L_f \Delta_{F_0}}{36000\pi^2} \epsilon^{-2} \right\rceil$ calls to 
					$\mathrm{ORACLE}_1$ to obtain a point that is $\omega$-close to an $\epsilon$-stationary point of the instance, where $\Delta_{F_0}=F_0(\vx^{(0)})-\inf_\vx F_0(\vx)$. %as defined in~\eqref{eq:def-Delta}.
					\end{theorem}
					
					\proof{Proof of Theorem \ref{thm:lower}.}
As we discussed below~\eqref{eq:t-th-iter}, 	we assume $\vx^{(0)}=\mathbf{0}$ without loss of generality. %,  in the algorithm because, if not, we can always change variable $\vx$ in instance~$\mathcal{P}$ to $\vx-\vx^{(0)}$ and process the rest of this proof with this new instance~$\mathcal{P}$. 
%						
						%Suppose the algorithm is applied to an instance~$\mathcal{P}$ and all notations below are defined as in an instance~$\mathcal{P}$. 
						%Let $\vx^{(t)}=(\vx^{(t)\top}_1, \ldots,\vx^{(t)\top}_m)\zz$ with $\vx^{(t)}_i\in\RR^{{\bar{d}}}$ and $\bar{\vx}^{(t)}:= \frac{1}{m} \sum_{i=1}^m \vx_i^{(t)}$. 
						Thus by the notation in~\eqref{eq:t-th-iter}, Proposition~\ref{thm:iterateguess} indicates that 	$\textnormal{supp}(\vx_i^{(t)})\subset\{1,\dots,\bar{d}-1\}$ for any $i \in [1,m]$ and any $t\leq 1+ m(\bar{d}-2)/6$, which means $[\bar{\vx}^{(t)}]_{\bar{d}}=0$ if  $t\leq 1+ m(\bar{d}-2)/6$, where $\bar{\vx}^{(t)}= \frac{1}{m} \sum_{i=1}^m \vx_i^{(t)}$. 
						
On the other hand, suppose $\vx^*$ with the structure as in~\eqref{eq:t-th-iter} is an $\epsilon$-stationary point of instance~$\mathcal{P}$. Then by 	Lemma~\ref{cor:kktvio2}, it must also be an $\epsilon$-stationary point of~\eqref{eq:model3}. Hence, by  Lemmas~\ref{lem:nablaf} and~\ref{lem:kktvio}, we have $\left|[\bar\vx^*]_j\right| \ge \frac{150\pi\epsilon}{\sqrt{m} L_f}$ for all $j=1,\ldots,\bar d$, where $\bar{\vx}^*= \frac{1}{m} \sum_{i=1}^m \vx_i^*$. Therefore, by the convexity of the square function, it follows that
$$\|\vx^{(t)} - \vx^*\|^2 \ge \sum_{i=1}^m \left([\vx_i^{(t)}]_{\bar d} - [\vx_i^*]_{\bar d}\right)^2 \ge m \left([\bar\vx^{(t)}]_{\bar d} - [\bar\vx^*]_{\bar d}\right)^2 \ge m \left(\frac{150\pi\epsilon}{\sqrt{m} L_f}\right)^2 > \omega^2,$$
and thus $\vx^{(t)}$ is not $\omega$-close to $\vx^*$ if  $t\leq 1+ m(\bar{d}-2)/6$.
%						$$
%						\max\left\{\left\|\vH\vx^{(t)}\right\|, \min _{\vgamma}\left\|\nabla  f_0(\vx^{(t)})+\vH^{\top} \vgamma\right\|\right\} \geq \frac{\sqrt{m}}{2}\left\|\frac{1}{m}\sum_{i=1}^m\nabla f_i(\bar{\vx}^{(t)})\right\|>\epsilon,\ \forall\, t\leq 1+ m(\bar{d}-2)/6.
%						$$ 						
%						%if it holds that $t\leq 1+ m(\bar{d}-2)/6$. 
%						Hence, $\vx^{(t)}$ is not an $\epsilon$-stationary point of~\eqref{eq:model3} if $t\leq m({\bar{d}}-1)/6$
%						and not of instant~$\mathcal{P}$ either according to Lemma~\ref{cor:kktvio2}. 
						
						Moreover, by 				
						Lemma~\ref{cor:boundf}(a) and the fact that $g(\vx^{(0)})=0$ and $g(\vx) \ge 0, \forall\, \vx$, it holds that
						$$ {\bar{d}}\geq \frac{L_f \left( F_0(\vx^{(0)})-\inf_{\vx} F_0\left(\vx\right)\right)}{3000\pi^2} \epsilon^{-2}=\frac{L_f\Delta_{F_0}}{3000\pi^2}\epsilon^{-2}.
						$$
						In other words, in order for $\vx^{(t)}$ to be $\omega$-close to an $\epsilon$-stationary point of instant~$\mathcal{P}$, the algorithm needs at least $t=2+m({\bar{d}}-2)/6$ calls to $\mathrm{ORACLE}_1$.
We complete the proof by noticing
						\begin{equation*}
							2+m({\bar{d}}-2)/6\geq m{\bar{d}}/12\geq\frac{m L_f \Delta_{F_0}}{36000\pi^2} \epsilon^{-2}>\frac{\kappa([\bar{\vA}; \vA]) L_f \Delta_{F_0}}{36000\pi^2\epsilon^2},
						\end{equation*}
						 where the first inequality is because ${\bar{d}}\geq 5$, and the last one is by Lemma~\ref{lem:condH}.
						\null\hfill$\square$
\endproof
%\begin{remark}
%Our new lower complexity bound of  $\left \lceil \frac{\kappa([\bar{\vA}; \vA]) L_f \Delta_{F_0}}{36000\pi^2} \epsilon^{-2} \right\rceil$ is established based on an instance in Problem Class~\ref{prob:class-2}. Because Problem Class~\ref{prob:class-1} is larger than Problem Class~\ref{prob:class-2}, this lower bound is also valid for the former. However, without the structure in Problem Class~\ref{prob:class-2}, we do not have a clear form of the lower bound for Problem Class~\ref{prob:class-1}. We conjecture that $\Omega\left( \mathrm{compl}(\vA, g,\epsilon)L_f \Delta_{F_0} \epsilon^{-2}\right)$ is a lower bound of FOMs in Algorithm Class~\ref{algo:class-1} for finding an $\epsilon$-stationary solution of instances in Problem Class~\ref{prob:class-1}, where $\mathrm{compl}(\vA, g,\epsilon)$ denotes the minimum number of calls to $\mathrm{ORACLE}_1$ to find an $O(\epsilon)$-stationary point of $\min \{\frac{L_f}{2}\|\vx\|^2 + g(\vx): \vA\vx+\vb = \vzero\}$. We leave this to the future work.
%\end{remark}					
					
By Lemma~\ref{lem:condH} and the block diagonal structure of $\vA$ in \eqref{eq:AandAbar}, we have $\kappa(\vA) \le m_1$ and thus $\frac{\kappa([\bar \vA; \vA])}{\kappa(\vA)} \ge \frac{m}{4 m_1} = \frac{3 m_2}{4}$, which can be arbitrarily large as there is no restriction on $m_2$. Therefore, by choosing appropriate $m_1$ and $m_2$, we have the following claim, which answers the question that we raise in the beginning.			

\begin{claim}\label{claim2}
Let $\epsilon>0$ and $L_f > 0$ be given. For any $C > 1$, there exists an instance %, i.e., instance~$\mathcal{P}$ in Definition~\ref{def:hardinstance}, 
in Problem Class~\ref{prob:class-2}, such that  the number of calls to $\mathrm{ORACLE}_1$  by an FOM in Algorithm Class~\ref{algo:class-1} is at least $\left\lceil C \kappa(\vA) L_f \Delta_{F_0} \epsilon^{-2} \right\rceil $  to produce an $\epsilon$-stationary point of the instance. Namely, the existence of a non-smooth regularizer $g$ in problem \eqref{eq:model} can make the problem {significantly} harder.
\end{claim}	

\begin{remark}[Comparison to the lower complexity bound in Sun and Hong~\cite{sun2019distributed}]
The class of affinely constrained smooth problems considered in Sun and Hong~\cite{sun2019distributed} is a strict subclass of Problem Class~\ref{prob:class-2} that we consider. Claim~\ref{claim2} indicates that our lower complexity bound for a broader problem class can be {significantly} higher than the lower bound shown in Sun and Hong~\cite{sun2019distributed} for a smaller problem class, even if the proximal mapping of the non-smooth term is accesible. %In addition, notice that the instance $\cP$ in Definition~\ref{def:hardinstance} has the regularizer $g(\vx) = \beta \|\bar\vA \vx\|_1$. For a fixed finite $\eta>0$, performing $\prox_{\eta g}(\hat\vx)$ is equivalent to projecting $\hat\vx$ onto the set $\{\vx: \bar\vA \vx\}$ as $\beta\to \infty$. That means for an affinely constrained smooth problem, it does not become easier even if the projection onto a sub-system of the affine constraint is allowed, i.e., a stronger oracle is permitted. 
\end{remark}

	\section{Lower Bound of Oracle Complexity for Algorithm Class~\ref{algo:class-2}}
				\label{sec:extension}
				%\comm{Is the word ``sophisticated'' proper?}
In Algorithm Class~\ref{algo:class-1}, an FOM is allowed to call the operator $\prox_{\eta g}(\cdot)$ that may be nontrivial when $g$ has the structure as that in Problem Class~\ref{prob:class-2}. Calculating $\prox_{\eta g}(\cdot)$ to a high accuracy or exactly may require many (or even infinitely many) calls to $\bar{\vA}$ and $\bar{\vb}$. In contrast, $\bar g$ %introduced in Assumption~\ref{assume:problemsetup} 
is simpler than $g$, making $\prox_{\eta \bar g}(\cdot)$ easier to compute such as when $\bar g(\cdot) = \lambda\|\cdot\|_1$ for some $\lambda >0$ as in instance $\cP$. %compared to $\prox_{\eta g}(\cdot)$.  
	These observations motivate us to reformulate \eqref{eq:model} into~\eqref{eq:model-spli} and seek an $\epsilon$-stationary point of \eqref{eq:model-spli} by an FOM in Algorithm Class~\ref{algo:class-2}. %which only requires the computation of $\prox_{\eta \bar g}(\cdot)$.
%Furthermore, the approach that is designed based on \eqref{eq:model-spli} exhibits the potential to find a near $\epsilon$-stationary point of problem~\eqref{eq:model}, as shown in the next section. 

Consequently, two intriguing questions arise: (i) whether finding an $\epsilon$-stationary point of \eqref{eq:model-spli} by an FOM in Algorithm Class~\ref{algo:class-2} is easier or more challenging compared to finding an $\epsilon$-stationary point of its corresponding problem \eqref{eq:model} in Problem Class~\ref{prob:class-2} by an FOM in Algorithm Class~\ref{algo:class-1}, and (ii) whether finding a near $\epsilon$-stationary point of problem \eqref{eq:model} in Problem Class~\ref{prob:class-2} by an FOM in Algorithm Class~\ref{algo:class-2} is easier or harder  compared to that  %finding a nearly $\epsilon$-stationary point using an algorithm that satisfies 
by an FOM in Algorithm Class~\ref{algo:class-1}. We provide an answer to the first question, by showing that the same-order lower bound of complexity in Theorem~\ref{thm:lower} holds for finding  an $\epsilon$-stationary point of~\eqref{eq:model-spli} %by an algorithm that satisfies 
by an FOM in Algorithm Class~\ref{algo:class-2}. Moreover, we provide an answer to the second question, by showing that the lower bound of oracle complexity for finding a near $\epsilon$-stationary point of problem \eqref{eq:model} in Problem Class~\ref{prob:class-2} is $O({\kappa([\bar{\vA}; \vA]) L_f \Delta_{F_0}} \epsilon^{-2})$ by an FOM in either Algorithm Class~\ref{algo:class-1} or Algorithm Class~\ref{algo:class-2}; see Theorem \ref{thm:lower} and Corollary \ref{cor:lower2}.

Before we give our main results in this section, we introduce an instance of~\eqref{eq:model-spli} that is a reformulation of instant~$\mathcal{P}$ in Definition~\ref{def:hardinstance}. 	We show that an $\epsilon$-stationary point of the reformulation~\eqref{eq:model-spli} of instant~$\mathcal{P}$ is a $2\epsilon$-stationary point of the auxiliary problem~\eqref{eq:model3} of instant~$\mathcal{P}$. This is formally stated in the lemma below. 

%{\color{blue}Before we give our main results in this section, we introduce an instance of~\eqref{eq:model-spli} that is a reformulation of instant~$\mathcal{P}$ in Definition~\ref{def:hardinstance}. 
%\begin{definition}[instance $\cQ$]
%	\label{def:hardinstance_SP}
%We refer to as \emph{instance $\cQ$} the instance of problem~\eqref{eq:model-spli}  where $f_0$, $\bar g$, and $(\vA, \bar\vA, \vb, \bar\vb)$ are from instance $\cP$.
%\end{definition}	
%	We  show that an $\epsilon$-stationary point of instant $\mathcal{Q}$ is a $2\epsilon$-stationary point of the auxiliary problem~\eqref{eq:model3} of instant~$\mathcal{P}$. This is formally stated in the lemma below. 
%}

%{\color{blue} The description of the following lemma may not be proper. The variables of problem~\eqref{eq:model-spli} are $\vx$ and $\vy$, while the  variables of problem~\eqref{eq:model3} is $\vx$.}
\begin{lemma}
	\label{cor:kktequiv3}
Let $\epsilon>0$ be given in Definition~\ref{def:hardinstance} for instance~$\mathcal{P}$, and let  $\widehat{\epsilon}\in[0,\epsilon]$. Suppose $(\vx^*,\vy^*)$ is an $\widehat{\epsilon}$-stationary point of the reformulation~\eqref{eq:model-spli} of instant~$\mathcal{P}$.
%{\color{blue}instant $\mathcal{Q}$}. 
Then $\vx^*$ is a $2\widehat{\epsilon}$-stationary point of the auxiliary problem~\eqref{eq:model3} of instant~$\mathcal{P}$. 
\end{lemma}
\proof{Proof of Lemma \ref{cor:kktequiv3}.}
By Definition~\ref{def:eps-pt-P}, there exist $\vz_1\in \RR^{\bar{n}}$ and $\vz_2\in\RR^n$ such that the conditions in \eqref{eq:kktviofgesub} hold.
 Hence, for some $\vxi\in \partial \bar{g}(\vy^*)$, we have  $\|\vxi - \vz_1\| \le \widehat{\epsilon}$, and
\begin{equation}
	\label{eq:dderivegeq0_changenew}
\|\nabla f_0(\vx^*) +\bar{\vA}\zz\vxi + \vA\zz\vz_2\|\leq\|\nabla f_0(\vx^*) + \bar{\vA}\zz \vz_1 + \vA\zz \vz_2\|+ \|\bar\vA\|\widehat{\epsilon}\leq \widehat{\epsilon}+ \|\bar\vA\|\widehat{\epsilon}.
\end{equation}
Moreover, \eqref{eq:kktviofgesub} implies that, for any $\vv\in\mathbf{Null}(\vA)$,  we have
\begin{equation}
	\label{eq:dderivegeq0epsilon_changenew}
	\begin{aligned}
		F^{\prime}(\vx^*,\vy^*;\vv, \bar{\vA}\vv)=& \vv\zz\nabla  f_0(\vx^*)+ \bar g^{\prime}(\vy^*;\bar\vA\vv) \\
		\geq&
		 \vv\zz\nabla  f_0(\vx^*) +\vv\zz \bar\vA\zz\vxi =\vv\zz\left(\nabla f_0(\vx^*) +\bar{\vA}\zz\vxi + \vA\zz\vz_2\right)\geq -\widehat{\epsilon}(1+\|\bar\vA\|)\|\vv\|,
	\end{aligned}
\end{equation}
where the first inequality follows from \eqref{eq:result-dir-der} with $g$ replaced by $\bar g$, %is because $\bar g^{\prime}(\vy^*;\vu)=\sup_{\vxi'\in\partial \bar g(\vx^*)}\vv\zz \bar\vA\zz \vxi'\geq \vv\zz \bar\vA\zz \vxi$ 
and the second one is by Cauchy-Schwarz inequality and~\eqref{eq:dderivegeq0_changenew}.

Below we prove $\vy^*=\mathbf{0}$. We write it into the block-structured form
$$\vy^*=(\vy_1^{*\top}, \ldots,\vy_{3m_2-1}^{*\top})\zz \text{ with }\vy^*_{i}\in\RR^{{\bar{d}}}, \forall\, i=1,2,\ldots,3m_2-1.$$ 
If  $\vy^*\neq \mathbf{0}$, then $\vy^*_{\bar i}\neq \vzero$ for some $\bar i\in \{1,2,\ldots,3m_2-1\}$. Let  
%$
%\bar{\vxi}:=\vy^*_{\bar i}\neq \mathbf{0} 
%$ 
%and 
$\vv^*=(\vv_1^{*\top},\vv_2^{*\top},\dots,\vv_m^{*\top})^\top$ where %$\vv_i^*\in\mathbb{R}^{\bar{d}}$ is defined as
$
\vv_i^*= \vy^*_{\bar i}/(mL_f) \text{ for }i\leq \bar{i}m_1,  \text{ and } \vv_i^*= \vzero \text{ otherwise.}
%\vv_i^*=\left\{
%\begin{array}{cl}
%	\bar{\vxi}/(mL_f)&\text{ if }i\leq \bar{i}m_1\\
%	\mathbf{0}	&\text{ if }i>\bar{i}m_1,
%\end{array}
%\right.
$
We then have $\vv^*_i=\vv^*_{i+1}$ for any $i\neq \bar{i}m_1$, so $\vA\vv^*=\mathbf{0}$ by Proposition~\ref{prop:consensus}(c). Moreover, let $\vu^*=\bar\vA\vv^*=(\vu_1^{*\top}, \ldots,\vu_{3m_2-1}^{*\top})\zz$ with $\vu^*_{i}\in\RR^{{\bar{d}}}$ for $i=1,2,\ldots,3m_2-1$. We must have 
$\vu_i^*= -\vy^*_{\bar i}$ for $i= \bar{i}$ and $\vu_i^*= \vzero$ for $i\neq\bar{i}$. %, 
Therefore by \eqref{eq:gbar}, for any $s\in (0,1)$, %we can show that
\begin{align}\label{eq:gtvchangey}
	%\nonumber
	 \bar g\left(\vy^*+s \vu^*\right) 
	= & \frac{\beta}{mL_f}  \sum_{i<\bar{i}}\left\|\vy^*_i\right\|_1+\frac{\beta}{mL_f} \|\vy^*_{\bar{i}}-s\vy^*_{\bar i}\|_1
	+\frac{\beta}{mL_f}\sum_{i>\bar{i}}\left\|\vy^*_i\right\|_1 
	=  \bar g\left(\vy^*\right)-\frac{s\beta}{mL_f}\|\vy^*_{\bar i}\|_1. 
\end{align}

%According to the inequality~
Now by \eqref{eq:fi_infinitynorm}, the choice of $\vv^*_i$, and the mean value theorem, we have for any $i=1,\dots,m$ that
$$
f_i(\vx^*_i+s \vv^*_i)- f_i(\vx^*_i)=s\nabla f_i(\vx^*_i+s' \vv^*_i)^\top\vv^*_i\leq s\|\nabla f_i(\vx^*_i+s' \vv^*_i)\|_\infty\|\vv^*_i\|_1\leq \frac{50s\pi{\epsilon}}{\sqrt{m}}\cdot\frac{\|\vy^*_{\bar i}\|_1}{mL_f},
$$
where $s'\in(0,s)$. The inequality above, together with \eqref{eq:gtvchangey} and the definition of $f_0$ in~\eqref{eq:f0}, implies
\begin{equation*}
	\begin{aligned}
		&\frac{1}{s}\bigg(F(\vx^*+s \vv^*,\vy^*+s \vu^*)-F(\vx^*,\vy^*)\bigg)=\frac{1}{s}\bigg( f_0(\vx^*+s \vv^*)- f_0(\vx^*)+\bar g(\vy^*+s \vu^*)-\bar g(\vy^*)\bigg)\\
		=&\frac{1}{s}\bigg( \sum_{i=1}^m\big(f_i(\vx_i^*+s \vv_i^*)- f_i(\vx_i^*)\big)+\bar g(\vy^*+s \vu^*)-\bar g(\vy^*)\bigg)\\
		\leq& \frac{1}{smL_f}\left(50\pi s{\epsilon}\sqrt{m}\|\vy^*_{\bar i}\|_1-\beta s\|\vy^*_{\bar i}\|_1\right)= \frac{\left(50\pi {\epsilon}\sqrt{m}-\beta \right)\|\vy^*_{\bar i}\|_1}{mL_f}.
	\end{aligned}
\end{equation*}
%Taking the limit of the left-hand side of the inequality above as $t$ approaching zero, 
Letting $s\downarrow 0$ in the inequality above gives $\frac{\left(50\pi {\epsilon}\sqrt{m}-\beta \right)\|\vy^*_{\bar i}\|_1}{mL_f}\geq F^{\prime}(\vx^*,\vy^*;\vv^*,\vu^*)$. Recall $\vA\vv^*=\mathbf{0}$ and $\vu^*=\bar\vA\vv^*$. Hence, we have from \eqref{eq:dderivegeq0epsilon_changenew} and the choice of $\vv^*$ that 
$$
\frac{\left(50\pi {\epsilon}\sqrt{m}-\beta \right)\|\vy^*_{\bar i}\|_1}{mL_f}\geq -\widehat{\epsilon}(1+\|\bar\vA\|)\|\vv^*\|\geq -\widehat{\epsilon}(1+\|\bar\vA\|)\frac{\sqrt{\bar{i}}\|\vy^*_{\bar i}\|}{mL_f}
\geq -\widehat{\epsilon}(1+\|\bar\vA\|)\frac{\sqrt{m}\|\vy^*_{\bar i}\|_1}{mL_f}.
$$
%where the second inequality is obtained by setting $\vv=\vv^*$ in relation~\eqref{eq:dderivegeq0epsilon_changenew}.
Since $\beta > (50\pi+1+\|\vA\|)\sqrt{m}{\epsilon}$, $\epsilon\ge \widehat{\epsilon}$ and $\|\vA\| \ge \|\bar \vA\|$, the inequalities above can hold only when $\vy^*_{\bar i} = \vzero$. This contradicts to the hypothesis $\vy^*_{\bar i}\neq \vzero$. Hence, $\vy^* = \vzero$, %leads to a contradiction as $\beta > (50\pi+1+\|\vA\|)\sqrt{m}\epsilon$.  Hence, it holds that $\vy^*=\mathbf{0}$ for $i\in\mathcal{M}$, 
which together with \eqref{eq:kktviofgesub} gives $\|\bar\vA\vx^*\|\leq \widehat{\epsilon}$. Furthermore, $\|\vA\vx^*\|\leq \widehat{\epsilon}$ from~\eqref{eq:kktviofgesub}. Thus %Additionally, by the last inequality in~\eqref{eq:kktviofgesub}, we have 
$\|\vH\vx^*\|\leq \|\bar\vA\vx^*\|+\|\vA\vx^*\|\leq 2\widehat{\epsilon}$, which, together with  $\|\nabla f_0(\vx^*) + \bar{\vA}\zz \vz_1 + \vA\zz \vz_2\| \le \widehat{\epsilon}$ from~\eqref{eq:kktviofgesub}, indicates that $\vx^*$ a $2\widehat{\epsilon}$-stationary point of~\eqref{eq:model3}.
	\null\hfill$\square$
\endproof

%\YX{We may give a complete proof in arXiv and refer to it. That would be more convincing and also necessary to show the dependence on constants. Or can we argue that the sequence under Assumption 3 satisfies Assumption 2 for our worst-case instance?} 
With Lemma~\ref{cor:kktequiv3}, we can establish the lower-bound complexity for solving problem~\eqref{eq:model-spli} in a similar way to show Theorem~\ref{thm:lower}. In particular, Lemma~\ref{cor:kktequiv3} indicates that a solution $(\vx^{(t)},\vy^{(t)})$ cannot be an $\epsilon/2$-stationary point of %{\color{blue} instance $\mathcal{Q}$}, %~\eqref{eq:model-spli}, 
the reformulation~\eqref{eq:model-spli} of instant~$\mathcal{P}$, 
if $\vx^{(t)}$ is not an $\epsilon$-stationary point of~\eqref{eq:model3}. By Lemmas~\ref{lem:nablaf} and~\ref{lem:kktvio}, if there exists $\bar{j}\in\{1,2,\ldots,{\bar{d}}\}$ such that $[\bar{\vx}^{(t)}]_{\bar{j}}=0$, where $\bar{\vx}^{(t)}= \frac{1}{m} \sum_{i=1}^m \vx_i^{(t)}$, then $\vx^{(t)}$ cannot be an $\epsilon$-stationary point of problem~\eqref{eq:model3} and thus cannot be an $\epsilon/2$-stationary point of the reformulation~\eqref{eq:model-spli} of instant~$\mathcal{P}$.
%{\color{blue} instance $\mathcal{Q}$}.
%problem~\eqref{eq:model-spli}. 
Finally, similar to Proposition~\ref{thm:iterateguess}, we can show that, for any FOM in Algorithm Class~\ref{algo:class-2} that is applied to 
%{\color{blue} instance $\mathcal{Q}$},
the reformulation~\eqref{eq:model-spli} of instance~$\mathcal{P}$, 
if it starts from $(\vx^{(0)},\vy^{(0)})=(\mathbf{0},\mathbf{0})$ and generates a sequence $\{({\vx}^{(t)}, {\vy}^{(t)} )\}_{t\ge 1}$, then $\textnormal{supp}({\vx}^{(t)})\subset\{1,\dots,\bar{j}-1\}$ for some $\bar{j}\in\{1,2,\ldots,{\bar{d}}\}$ if $t$ is not large enough. Hence, ${\vx}^{(t)}$ cannot be $\omega$-close to any $\epsilon$-stationary point of instance~$\mathcal{P}$ for any $\omega\in [0, \frac{150\pi \epsilon}{L_f})$, according to the proof of Theorem \ref{thm:lower}. This way, we can obtain a lower bound of oracle complexity to produce an $\epsilon$-stationary point of \eqref{eq:model-spli} and also a lower bound to obtain a near $\epsilon$-stationary point of \eqref{eq:model}. %Note that, we assume $(\vx^{(0)},\vy^{(0)})=(\mathbf{0},\mathbf{0})$ without loss of generality because, otherwise, we can change variable $\vx$ and $\vy$ in  problem~\eqref{eq:model-spli} to $\vx-\vx^{(0)}$ and $\vy-\vy^{(0)}$, respectively, and process the rest of the proof with this new instance.~$\mathcal{P}$. 
Since the aforementioned arguments are similar to those for proving Theorem~\ref{thm:lower}, we simply present the lower complexity bounds in the theorem and the corollary below and put the proofs in Appendix~\ref{sec:appen1}.  

%\YX{Assumption~\ref{assume:problemsetup} is not for problem~\eqref{eq:model-spli}. }

\begin{theorem}
						\label{thm:lbcomposite}
Let $\epsilon>0$ and $L_f > 0$ be given.								%Suppose an algorithm  is applied to problem~\eqref{eq:model-spli} that satisfies Assumption~\ref{assume:problemsetup}\textnormal{(a, c)} and $\inf_{\vx, \vy} F(\vx, \vy) > -\infty$.
						%Suppose the generated sequence $\{(\vx^{(t)},\vy^{(t)})\}_{t\geq0}$ satisfies Assumption~\ref{ass:linearspan3}. 
						Then there exists an instance in Problem Class~\ref{prob:class-2}, i.e., instance~$\mathcal{P}$ in Definition~\ref{def:hardinstance}, such that any FOM in Algorithm Class~\ref{algo:class-2} requires at least $\left\lceil \frac{\kappa([\bar{\vA}; \vA]) L_f \Delta_{F_0}}{72000\pi^2} \epsilon^{-2}\right\rceil$ calls to $\mathrm{ORACLE}_2$ to obtain an $\epsilon$-stationary point of the reformulation~\eqref{eq:model-spli} of that instance, where $\Delta_{F_0}=F_0(\vx^{(0)})-\inf_\vx F_0(\vx)$. %where $\Delta_F =F(\vx^{(0)}, \vy^{(0)})-\inf_{\vx, \vy} F(\vx, \vy)$.
\end{theorem}

\begin{corollary}
	\label{cor:lower2}
	Let $\epsilon>0$ and $L_f > 0$ be given. %Suppose an algorithm is applied to problem~\eqref{eq:model} that satisfies Assumption~\ref{assume:problemsetup} 
	%and generates a sequence $\{(\vx^{(t)},\vy^{(t)})\}_{t\geq0}$ that satisfies Assumption~\ref{ass:linearspan3}. 
	Then for any $\omega\in [0, \frac{150\pi \epsilon}{L_f})$, there exists an instance of problem~\eqref{eq:model}, i.e., instance~$\mathcal{P}$ in Definition~\ref{def:hardinstance}, such that any FOM in Algorithm Class~\ref{algo:class-2} requires at least 
	$\left \lceil \frac{\kappa([\bar{\vA}; \vA]) L_f \Delta_{F_0}}{18000\pi^2} \epsilon^{-2} \right\rceil$ 
	calls to $\mathrm{ORACLE}_2$ to obtain a point $\vx^{(t)}$ that is $\omega$-close to an $\epsilon$-stationary point of instance~$\mathcal{P}$, where $\Delta_{F_0}=F_0(\vx^{(0)})-\inf_\vx F_0(\vx)$. %as defined in~\eqref{eq:def-Delta}.
\end{corollary}

\subsection*{Tightness of Lower Complexity Bounds}\label{sec:ub}
For Problem Class~\ref{prob:class-2}, we give in the extended arXiv version (see Liu et al.~\cite{liu2023first}) of this paper  an FOM that belongs to Algorithm Class~\ref{algo:class-2}. When $[\bar\vA; \vA]$ has full row rank and $\bar g$ is Lipschitz continuous, the FOM in Liu et al.~\cite{liu2023first} can find a point that is $\epsilon$-close to an $\epsilon$-stationary point of Problem Class~\ref{prob:class-2} by $\tilde O\left(\kappa([\bar{\vA}; \vA]) L_f \Delta_{F_0} \epsilon^{-2}\right)$ calls to $\mathrm{ORACLE}_2$. It can also find an $\epsilon$-stationary point of the reformulation~\eqref{eq:model-spli} by the same-order complexity. Notice that the instance $\cP$ we use to establish our lower complexity bounds satisfies the full row-rankness condition and the Lipschitz continuity of $\bar g$. Thus comparing the upper complexity bound in Liu et al.~\cite{liu2023first} to the lower bounds we establish in Theorem~\ref{thm:lbcomposite} and Corollary~\ref{cor:lower2}, we claim that our lower complexity bounds of Algorithm Class~\ref{algo:class-2}  are tight up to a logarithmic factor under the additional restriction that $\bar g$ is Lipschitz continuous and $[\bar\vA; \vA]$ has full row rank, {which indicates $\bar n + n \le d$}. We are unable to claim the tightness of our lower complexity bound of Algorithm Class~\ref{algo:class-1} in Theorem~\ref{thm:lower} for solving Problem Class~\ref{prob:class-2}. However, notice that $\prox_{\eta g}$ is harder than $\prox_{\eta \bar g}$ to evaluate. In this sense,  $\mathrm{ORACLE}_1$ is stronger than $\mathrm{ORACLE}_2$. Hence, we conjecture that our lower complexity bound of Algorithm Class~\ref{algo:class-1} in Theorem~\ref{thm:lower} is also nearly tight under the additional conditions mentioned above. 

Before concluding this section, we make a few more discussions on the challenges of claiming a (nearly) tight lower complexity bound for Algorithm Class~\ref{algo:class-1}. In literature,  {for a certain special $g$, some existing FOMs in Algorithm Class~\ref{algo:class-1} can achieve a complexity of $O(\epsilon^{-2})$, where the hidden key constants are the same as but in a higher order than those in our lower bound. For example, when $g$ is the indicator function of a polyhedral set $\cC$ and the exact projection onto $\cC$ can be computed, Zhang and Luo~\cite{zhang2022global} give an ALM-based FOM  for problem~\eqref{eq:model} and obtain a complexity of $O(\widehat\kappa^{2}L_f^3\Delta_{F_0} \epsilon^{-2})$, where $\widehat\kappa$ is a joint condition number of the equality $\vA \vx + \textbf{b} =\vzero$  and $g=0$.}  When $g$ is a general convex function, the best-known complexity bound of FOMs in Algorithm Class~\ref{algo:class-1} is $O(\epsilon^{-5/2})$ to find an $\epsilon$-KKT point of problem~\eqref{eq:model}, e.g., in Lin et al.~\cite{lin2022complexity}, Li et al.~\cite{li2021augmented}, and Kong et al.~\cite{kong2023iteration}.  The dependence on $\epsilon$ is in a higher order than our lower bound.  {In addition, these upper bound results have a higher-order dependence on $\Delta_{F_0}$. For example, the result in Lin et al.~\cite{lin2022complexity} has a dependence of $\Delta_{F_0}^{5/4}$. However, some other key constants concealed in $O(\epsilon^{-5/2})$ are not revealed in our lower complexity bound, such as a uniform upper bound on the dual variables in Lin et al.~\cite{lin2022complexity}.} This can be one main reason for the discrepancy between the $O(\epsilon^{-5/2})$ upper complexity bound and our established lower bound.  {Though we may express our lower bound in a different formulation to reveal the dependence on the key constants in these existing upper bound results, we prefer our formulation because it more clearly shows the difficulty caused by the nonsmooth convex regularizer $g$, while the existing $O(\epsilon^{-5/2})$ results remain in the same order even if $g\equiv 0$.} It is possible to have an upper bound that has a lower-order than $O(\epsilon^{-5/2})$, if the upper bound reveals a dependence on a certain ``condition number'' that quantifies the difficulty caused by the interaction of $g$ and the linear constraint $\vA\vx+\vb=\vzero$,  {e.g., the difficulty of solving $\min_\vx \left\{\frac{1}{2}\|\vx-\hat\vx\|^2+ \alpha\cdot g(\vx), \st \vA\vx+\vb=\vzero\right\}$ to a high accuracy for any $\hat\vx$ and $\alpha>0$, by using $\mathrm{ORACLE}_1$ in \eqref{oracle-1}}. However, it is unclear how to define such a ``condition number'' for a general convex function~$g$.  
 	Our results, though not providing a complete solution here, offer insights into the complexity of FOMs under the imposed conditions, thereby enlightening the field and paving the way for future research.

\section{Conclusion and Open Questions}\label{sec:conclusion}
We have studied lower bounds of the oracle complexity of first-order methods for finding a (near) $\epsilon$-stationary point of a non-convex composite non-smooth optimization problem with affine equality constraints. We have shown that the non-smooth term can cause the problem {significantly harder} even if its proximal mapping is accessible. This is fundamentally different from composite non-smooth optimization without affine constraints. When the non-smooth term has a special structure, we give an explicit formula of the lower bound, and we show that the same-order lower bound holds for first-order methods applied to a splitting reformulation of the original problem. %using auxiliary variables. 
In addition, by relating our lower bounds to an upper complexity bound result in the extended arXiv version see (Liu et al.~\cite{liu2023first}) of this paper, we show that our lower bounds are tight, up to a logarithmic factor, when the non-smooth term has a special structure and a certain full row-rankness holds. %design an inexact proximal gradient method for the reformulation. To find an $\epsilon$-stationary point (that is also a near $\epsilon$-stationary point of the original problem), our designed method has an oracle complexity that matches the established lower bounds under two scenarios, . This shows that both the lower bound and the oracle complexity of the  inexact proximal gradient method are nearly not improvable. 

%\comm{Shall we mention inequality linear constrained non-convex non-smooth problems a little bit here?}

Though we make the first attempt on establishing lower complexity bounds of first-order methods for solving affine-constrained non-convex composite non-smooth problems, there are still a few open questions that are worth further exploration. %Given the results in this paper, 
First, for Algorithm Class~\ref{algo:class-1} on solving Problem Class~\ref{prob:class-2}, can the lower bound $O({\kappa([\bar{\vA};\vA]) L_f \Delta_{F_0}} \epsilon^{-2})$ be achieved, i.e., is it tight? The second question is how to formulate a lower complexity bound of Algorithm Class~\ref{algo:class-1} on solving problem~\eqref{eq:model} if it does have the structure of $g(\vx) = \bar g(\bar\vA \vx + \bar\vb)$ as in Problem Class~\ref{prob:class-2}. %whether the lower bound $O({\kappa([\bar{\vA};\vA]) L_f \Delta_{F}} \epsilon^{-2})$ for problem~\eqref{eq:model-spli} can be achieved by an algorithm satisfying Assumption~\ref{ass:linearspan3} without Assumptions~\ref{ass:222},~\ref{ass:dual2} or~\ref{ass:polyhedralg}. 
Third, what will the lower bound look like if there are convex nonlinear inequality constraints? %Another open question is either the lower bound or the iteration complexity~\eqref{eq:complexity_polyhedral} can be improved under Assumptions~\ref{ass:222} and~\ref{ass:polyhedralg} when $\bar{\vb}\neq\mathbf{0}$ or $\vb\neq\mathbf{0}$. \YX{How about remove the 3rd question?}		

\begin{APPENDICES}
	\section{Proofs of Some Lemmas}\label{sec:proof-lem}
In this section, we provide the proofs of some lemmas that we use in Section~\ref{sec:lb1}.

\proof{Proof of Lemma~\ref{lem:nablaf}.}	Let $\bar{\vz}=\frac{\sqrt{m} L_f \vz}{150\pi\epsilon}$. We first consider the  case where  $|[\vz]_j|<\frac{150\pi\epsilon}{\sqrt{m} L_f}$ for all $j=1,2,\ldots,\bar{j}$. In this case,  $| [\bar \vz]_1|=\frac{\sqrt{m} L_f |[\vz]_1|}{150\pi\epsilon} <1$. By \eqref{eq:fi}, we have 
	\begin{equation}
		\label{eq:nablaf1}
		\begin{aligned}
			\left\|\frac{1}{m}\sum_{i=1}^m\nabla f_i(\vz)\right\|  \geq \left| \frac{1}{m}\sum_{i=1}^m\left[\nabla f_i(\vz)\right]_1\right| = \left| \frac{2\epsilon}{m\sqrt{m}}\sum_{i=1}^m \left[\nabla {h}_i\left(\bar{\vz}\right)\right]_1\right|.
		\end{aligned}
	\end{equation}
	In addition, according to \eqref{eq:deriveh3}, we have
	\begin{equation}
		\label{eq:nablaf2}
		\begin{aligned}
			\frac{1}{m}\sum_{i=1}^m \left[\nabla {h}_i\left(\bar{\vz}\right)\right]_1 &=-\Psi(1) \Phi^{\prime}([\bar{\vz}]_1)+ \left[-\Psi^{\prime}(-[\bar \vz]_1) \Phi(-[\bar \vz]_2)-\Psi^{\prime}([\bar \vz]_1) \Phi([\bar \vz]_2)\right] \leq-\Psi(1) \Phi'([\bar{\vz}]_1)< -1,
		\end{aligned}
	\end{equation}
	where the first inequality comes from the non-negativity of $\Psi^{\prime}$ and $\Phi$ by Lemma~\ref{lem:functions}(a), and the second inequality is  by  Lemma~\ref{lem:functions}(b) and $|\bar \vz_1|<1$. Combing~\eqref{eq:nablaf1} and~\eqref{eq:nablaf2}
	yields the desired inequality.
	
	Second, we consider the case where there exists $j\in\{2,\ldots,\bar{j}\}$ such that $|[\vz]_j| <\frac{150\pi\epsilon}{\sqrt{m} L_f} \leq |[\vz]_{j-1}|$. In this case,  $|[\bar \vz]_j| <1 \leq |[\bar  \vz]_{j-1}|$. By~\eqref{eq:fi} again, we have 
	\begin{equation}
		\label{eq:nablaf1new}
		\left\|\frac{1}{m}\sum_{i=1}^m\nabla f_i(\vz)\right\|  \geq \left| \frac{1}{m}\sum_{i=1}^m\left[\nabla f_i(\vz)\right]_j\right| = \left| \frac{2\epsilon}{m\sqrt{m}}\sum_{i=1}^m \left[\nabla {h}_i\left(\bar{\vz}\right)\right]_j\right|.
	\end{equation}		
	According to \eqref{eq:deriveh1} and~\eqref{eq:deriveh2}, we have
	\begin{equation}
		\label{eq:nablaf3}
		\begin{aligned}
			&\frac{1}{m}\sum_{i=1}^m \left[\nabla {h}_i\left(\bar{\vz}\right)\right]_j \\=& -\Psi(-[\bar{\vz}]_{j-1}) \Phi^{\prime}(-[\bar{\vz}]_{j})-\Psi([\bar{\vz}]_{j-1}) \Phi^{\prime}([\bar{\vz}]_{j})-\Psi^{\prime}(-[\bar{\vz}]_{j}) \Phi(-[\bar{\vz}]_{j+1})-\Psi^{\prime}([\bar{\vz}]_{j}) \Phi([\bar{\vz}]_{j+1})\\
			\leq &-\Psi(-[\bar{\vz}]_{j-1}) \Phi^{\prime}(-[\bar{\vz}]_{j})-\Psi([\bar{\vz}]_{j-1}) \Phi^{\prime}([\bar{\vz}]_{j}) = -\Psi(|[\bar{\vz}]_{j-1}|) \Phi^{\prime}([\bar{\vz}]_{j}) < -1,
		\end{aligned}
	\end{equation}
	where the first inequality comes from the nonnegativity of $\Psi^{\prime}$ and $\Phi$ by  Lemma~\ref{lem:functions}(a), the second equality holds by the fact that $\Phi^{\prime}(u)=\Phi^{\prime}(-u)$ and $\Psi(u)=0$ for all $u\leq 0$ from \eqref{eq:PsiPhi} and Lemma~\ref{lem:functions}(a), and the second inequality is  by  Lemma~\ref{lem:functions}(b) and the fact that $|[\bar{\vz}]_{j-1}|\geq 1$ and $|[\bar{\vz}]_{j}|<1$. Combining~\eqref{eq:nablaf1new} and~\eqref{eq:nablaf3} 	yields the desired inequality and completes the proof.		
	\null\hfill$\square$
\endproof

\proof{Proof of Lemma~\ref{lem:kktvio}.}
	By simple calculation and the fact $\mathbf{Null}(\vH)=\left\{ \mathbf{1}_{m} \otimes \vu:\vu\in\RR^{\bar{d}}\right\}$, we obtain  		
	\begin{align}
		\min_{\vgamma}\left\|\nabla  f_0(\vx)+\vH^{\top} \vgamma\right\|^2 
		= &\left\|\mathbf{Proj}_{\mathbf{Null}(\vH)}(\nabla  f_0(\vx))\right\|^2 \textstyle  =\frac{1}{m}\left\|\sum_{i=1}^m \nabla f_i\left(\vx_i\right)\right\|^2, \label{eq:kktvio1} \\
		\left\|\vH\vx\right\|^2=&m^2L_f^2\sum_{i=1}^{m-1} \left\|\vx_{i}-\vx_{i+1}\right\|^2. \label{eq:kktvio1-Hx}
		%&\left\|\vv^* -\mathbf{Proj}_{\vH\vx=0}(\vx^*)\right\|^2=\frac{1}{m} \sum_{i=1}^m \left\|\frac{1}{m} \sum_{j=1}^m \vx_j-\vx_i\right\|^2.
	\end{align}
	By the $L_f$-Lipschitz continuity of $\nabla  f_0$, 
	we  have %\YX{Give explanation of the last three inequalities.}
	\begin{align*}
		&\textstyle \frac{1}{2}\left\|\frac{1}{m}\sum_{i=1}^m\nabla f_i(\bar{\vx})\right\|^2 - \left\|\frac{1}{m} \sum_{i=1}^m \nabla f_i\left(\vx_i\right)\right\|^2\leq \left\|\frac{1}{m} \sum_{i=1}^m\left(\nabla f_i(\bar{\vx})-\nabla f_i\left(\vx_i\right)\right)\right\|^2 \\
		\leq &\frac{1}{m} \sum_{i=1}^m\left\|\nabla f_i\left(\bar{\vx}\right)-\nabla f_i\left(\vx_i\right)\right\|^2
		\leq \frac{1}{m} \sum_{i=1}^m L_f^2\left\|\bar{\vx}-\vx_i\right\|^2 
		%\leq &\frac{2}{m}\max\left\{\left\|\vv^* -\mathbf{Proj}_{\vH\vx=0}(\vx^*)\right\|^2, \left(\min _{\gamma}\left\|\nabla  f_0(\vx)+\vH^{\top} \gamma\right\|\right)^2\right\}.\\
		\stackrel{(a)}{\leq} \frac{L_f^2}{m^2} \sum_{i=1}^{m} \sum_{j=1}^m\left\|\vx_j-\vx_i\right\|^2\\
		\stackrel{(b)}{\leq} &\frac{L_f^2}{m^2} \sum_{i=1}^{m} \sum_{j=1}^m \left[|{j-i}|\sum_{k=\min\{i,j\}}^{\max\{i,j\}-1}\left\|\vx_k-\vx_{k+1}\right\|^2\right]\stackrel{(c)}{\leq} \frac{L_f^2}{m^2} \sum_{i=1}^{m} \sum_{j=1}^m m\sum_{k=1}^{m-1} \left\|\vx_{k}-\vx_{k+1}\right\|^2
\\=   &mL_f^2 \sum_{i=1}^{m-1} \left\|\vx_{i}-\vx_{i+1}\right\|^2, % = \frac{1}{m} \left\|\vH\vx\right\|^2 \leq \frac{2}{m}\max\left\{\left\|\vH\vx\right\|^2, \min _{\vgamma}\left\|\nabla  f_0(\vx)+\vH^{\top} \vgamma\right\|^2\right\},  %(\text{ multiplying $\vH$ by $m$})
	\end{align*}
	where (a) comes from $\|\bar{\vx}-\vx_i\|^2 = \frac{1}{m^2}\|\sum_{j=1}^m \vx_j- \vx_i\|^2\leq \frac{1}{m}  \sum_{j=1}^m \|\vx_j- \vx_i\|^2 $, (b) results from the fact that $\|\vx_j- \vx_i\|^2 \leq |j-i| \sum_{k=\min\{i,j\}}^{\max\{i,j\}-1}\left\|\vx_k-\vx_{k+1}\right\|^2$, (c) holds by $|j-i|<m$, $\max\{i,j\}\leq m$ and $\min\{i,j\}\ge 1$.
	
	%where the last inequality is from equality~\eqref{eq:kktvio1}. 
	Hence, by~\eqref{eq:kktvio1} and~\eqref{eq:kktvio1-Hx} and the fact $a+b \le 2\max\{a,  b\}$ for any $a, b\in \RR$, we obtain the desired result from the inequality above and  complete the proof. 
	\null\hfill$\square$
\endproof

\section{Proofs of Theorem~\ref{thm:lbcomposite} and Corollary \ref{cor:lower2}}		
\label{sec:appen1}
In this section, we give a complete proof of Theorem~\ref{thm:lbcomposite} and Corollary \ref{cor:lower2}.
%Before we give the proof to Theorem~\ref{thm:lbcomposite}, 
We first show a lemma and a proposition.			
According to the structure of $\vA$ and $\bar{\vA}$ given in~\eqref{eq:AandAbar},  $\textnormal{supp}((\vA^{\top}\vA\vx)_i)$ and $\textnormal{supp}((\bar{\vA}^{\top}\bar{\vA}\vx)_i)$ are determined by $\textnormal{supp}(\vx_{i-1})$, $\textnormal{supp}(\vx_{i})$ and $\textnormal{supp}(\vx_{i+1})$. Also, $\textnormal{supp}(\prox_{\eta \bar{g}}(\vy))$, $\bar{\vA}\zz\vy$ and $\vA\vx$ have a similar property. They are stated in the following lemma.

\begin{lemma}
	\label{lem:supp2}
	Let $\vx$ be the structured vector given in~\eqref{eq:xblock}, $\vA$ {and $\bar{\vA}$} in~\eqref{eq:AandAbar}, and $\bar{g}$ be given in~\eqref{eq:gbar}. Define $\vx_0=\vx_{m+1}=\mathbf{0}\in\mathbb{R}^{\bar d}$. The following statements hold:
	\begin{enumerate}
		\item[\textnormal{(a)}] Let  $\widehat\vx=(\widehat\vx_1\zz, \ldots,\widehat\vx_m\zz)\zz\in\textnormal{\textbf{span}}\{\vA^{\top}\vA\vx, \bar{\vA}^{\top}\bar{\vA}\vx\}$ with $\widehat\vx_i\in\RR^{{\bar{d}}}$.  Then %It then yields that 
		\begin{equation}\label{eq:supp-relation2}
			\textnormal{supp}(\widehat\vx_i)\subset \textnormal{supp}(\vx_{i-1})\cup\textnormal{supp}(\vx_{i})\cup\textnormal{supp}(\vx_{i+1}),\, \forall\, i\in [1, m].
		\end{equation}		
		
		\item[\textnormal{(b)}] Let  $\vy=(\vy^{\top}_1, \ldots,\vy^{\top}_{3m_2-1})\zz$ with $\vy_j\in\RR^{{\bar{d}}}$ and $\widetilde \vx=\bar\vA^\top \vy$. Then 
		\begin{eqnarray}
			\label{eq:suppAy}
			\textnormal{supp}(\widetilde\vx_i)&\subset&
			\left\{
			\begin{array}{ll}
				\emptyset&\text{ if }i-1,i\notin\mathcal{M},\\
				\textnormal{supp}(\vy_{j})&\text{ if }i-1=jm_1\in\mathcal{M},\\
				\textnormal{supp}(\vy_{j})&\text{ if }i=jm_1\in\mathcal{M},
			\end{array}
			\right.\,\,\forall i\in[1,m].
		\end{eqnarray} 
		
		\item[\textnormal{(c)}]	Let $\widehat\vy=\bar\vA\vx$ and $\widehat\vy=(\widehat\vy_1\zz, \ldots,\widehat\vy_{3m_2-1}\zz)\zz$ with $\widehat\vy_j\in\RR^{{\bar{d}}}$. Then 
		\begin{eqnarray}
			\label{eq:suppAbar3}
			\textnormal{supp}(\widehat\vy_j)
			&\subset&\textnormal{supp}(\vx_{jm_1})\cup\textnormal{supp}(\vx_{jm_1+1}),\,\,\forall j\in[1,3m_2-1].
		\end{eqnarray} 
		
		\item[\textnormal{(d)}] It holds that 
		\begin{eqnarray}
			\label{eq:suppAbarAbary}
			\textnormal{supp}(\vy)&=&\textnormal{supp}(\bar\vA\bar\vA\zz\vy),\,\,\forall \vy\in\mathbb{R}^{\bar n}.
		\end{eqnarray}
		
		\item[\textnormal{(e)}] For any given $\eta>0$, let $\widetilde\vy=\prox_{\eta \bar g}(\vy)=(\widetilde\vy_1\zz, \ldots,\widetilde\vy_{3m_2-1}\zz)\zz$ with $\widetilde\vy_j\in\RR^{{\bar{d}}}$.   Then 
		\begin{eqnarray}
			\label{eq:suppg3}
			\textnormal{supp}(\widetilde\vy_j)&\subset&\textnormal{supp}(\vy_j),\,\,\forall j\in[1,3m_2-1].
		\end{eqnarray}
		
	\end{enumerate}
\end{lemma}

\proof{Proof of Lemma \ref{lem:supp2}.}
	(a) %Note that $\vA^{\top}\vb=\mathbf{0}$ for instance~$\mathcal{P}$. Let $\vx=(\vx_1\zz, \ldots,\vx_m\zz)\zz$ with $\vx_i\in\RR^{{\bar{d}}}$. Let  $\widehat\vx=\vA^{\top}\vA\vx=(\widehat\vx_1\zz, \ldots,\widehat\vx_m\zz)\zz$ with $\widehat\vx_i\in\RR^{{\bar{d}}}$. We define $\vx_0=\vx_{m+1}=\mathbf{0}$. By the definition of $\vA$ in~\eqref{eq:AandAbar}, we have 
	Recall that $\vH$ are split into $\bar{\vA}$ and $\vA$ in rows. The relation in~\eqref{eq:supp-relation2} immediately follows from \eqref{eq:ata-struc} the observation
	\begin{equation}\label{eq:ata-struc2}
		\newcommand{\Hmat}{\left[ 
			\begin{array}{rrrrr} 
				\vI_{\bar d} & - \vI_{\bar d} & & & \\
				- \vI_{\bar d} & 2 \vI_{\bar d} & - \vI_{\bar d} & \\
				& \ddots & \ddots & \ddots & \\
				&  &	- \vI_{\bar d} & 2 \vI_{\bar d} & - \vI_{\bar d}  \\
				& & & - \vI_{\bar d} & \vI_{\bar d}
			\end{array} \right]
		}
		\vH\zz\vH = m^2 L_f^2 
		\left.
		\,\smash[b]{\underbrace{\!\Hmat\!}_{\textstyle\text{$m-1$ blocks}}}\,
		\right\}\text{$m-1$ blocks}
		\vphantom{\underbrace{\Hmat}_{\textstyle\text{$m$ blocks}}}
		~\text{  and  }~\bar{\vA}^{\top}\bar{\vA}=	\vH\zz\vH- 	\vA^{\top}\vA.
	\end{equation}
	
	%Hence,~\eqref{eq:supp-relation2} immediately follows.	
	%	\begin{eqnarray}
		%		\nonumber
		%		\textnormal{supp}(\widehat\vx_i)&\subset&
		%		\left\{
		%		\begin{array}{ll}
			%			\textnormal{supp}(\vx_{i-1})\cup\textnormal{supp}(\vx_{i})\cup\textnormal{supp}(\vx_{i+1})&\text{ if }i-1,i\notin\mathcal{M},\\
			%			\textnormal{supp}(\vx_{i})\cup\textnormal{supp}(\vx_{i+1})&\text{ if }i-1\in\mathcal{M},\\
			%			\textnormal{supp}(\vx_{i-1})\cup\textnormal{supp}(\vx_{i})&\text{ if }i\in\mathcal{M}.
			%		\end{array}
		%		\right.%\\\label{eq:suppA}
		%		%&\subset&
		%	\end{eqnarray}
	(b)  The relation in~\eqref{eq:suppAy} immediately follows from the definitions of $\bar\vA$ and $\mathcal{M}$ in~\eqref{eq:indexsetM}.
	
	(c) The relation in~\eqref{eq:suppAbar3} immediately follows from the definition of $\bar\vA$ and $\mathcal{M}$ in~\eqref{eq:indexsetM}.
	
	(d) The relation in~\eqref{eq:suppAbarAbary} immediately follows from the fact that $\bar\vA\bar\vA\zz=2m^2L_f^2\vI_{\bar{n}}$ by
	the definition of $\bar\vA$ and $\mathcal{M}$ in~\eqref{eq:indexsetM}.
	
	(e) Given any $y\in\mathbb{R}$ and any $c>0$, consider the following optimization problem in $\mathbb{R}$:
	\begin{align}
		\label{eq:2dproxg3}
		\widetilde y=\argmin_{z\in\mathbb{R}}\frac{1}{2}(z-y)^2+c|z|=\text{sign}(y)\cdot (|y|-c)_+.
	\end{align}
	Recall the definition of $\bar{g}$ in \eqref{eq:gbar}, we obtain that 
	$$
	\widetilde\vy=\prox_{\eta \bar g}(\vy) = \argmin_{\vy'} \frac{\beta}{mL_f} \|\vy'\|_1 + \frac{1}{2}\|\vy'-\vy\|^2.
	$$	
	Applying \eqref{eq:2dproxg3} to each coordinate of $\widetilde\vy$ above, we have~\eqref{eq:suppg3}
	and complete the proof. 
	\null\hfill$\square$
\endproof

Now we are ready to show the following result on how fast $\textnormal{supp}(\bar{\vx}^{(t)})$ and $\textnormal{supp}(\bar{\vy}^{(t)})$ can expand with $t$.
%Using the lemma above, we can 
%show that $\textnormal{supp}(\bar{\vx}^{(t)})\subset\{1,\dots,\bar{j}-1\}$ for some $\bar{j}\in\{1,2,\ldots,{\bar{d}}\}$ if $t$ is not large enough. 

\begin{proposition}
	\label{thm:iterateguess3}
	Suppose an FOM in Algorithm Class~\ref{algo:class-2} is applied to the reformulation \eqref{eq:model-spli} of instance~$\mathcal{P}$
	%{\color{blue}instance $\mathcal{Q}$} 
	from an initial solution $\vx^{(0)}=\mathbf{0}$ and $\vy^{(0)}=\mathbf{0}$, and generates a sequence $\{(\vx^{(t)},\vy^{(t)})\}_{t\geq1}$.  By notations in~\eqref{eq:t-th-iter} and 
	$\vy^{(t)}=(\vy^{(t)\top}_1, \ldots,\vy^{(t)\top}_{3m_2-1})\zz$ with $\vy^{(t)}_j\in\RR^{{\bar{d}}}$. It holds, for any $\bar{j}\in\{2,3,\dots,\bar{d}\}$,  that 
	\begin{align}
		\label{eq:outerinduction3}
		\textnormal{supp}(\vx_i^{(t)})\subset\{1,\dots,\bar{j}-1\} \text{ and }	\textnormal{supp}(\vy_j^{(t)})\subset\{1,\dots,\bar{j}-1\}
	\end{align}
	$\text{ for }i=1,\dots,m,~j=1,\dots,3m_2-1 \text{ and }t\leq 1+ m(\bar{j}-2)/3.$
\end{proposition}

\proof{Proof of Proposition \ref{thm:iterateguess3}.}
	We prove this claim by induction on $\bar{j}$. %Before  that, we derived some properties of the operations allowed by Assumption~\ref{ass:linearspan3}.
	Let $\vxi^{(t)}=(\vxi^{(t)\top}_1, \ldots,\vxi^{(t)\top}_{3m_2-1})\zz$ with $\vxi^{(t)}_j\in\RR^{{\bar{d}}}$ and $\vzeta^{(t)}=\prox_{\eta_t \bar{g}}(\vxi^{(t)})=(\vzeta^{(t)\top}_1, \ldots,\vzeta^{(t)\top}_{3m_2-1})\zz$  with $\vzeta^{(t)}_j\in\RR^{{\bar{d}}},\forall\, j$ for $t\geq1$.
	Since the algorithm is initialized with $\vx^{(0)}=\mathbf{0}$ and $\vy^{(0)}=\mathbf{0}$, we have $\textnormal{supp}(\nabla f_i(\vx^{(0)}_i))\subset\{1\}$ for any $i$ according to Lemma~\ref{lem:iterateguess}. Notice $\vb=\textbf{0}$ and $\bar{\vb}=\textbf{0}$. By Algorithm Class~\ref{algo:class-2} and \eqref{eq:suppg3},  we have $\textnormal{supp}(\vx_i^{(1)})\subset\{1\}$ for any $i$. Meanwhile, we have $\vxi^{(1)}_j=\vzero$ and $\vzeta^{(1)}_j=\vzero$ for any $j$. This implies $\vy_j^{(1)}=\mathbf{0}$.  Thus the claim in \eqref{eq:outerinduction3} holds for $\bar{j}=2$. Suppose that we have proved the claim in \eqref{eq:outerinduction3}  for all $\bar{j}\geq 2$. We next prove it for $\bar{j}+1$. 
	According to the hypothesis of the induction, we have
	\begin{equation}
		\begin{aligned}
			\label{eq:outerinduction3-inner}
			&\textnormal{supp}(\vx_i^{(t)})\subset\{1,\dots,\bar{j}-1\} \text{ and }	\textnormal{supp}(\vy_j^{(t)})\subset\{1,\dots,\bar{j}-1\},
			\\&\forall\, i\in[1,m],\,\,\forall j\in[1,3m_2-1] \text{ and }r\leq \bar{t}:= 1+ m(\bar{j}-2)/3.
		\end{aligned}
	\end{equation}
	
	Below we let $\widehat\vx^{(s)}$ be any vector in $\mathbf{span}\left\{\vA^{\top}\vA\vx^{(s)},\bar\vA^{\top}\bar\vA\vx^{(s)}\right\}$, %for any $s\geq0$,
	$\widetilde \vx^{(s)}=\bar\vA^\top \vy^{(s)}$ and $\widehat\vy^{(s)}=\bar\vA\vx^{(s)}$ for any $s\geq0$, %Consider iteration $t$ with $t\leq 1+ m(\bar{j}-2)/3$. %According to the hypothesis of the induction, we have $\textnormal{supp}(\vx_i^{(s)})\subset\{1,\dots,\bar{j}-1\}$ and $\textnormal{supp}(\vy_j^{(s)})\subset\{1,\dots,\bar{j}-1\}$  for any $i$ and $s\leq t$. 
	and	we  consider two cases: $\bar{j}$ is even and $\bar{j}$ is odd. 
	
	\textbf{Case 1}: Suppose $\bar{j}$ is even.  We claim that, for $s=0,1,\dots,\frac{m}{3}$,
	\begin{eqnarray}
		\label{eq:innerinduction13}
		&\textnormal{supp}(\vx_i^{(r)})\subset\left\{
		\begin{array}{ll}
			\{1,\dots,\bar{j}\},& \text{ if }i\in\left[1, \frac{m}{3}+s\right],\\
			\{1,\dots,\bar{j}-1\},&\text{ if } i\in\left[\frac{m}{3}+s+1, m\right],
		\end{array}
		\right.\,\, \forall r\in[\bar{t}+s]
		\text{ and } \\ \label{eq:innerinduction23}
		&\textnormal{supp}( \vy_j^{(r)})\subset\left\{
		\begin{array}{ll}
			\{1,\dots,\bar{j}\},& \text{ if }j\in\left[1, m_2+\lfloor \frac{s}{m_1}\rfloor \right],\\
			\{1,\dots,\bar{j}-1\},& \text{ if }j\in\left[m_2+\lfloor \frac{s}{m_1}\rfloor+1, 3m_2-1\right],
		\end{array}
		\right.\,\, \forall r\in[\bar{t}+s].
	\end{eqnarray} 
	Notice~\eqref{eq:outerinduction3-inner} implies~\eqref{eq:innerinduction13} and~\eqref{eq:innerinduction23} for $s=0$. 
	%Since $\textnormal{supp}(\vx_i^{(r)})\subset\{1,\dots,\bar{j}-1\}$ and $\textnormal{supp}(\vy_j^{(r)})\subset\{1,\dots,\bar{j}-1\}$ for any $i$, $j$ and $r\leq t$, this claim holds for $s=0$. 
	Suppose~\eqref{eq:innerinduction13} and~\eqref{eq:innerinduction23} hold for an integer  $s$ satisfying $0\leq s\leq\frac{m}{3}$. According to Lemma~\ref{lem:iterateguess} and $\frac{m}{3}+s \leq \frac{2m}{3}$, 
	$$
	\textnormal{supp}(\nabla f_i(\vx_i^{(r)}))\subset\left\{
	\begin{array}{ll}
		\{1,\dots,\bar{j}\},& \text{ if }i\in\left[1, \frac{m}{3}+s\right],\\
		\{1,\dots,\bar{j}-1\},&\text{ if } i\in\left[\frac{m}{3}+s+1, m\right],
	\end{array}
	\right.\,\, \forall r\in[\bar{t}+s].
	$$
	In addition, by Lemma~\ref{lem:supp2}(a), we have from~\eqref{eq:suppAy} that
	$$
	\textnormal{supp}(\widehat\vx_i^{(r)}),~\textnormal{supp}(\widetilde\vx_i^{(r)}) \subset\left\{
	\begin{array}{ll}
		\{1,\dots,\bar{j}\},&\text{ if } i\in\left[1, \frac{m}{3}+s+1\right],\\
		\{1,\dots,\bar{j}-1\},&\text{ if } i\in\left[\frac{m}{3}+s+2, m\right],
	\end{array}
	\right.\,\, \forall r\in[\bar{t}+s].
	%\text{ and }
	%\textnormal{supp}(\widetilde\vx_i^{(r)})\subset\left\{
	%\begin{array}{ll}
	%	\{1,\dots,\bar{j}\}& i\in\left[1, \frac{m}{3}+s+1\right]\\
	%	\{1,\dots,\bar{j}-1\}& i\in\left[\frac{m}{3}+s+2, m\right]
	%\end{array}
	%\right.
	$$
	Hence, by Algorithm Class~\ref{algo:class-2}, we have
	$$
	\textnormal{supp}(\vx_i^{(\bar{t}+s+1)})\subset\left\{
	\begin{array}{ll}
		\{1,\dots,\bar{j}\},&\text{ if } i\in\left[1, \frac{m}{3}+s+1\right],\\
		\{1,\dots,\bar{j}-1\},&\text{ if } i\in\left[\frac{m}{3}+s+2, m\right].
	\end{array}
	\right.
	$$
	This means the claim in \eqref{eq:innerinduction13} holds for $s+1$ as well. 
	
	In addition, by the relation in \eqref{eq:suppAbar3}, we have
	$$
	\textnormal{supp}(\widehat\vy_j^{(r)})\subset\left\{
	\begin{array}{ll}
		\{1,\dots,\bar{j}\},& \text{ if }j\in\left[1, m_2+\lfloor \frac{s}{m_1}\rfloor \right],\\
		\{1,\dots,\bar{j}-1\},&\text{ if } j\in\left[m_2+\lfloor \frac{s}{m_1}\rfloor+1, 3m_2-1\right],
	\end{array}
	\right.\,\, \forall r\in[\bar{t}+s].
	$$
	Together with \eqref{eq:innerinduction23} and \eqref{eq:suppAbarAbary}, the inclusion above implies that 
	$$
	\textnormal{supp}(\vxi_j^{(\bar{t}+s+1)})\subset\left\{
	\begin{array}{ll}
		\{1,\dots,\bar{j}\},&\text{ if } j\in\left[1, m_2+\lfloor \frac{s}{m_1}\rfloor \right],\\
		\{1,\dots,\bar{j}-1\},&\text{ if } j\in\left[m_2+\lfloor \frac{s}{m_1}\rfloor+1, 3m_2-1\right].
	\end{array}
	\right.
	$$
	It then follows from~\eqref{eq:suppg3} that
	$$
	\textnormal{supp}(\vzeta_j^{(\bar{t}+s+1)})\subset\left\{
	\begin{array}{ll}
		\{1,\dots,\bar{j}\},&\text{ if } j\in\left[1, m_2+\lfloor \frac{s}{m_1}\rfloor \right],\\
		\{1,\dots,\bar{j}-1\},&\text{ if } j\in\left[m_2+\lfloor \frac{s}{m_1}\rfloor+1, 3m_2-1\right].
	\end{array}
	\right.
	$$
	By Algorithm Class~\ref{algo:class-2}, we have
	$$
	\textnormal{supp}(\vy_j^{(\bar{t}+s+1)})\subset\left\{
	\begin{array}{ll}
		\{1,\dots,\bar{j}\},&\text{ if } j\in\left[1, m_2+\lfloor \frac{s}{m_1}\rfloor \right],\\
		\{1,\dots,\bar{j}-1\},&\text{ if } j\in\left[m_2+\lfloor \frac{s}{m_1}\rfloor+1, 3m_2-1\right].
	\end{array}
	\right.
	$$
	This means the claim in \eqref{eq:innerinduction23} holds for $s+1$ as well. By induction, \eqref{eq:innerinduction13} and~\eqref{eq:innerinduction23} hold for $s=0,1,\dots,\frac{m}{3}$. Let $s=\frac{m}{3}$ in them. We have $\textnormal{supp}(\vx_i^{(r)})\subset\{1,\dots,\bar{j}\}$ and $\textnormal{supp}(\vy_j^{(r)})\subset\{1,\dots,\bar{j}\}$ for any $i$, $j$ and $r\leq \bar{t}+\frac{m}{3}= 1+ m(\bar{j}-2)/3+\frac{m}{3}=1+m(\bar{j}-1)/3$. 
	
	\textbf{Case 2}: Suppose $\bar{j}$ is odd.  We claim that, for $s=0,1,\dots,\frac{m}{3}$,
	\begin{eqnarray}
		\label{eq:innerinduction33}
		\textnormal{supp}(\vx_i^{(r)})\subset\left\{
		\begin{array}{ll}
			\{1,\dots,\bar{j}-1\},&\text{ if } i\in\left[1, \frac{2m}{3}-s\right],\\
			\{1,\dots,\bar{j}\},&\text{ if } i\in\left[\frac{2m}{3}-s+1, m\right],
		\end{array}
		\right.\,\, \forall r\in[\bar{t}+s],
		\text{ and }\\\label{eq:innerinduction43}
		\textnormal{supp}( \vy_j^{(r)})\subset\left\{
		\begin{array}{ll}
			\{1,\dots,\bar{j}-1\},&\text{ if } j\in\left[1, 2m_2-\lceil \frac{s}{m_1}\rceil \right],\\
			\{1,\dots,\bar{j}\},&\text{ if } j\in\left[2m_2-\lceil \frac{s}{m_1}\rceil+1, 3m_2-1\right],
		\end{array}
		\right.\,\, \forall r\in[\bar{t}+s].
	\end{eqnarray} 
	Again~\eqref{eq:outerinduction3-inner} implies~\eqref{eq:innerinduction33} and~\eqref{eq:innerinduction43} for $s=0$. 
	%Since $\textnormal{supp}(\vx_i^{(r)})\subset\{1,\dots,\bar{j}-1\}$ for any $i$ and $r\leq t$, this claim holds for $s=0$. 
	Suppose~\eqref{eq:innerinduction33} and~\eqref{eq:innerinduction43} hold for an integer $s$ satisfying $0\leq s\leq\frac{m}{3}$. According to Lemma~\ref{lem:iterateguess} and $\frac{2m}{3}-s\geq \frac{m}{3}$, it holds
	$$
	\textnormal{supp}(\nabla f_i(\vx_i^{(r)}))\subset\left\{
	\begin{array}{ll}
		\{1,\dots,\bar{j}-1\},&\text{ if } i\in\left[1, \frac{2m}{3}-s\right],\\
		\{1,\dots,\bar{j}\},&\text{ if } i\in\left[\frac{2m}{3}-s+1, m\right],
	\end{array}
	\right.\,\, \forall r\in[\bar{t}+s].
	$$
	In addition, by Lemma~\ref{lem:supp2}(a) and~\eqref{eq:suppAy}, we have 
	$$
	\textnormal{supp}(\widehat\vx_i^{(r)}),~\textnormal{supp}(\widetilde\vx_i^{(r)}) \subset\left\{
	\begin{array}{ll}
		\{1,\dots,\bar{j}-1\},&\text{ if } i\in\left[1, \frac{2m}{3}-s-1\right],\\
		\{1,\dots,\bar{j}\},&\text{ if } i\in\left[\frac{2m}{3}-s, m\right],
	\end{array}
	\right.\,\, \forall r\in[\bar{t}+s].
	$$
	Hence, by Algorithm Class~\ref{algo:class-2}, we have
	$$
	\textnormal{supp}(\vx_i^{(\bar{t}+s+1)})\subset\left\{
	\begin{array}{ll}
		\{1,\dots,\bar{j}-1\},&\text{ if } i\in\left[1, \frac{2m}{3}-s-1\right],\\
		\{1,\dots,\bar{j}\},&\text{ if } i\in\left[\frac{2m}{3}-s, m\right].
	\end{array}
	\right.
	$$
	This means claim~\eqref{eq:innerinduction33} holds for $s+1$ as well. 
	
	In addition, by \eqref{eq:suppAbar3}, we have
	$$
	\textnormal{supp}(\widehat\vy_j^{(\bar{t}+s+1)})\subset\left\{
	\begin{array}{ll}
		\{1,\dots,\bar{j}-1\},&\text{ if } j\in\left[1, 2m_2-\lceil \frac{s+1}{m_1}\rceil \right],\\
		\{1,\dots,\bar{j}\},&\text{ if } j\in\left[2m_2-\lceil \frac{s+1}{m_1}\rceil+1, 3m_2-1\right].
	\end{array}
	\right.
	$$
	Together with \eqref{eq:innerinduction43}, the inclusion above implies that 
	$$
	\textnormal{supp}(\vxi_j^{(\bar{t}+s+1)})\subset\left\{
	\begin{array}{ll}
		\{1,\dots,\bar{j}-1\},&\text{ if } j\in\left[1, 2m_2-\lceil \frac{s+1}{m_1}\rceil \right],\\
		\{1,\dots,\bar{j}\},&\text{ if } j\in\left[2m_2-\lceil \frac{s+1}{m_1}\rceil+1, 3m_2-1\right].
	\end{array}
	\right.
	$$
	It then follows from~\eqref{eq:suppg3} that
	$$
	\textnormal{supp}(\vzeta_j^{(\bar{t}+s+1)})\subset\left\{
	\begin{array}{ll}
		\{1,\dots,\bar{j}-1\},&\text{ if } j\in\left[1, 2m_2-\lceil \frac{s+1}{m_1}\rceil \right],\\
		\{1,\dots,\bar{j}\},&\text{ if } j\in\left[2m_2-\lceil \frac{s+1}{m_1}\rceil+1, 3m_2-1\right].
	\end{array}
	\right.
	$$
	By Algorithm Class~\ref{algo:class-2}, we have
	$$
	\textnormal{supp}(\vy_j^{(\bar{t}+s+1)})\subset\left\{
	\begin{array}{ll}
		\{1,\dots,\bar{j}-1\},&\text{ if } j\in\left[1, 2m_2-\lceil \frac{s+1}{m_1}\rceil \right],\\
		\{1,\dots,\bar{j}\},&\text{ if } j\in\left[2m_2-\lceil \frac{s+1}{m_1}\rceil+1, 3m_2-1\right],
	\end{array}
	\right.
	$$
	which means \eqref{eq:innerinduction43} holds for $s+1$ as well. By induction, \eqref{eq:innerinduction33} and~\eqref{eq:innerinduction43} holds for $s=0,1,\dots,\frac{m}{3}$. Let $s=\frac{m}{3}$ in~\eqref{eq:innerinduction43}. We have $\textnormal{supp}(\vx_i^{(r)})\subset\{1,\dots,\bar{j}\}$ and $\textnormal{supp}(\vy_j^{(r)})\subset\{1,\dots,\bar{j}\}$ for any $i$, $j$ and $r\leq \bar{t}+\frac{m}{3}= 1+ m(\bar{j}-2)/3+\frac{m}{3}=1+m(\bar{j}-1)/3$. 
	
	Therefore, we have proved that~\eqref{eq:outerinduction3-inner} holds for $\bar{j}+1$, when $\bar{j}$ is either even or odd. By induction,~\eqref{eq:outerinduction3-inner} holds for any integer $\bar{j}\in[2,\bar{d}]$, and we complete the proof.	
	\null\hfill$\square$
\endproof

Now, we are ready to prove Theorem \ref{thm:lbcomposite}.

%\begin{theorem}
%	Let $\epsilon>0$ and $L_f > 0$ be given.				Suppose an algorithm  is applied to problem~\eqref{eq:model-spli} that satisfies Assumption~\ref{assume:problemsetup} 
%	and generates a sequence $\{(\vx^{(t)},\vy^{(t)})\}_{t\geq0}$ that satisfies Assumption~\ref{ass:linearspan3}. Then there exists an instance of problem~\eqref{eq:model-spli}, i.e., instance~$\mathcal{P}$ given in Definition~\ref{def:hardinstance}, such that the algorithm requires at least 
%	$\left \lceil \frac{\kappa([\bar{\vA}; \vA]) L_f \Delta_F}{72000\pi^2} \epsilon^{-2} \right\rceil$ 
%	iterations to obtain an $\epsilon$-stationary point of instance~$\mathcal{P}$. Here, it holds that $\Delta_F~=~F(\vx^{(0)},\vy^{(0)})-\inf_{\vx,\vy} F(\vx,\vy)$. %as defined in~\eqref{eq:def-Delta}.
%\end{theorem}

\proof{Proof of Theorem \ref{thm:lbcomposite}.}
	As we discussed below~\eqref{eq:t-th-iter}, we assume $\vx^{(0)}=\mathbf{0}$ and $\vy^{(0)}=\mathbf{0}$ without loss of generality.  Thus by notation in~\eqref{eq:t-th-iter}, Proposition~\ref{thm:iterateguess3} indicates that 	$\textnormal{supp}(\vx_i^{(t)})\subset\{1,\dots,\bar{d}-1\}$ and $\textnormal{supp}(\vy_j^{(t)})\subset\{1,\dots,\bar{d}-1\}$ for $i=1,\dots,m$ and $j=1,2,\ldots,3m_2-1$ for all $t\leq 1+ m(\bar{d}-2)/3$, which means $[\bar{\vx}^{(t)}]_{\bar{d}}=0$ if  $t\leq 1+ m(\bar{d}-2)/3$. 
	Hence, by  Lemmas~\ref{lem:nablaf} and~\ref{lem:kktvio}, we have
	$$
	\max\left\{\left\|\vH\vx^{(t)}\right\|, \min _{\vgamma}\left\|\nabla  f_0(\vx^{(t)})+\vH^{\top} \vgamma\right\|\right\} \geq \frac{\sqrt{m}}{2}\left\|\frac{1}{m}\sum_{i=1}^m\nabla f_i(\bar{\vx}^{(t)})\right\|>\epsilon,\,\,\forall t\leq 1+ m(\bar{d}-2)/3.
	$$ 						
	Hence, $\vx^{(t)}$ is not an $\epsilon$-stationary point of problem~\eqref{eq:model3} if $t\leq m({\bar{d}}-1)/3$. Thus, $(\vx^{(t)},\vy^{(t)})$
	cannot be an $\epsilon/2$-stationary point of the reformulation~\eqref{eq:model-spli} of instance $\cP$ according to Corollary~\ref{cor:kktequiv3}. Moreover, by 				
	Lemma~\ref{cor:boundf}(a) and the facts that $g(\vx^{(0)})=0$ and $g(\vx) \ge0$, %$\inf_\vy\bar g(\vy)=\bar g(\vy^{(0)})=0$, $\inf_{\vx,\vy} [f_0\left(\vx\right)+\bar g\left(\vy\right)]\geq \inf_{\vx} f_0\left(\vx\right)$, 
	it holds that
	\begin{equation}\label{eq:l-bd-d-bar} {\bar{d}}\geq \frac{L_f \left( f_0(\vx^{(0)})-\inf_{\vx} f_0\left(\vx\right)\right)}{3000\pi^2} \epsilon^{-2}
		\geq \frac{L_f \left( F_0(\vx^{(0)})-\inf_{\vx} F_0(\vx) \right)}{3000\pi^2} \epsilon^{-2}
		=\frac{L_f\Delta_{F_0}}{3000\pi^2}\epsilon^{-2}.
	\end{equation}
	In other words, in order for $(\vx^{(t)},\vy^{(t)})$ to be an $\epsilon/2$-stationary point of~\eqref{eq:model-spli}, %according to an instance~$\mathcal{P}$, 
	the algorithm needs at least $t= 2+m({\bar{d}}-2)/3$ oracles. Notice
	\begin{equation}
		\label{eq:statement}
		\textstyle 	2+m({\bar{d}}-2)/3\geq m{\bar{d}}/6 \overset{\eqref{eq:l-bd-d-bar}}\geq\frac{m L_f \Delta_{F_0}}{18000\pi^2} \epsilon^{-2}%=\frac{m_1m_2 L_f \Delta_F}{6000\pi^2} \epsilon^{-2}
		>\frac{\kappa([\bar{\vA}; \vA]) L_f \Delta_{F_0}}{18000\pi^2}\epsilon^{-2},
	\end{equation}
	where the second inequality is because ${\bar{d}}\geq 5$ and the last inequality is by Lemma~\ref{lem:condH}. The conclusion is then proved by replacing $\epsilon$ in \eqref{eq:statement} to $2\epsilon$.
	\null\hfill$\square$
\endproof

Finally, we give the proof to Corollary \ref{cor:lower2}.

%\begin{corollary}
%	Let $\epsilon>0$ and $L_f > 0$ be given.				Suppose an algorithm  is applied to problem~\eqref{eq:model} that satisfies Assumption~\ref{assume:problemsetup} 
%	and generates a sequence $\{(\vx^{(t)},\vy^{(t)})\}_{t\geq0}$ that satisfies Assumption~\ref{ass:linearspan3}. Then for any $\omega\in [0, \frac{150\pi \epsilon}{L_f})$, there exists an instance of problem~\eqref{eq:model}, i.e., instance~$\mathcal{P}$ in Definition~\ref{def:hardinstance}, such that the algorithm requires at least 
%	$\left \lceil \frac{\kappa([\bar{\vA}; \vA]) L_f \Delta_{F_0}}{18000\pi^2} \epsilon^{-2} \right\rceil$ 
%	iterations to obtain a point $\vx^{(t)}$ that is $\omega$-close to an $\epsilon$-stationary point of instance~$\mathcal{P}$, where $\Delta_{F_0}=F_0(\vx^{(0)})-\inf_\vx F_0(\vx)$. %as defined in~\eqref{eq:def-Delta}.
%\end{corollary}

\proof{Proof of Corollary \ref{cor:lower2}.}
	As we discussed in Section \ref{sec:extension}, 	we assume $\vx^{(0)}=\mathbf{0}$ and $\vy^{(0)}=\mathbf{0}$ without loss of generality. 
	Thus by notation in~\eqref{eq:t-th-iter}, Proposition~\ref{thm:iterateguess3} indicates that 	$\textnormal{supp}(\vx_i^{(t)})\subset\{1,\dots,\bar{d}-1\}$ for any $i \in [1,m]$ and any $t\leq 1+ m(\bar{d}-2)/3$, which means $[\bar{\vx}^{(t)}]_{\bar{d}}=0$ if  $t\leq 1+ m(\bar{d}-2)/3$, where $\bar{\vx}^{(t)}= \frac{1}{m} \sum_{i=1}^m \vx_i^{(t)}$.  
	
	On the other hand, suppose $\vx^*$ with the structure as in~\eqref{eq:t-th-iter} is an $\epsilon$-stationary point of instance~$\mathcal{P}$. Then by 	Lemma~\ref{cor:kktvio2}, it must also be an $\epsilon$-stationary point of~\eqref{eq:model3}. Hence, by  Lemmas~\ref{lem:nablaf} and~\ref{lem:kktvio}, we have $\left|[\bar\vx^*]_j\right| \ge \frac{150\pi\epsilon}{\sqrt{m} L_f}$ for all $j=1,\ldots,\bar d$, where $\bar{\vx}^*= \frac{1}{m} \sum_{i=1}^m \vx_i^*$.  Therefore, by the convexity of the square function, it follows that
	$$\textstyle \|\vx^{(t)} - \vx^*\|^2 \ge \sum_{i=1}^m \left([\vx_i^{(t)}]_{\bar d} - [\vx_i^*]_{\bar d}\right)^2 \ge m \left([\bar\vx^{(t)}]_{\bar d} - [\bar\vx^*]_{\bar d}\right)^2 \ge m \left(\frac{150\pi\epsilon}{\sqrt{m} L_f}\right)^2 > \omega^2,$$
	and thus $\vx^{(t)}$ is not $\omega$-close to $\vx^*$ if  $t\leq 1+ m(\bar{d}-2)/3$.
	%
	%	Moreover, by 				
	%	Lemma~\ref{cor:boundf}(a) and the fact that $g(\vx^{(0)})=0$ and $g(\vx) \ge 0, \forall\, \vx$, it holds that
	%	$$ {\bar{d}}\geq \frac{L_f \left( F_0(\vx^{(0)})-\inf_{\vx} F_0\left(\vx\right)\right)}{3000\pi^2} \epsilon^{-2}=\frac{L_f\Delta_{F_0}}{3000\pi^2}\epsilon^{-2}.
	%	$$
	In other words, in order for $\vx^{(t)}$ to be $\omega$-close to an $\epsilon$-stationary point of instant~$\mathcal{P}$, the algorithm needs at least $t=2+m({\bar{d}}-2)/3$ oracles.
	The proof is then completed by using \eqref{eq:statement}. %observing	
	%					\begin{equation*}
		%		2+m({\bar{d}}-2)/3\geq m{\bar{d}}/6\geq\frac{m L_f \Delta_{F_0}}{18000\pi^2} \epsilon^{-2}>\frac{\kappa([\bar{\vA}; \vA]) L_f \Delta_{F_0}}{18000\pi^2},
		%	\end{equation*}
	%	where the first inequality is because ${\bar{d}}\geq 5$, and the last one is by Lemma~\ref{lem:condH}. 
	\null\hfill$\square$
\endproof
\end{APPENDICES}

\section*{Acknowledgements}

This work is partly supported by NSF awards DMS-2053493 and IIS-2147253 and the ONR award N00014-22-1-2573.

	\bibliographystyle{plain}
	\bibliography{optim}

%The funding sources have no role in the study design, data collection, analysis, interpretation, or in the decision to publish the manuscript.

%Besides, the authors declare no financial or non-financial or competing interests that could be perceived as influencing the outcome of this research or the content of this manuscript. Human participants and animals are not involved in this study.
				
				%\thispagestyle{style2}

				%\input LBN_cvx2.bbl
			\end{document}

%% file: macros.tex
\usepackage{mathtools}
\usepackage[normalem]{ulem} % to use \sout

\newcommand{\vb}{{\mathbf{b}}}

\newcommand{\vu}{{\mathbf{u}}}
\newcommand{\vv}{{\mathbf{v}}}

\newcommand{\vx}{{\mathbf{x}}}
\newcommand{\vy}{{\mathbf{y}}}
\newcommand{\vz}{{\mathbf{z}}}

\newcommand{\vA}{{\mathbf{A}}}
\newcommand{\vB}{{\mathbf{B}}}

\newcommand{\vH}{{\mathbf{H}}}
\newcommand{\vI}{{\mathbf{I}}}
\newcommand{\vJ}{{\mathbf{J}}}

\newcommand{\cC}{{\mathcal{C}}}

\newcommand{\cP}{{\mathcal{P}}}

\newcommand{\cS}{{\mathcal{S}}}

\newcommand{\cX}{{\mathcal{X}}}

\newcommand{\RR}{\mathbb{R}} % real
\newcommand{\zz}{^{\top}} % integer
\newcommand{\vzero}{\mathbf{0}} % 0 vector
\newcommand{\dist}{\mathrm{dist}}    % distance
\newcommand{\prox}{{\mathbf{prox}}} % proximal map
\newcommand{\dom}{{\mathrm{dom}}} % domain
\newcommand{\st}{\mbox{ s.t. }}

\newcommand{\bc}{\begin{center}}
\newcommand{\ec}{\end{center}}

\newcommand{\bdm}{\begin{displaymath}}
\newcommand{\edm}{\end{displaymath}}

\newcommand{\beq}{\begin{equation}}
\newcommand{\eeq}{\end{equation}}

\newcommand{\bfl}{\begin{flushleft}}
\newcommand{\efl}{\end{flushleft}}

\newcommand{\bt}{\begin{tabbing}}
\newcommand{\et}{\end{tabbing}}

\newcommand{\beqn}{\begin{eqnarray}}
\newcommand{\eeqn}{\end{eqnarray}}

\newcommand{\vgamma}{{\boldsymbol{\gamma}}}
\newcommand{\vxi}{{\boldsymbol{\xi}}}
\newcommand{\vzeta}{{\boldsymbol{\zeta}}}

\newcommand{\beqs}{\begin{align*}} % no equation numbers
\newcommand{\eeqs}{\end{align*}}  % no equation numbers

\newtheorem{algo}{Algorithm Class}
\newtheorem{prob}{Problem Class}

%% file: LB_NCVX_affineConstraint.bbl
\begin{thebibliography}{50}
\providecommand{\natexlab}[1]{#1}
\providecommand{\url}[1]{\texttt{#1}}
\providecommand{\urlprefix}{URL }

\bibitem[{Bian \protect\BIBand{} Chen(2015)}]{bian2015linearly}
Bian W, Chen X (2015) Linearly constrained non-{Lipschitz} optimization for
  image restoration. \emph{SIAM Journal on Imaging Sciences} 8(4):2294--2322.

\bibitem[{Carmon et~al.(2020)Carmon, Duchi, Hinder, \protect\BIBand{}
  Sidford}]{carmon2020lower}
Carmon Y, Duchi JC, Hinder O, Sidford A (2020) Lower bounds for finding
  stationary points {I}. \emph{Mathematical Programming} 184(1-2):71--120.

\bibitem[{Carmon et~al.(2021)Carmon, Duchi, Hinder, \protect\BIBand{}
  Sidford}]{carmon2021lower}
Carmon Y, Duchi JC, Hinder O, Sidford A (2021) Lower bounds for finding
  stationary points {II}: first-order methods. \emph{Mathematical Programming}
  185(1-2):315--355.

\bibitem[{Cartis et~al.(2012)Cartis, Gould, \protect\BIBand{}
  Toint}]{Cartis2012how}
Cartis C, Gould N, Toint P (2012) How much patience do you have? a worst-case
  perspective on smooth nonconvex optimization. \emph{Optima. Mathematical
  Optimization Society Newsletter} (88):1--10.

\bibitem[{Cartis et~al.(2014)Cartis, Gould, \protect\BIBand{}
  Toint}]{cartis2014complexity}
Cartis C, Gould N, Toint P (2014) On the complexity of finding first-order
  critical points in constrained nonlinear optimization. \emph{Mathematical
  Programming} 144(1-2):93--106.

\bibitem[{Cartis et~al.(2017)Cartis, Gould, \protect\BIBand{}
  Toint}]{cartis2017corrigendum}
Cartis C, Gould N, Toint P (2017) Corrigendum: on the complexity of finding
  first-order critical points in constrained nonlinear optimization.
  \emph{Mathematical Programming} 161:611--626.

\bibitem[{Clarke(1990)}]{clarke1990optimization}
Clarke FH (1990) \emph{{Optimization and Nonsmooth Analysis}}, volume~5
  (Philadelphia: SIAM).

\bibitem[{Dahal et~al.(2023)Dahal, Liu, \protect\BIBand{} Xu}]{dahal2023damped}
Dahal H, Liu W, Xu Y (2023) Damped proximal augmented {Lagrangian} method for
  weakly-convex problems with convex constraints.
  \urlprefix\url{https://arxiv.org/pdf/2311.09065.pdf}.

\bibitem[{Goncalves et~al.(2017)Goncalves, Melo, \protect\BIBand{}
  Monteiro}]{goncalves2017convergence}
Goncalves MLN, Melo JG, Monteiro RDC (2017) Convergence rate bounds for a
  proximal {ADMM} with over-relaxation stepsize parameter for solving nonconvex
  linearly constrained problems.
  \urlprefix\url{https://arxiv.org/pdf/1702.01850.pdf}.

\bibitem[{Grapiglia \protect\BIBand{} Yuan(2021)}]{grapiglia2021complexity}
Grapiglia GN, Yuan Y (2021) On the complexity of an augmented {Lagrangian}
  method for nonconvex optimization. \emph{IMA Journal of Numerical Analysis}
  41(2):1546--1568.

\bibitem[{Gray(2006)}]{gray2006toeplitz}
Gray RM (2006) Toeplitz and circulant matrices: a review. \emph{Foundations and
  Trends{\textregistered} in Communications and Information Theory}
  2(3):155--239.

\bibitem[{Haeser et~al.(2019)Haeser, Liu, \protect\BIBand{}
  Ye}]{haeser2019optimality}
Haeser G, Liu H, Ye Y (2019) Optimality condition and complexity analysis for
  linearly-constrained optimization without differentiability on the boundary.
  \emph{Mathematical Programming} 178:263--299.

\bibitem[{Hajinezhad \protect\BIBand{} Hong(2019)}]{hajinezhad2019perturbed}
Hajinezhad D, Hong M (2019) Perturbed proximal primal--dual algorithm for
  nonconvex nonsmooth optimization. \emph{Mathematical Programming}
  176(1-2):207--245.

\bibitem[{Hong(2016)}]{hong2016decomposing}
Hong M (2016) Decomposing linearly constrained nonconvex problems by a proximal
  primal dual approach: algorithms, convergence, and applications.
  \urlprefix\url{https://arxiv.org/pdf/1604.00543.pdf}.

\bibitem[{Hong et~al.(2016)Hong, Luo, \protect\BIBand{}
  Razaviyayn}]{hong2016convergence}
Hong M, Luo Z, Razaviyayn M (2016) Convergence analysis of alternating
  direction method of multipliers for a family of nonconvex problems.
  \emph{SIAM Journal on Optimization} 26(1):337--364.

\bibitem[{Huang \protect\BIBand{} Lin(2023)}]{huang2023single}
Huang Y, Lin Q (2023) Oracle complexity of single-loop switching subgradient
  methods for non-smooth weakly convex functional constrained optimization.
  \emph{Advances in Neural Information Processing Systems} 36:61327--61340.

\bibitem[{Jiang et~al.(2019)Jiang, Lin, Ma, \protect\BIBand{}
  Zhang}]{jiang2019structured}
Jiang B, Lin T, Ma S, Zhang S (2019) Structured nonconvex and nonsmooth
  optimization: algorithms and iteration complexity analysis.
  \emph{Computational Optimization and Applications} 72(1):115--157.

\bibitem[{Kong et~al.(2019)Kong, Melo, \protect\BIBand{}
  Monteiro}]{kong2019complexity}
Kong W, Melo JG, Monteiro RDC (2019) Complexity of a quadratic penalty
  accelerated inexact proximal point method for solving linearly constrained
  nonconvex composite programs. \emph{SIAM Journal on Optimization}
  29(4):2566--2593.

\bibitem[{Kong et~al.(2023)Kong, Melo, \protect\BIBand{}
  Monteiro}]{kong2023iteration}
Kong W, Melo JG, Monteiro RDC (2023) Iteration complexity of a proximal
  augmented {Lagrangian} method for solving nonconvex composite optimization
  problems with nonlinear convex constraints. \emph{Mathematics of Operations
  Research} 48(2):1066--1094.

\bibitem[{Lan(2020)}]{lan2020first}
Lan G (2020) \emph{First-order and stochastic optimization methods for machine
  learning}, volume~1 (Springer).

\bibitem[{Li et~al.(2021)Li, Tian, Zhang, \protect\BIBand{}
  Jadbabaie}]{li2021complexity}
Li H, Tian Y, Zhang J, Jadbabaie A (2021) Complexity lower bounds for
  nonconvex-strongly-concave min-max optimization. \emph{Advances in Neural
  Information Processing Systems} 34:1792--1804.

\bibitem[{Li \protect\BIBand{} Xu(2021)}]{li2021augmented}
Li Z, Xu Y (2021) Augmented {Lagrangian}--based first-order methods for
  convex-constrained programs with weakly convex objective. \emph{INFORMS
  Journal on Optimization} 3(4):373--397.

\bibitem[{Lin et~al.(2022)Lin, Ma, \protect\BIBand{} Xu}]{lin2022complexity}
Lin Q, Ma R, Xu Y (2022) Complexity of an inexact proximal-point penalty method
  for constrained smooth non-convex optimization. \emph{Computational
  Optimization and Applications} 82(1):175--224.

\bibitem[{Lin et~al.(2020)Lin, Jin, \protect\BIBand{} Jordan}]{lin2020gradient}
Lin T, Jin C, Jordan M (2020) On gradient descent ascent for nonconvex-concave
  minimax problems. \emph{International Conference on Machine Learning},
  6083--6093 (PMLR).

\bibitem[{Liu et~al.(2023{\natexlab{a}})Liu, Lin, \protect\BIBand{}
  Xu}]{liu2023first}
Liu W, Lin Q, Xu Y (2023{\natexlab{a}}) First-order methods for affinely
  constrained composite non-convex non-smooth problems: Lower complexity bound
  and near-optimal methods.
  \urlprefix\url{https://arxiv.org/pdf/2307.07605.pdf}.

\bibitem[{Liu et~al.(2022)Liu, Liu, \protect\BIBand{} Chen}]{liu2022linearly}
Liu W, Liu X, Chen X (2022) Linearly constrained nonsmooth optimization for
  training autoencoders. \emph{SIAM Journal on Optimization} 32(3):1931--1957.

\bibitem[{Liu et~al.(2023{\natexlab{b}})Liu, Liu, \protect\BIBand{}
  Chen}]{liu2022inexact}
Liu W, Liu X, Chen X (2023{\natexlab{b}}) An inexact augmented {Lagrangian}
  algorithm for training leaky {ReLU} neural network with group sparsity.
  \emph{Journal of Machine Learning Research} 24(212):1--43.

\bibitem[{Lu et~al.(2020)Lu, Tsaknakis, Hong, \protect\BIBand{}
  Chen}]{lu2020hybrid}
Lu S, Tsaknakis I, Hong M, Chen Y (2020) Hybrid block successive approximation
  for one-sided non-convex min-max problems: algorithms and applications.
  \emph{IEEE Transactions on Signal Processing} 68:3676--3691.

\bibitem[{Melo \protect\BIBand{} Monteiro(2017)}]{melo2017iteration}
Melo JG, Monteiro RDC (2017) Iteration-complexity of a {Jacobi}-type
  non-{Euclidean} {ADMM} for multi-block linearly constrained nonconvex
  programs. \urlprefix\url{https://arxiv.org/pdf/1705.07229.pdf}.

\bibitem[{Melo et~al.(2024)Melo, Monteiro, \protect\BIBand{}
  Wang}]{melo2020iteration1}
Melo JG, Monteiro RDC, Wang H (2024) A proximal augmented lagrangian method for
  linearly constrained nonconvex composite optimization problems. \emph{Journal
  of Optimization Theory and Applications} 202(1):388--420.

\bibitem[{Murty \protect\BIBand{} Kabadi(1985)}]{murty1985some}
Murty KG, Kabadi SN (1985) Some {NP}-complete problems in quadratic and
  nonlinear programming. {Technical} report.
  \urlprefix\url{https://deepblue.lib.umich.edu/bitstream/handle/2027.42/6740/bam7058.0001.001.pdf?sequence=5}.

\bibitem[{Nemirovsky \protect\BIBand{} Yudin(1983)}]{nemirovskij1983problem}
Nemirovsky AS, Yudin DB (1983) \emph{Problem complexity and method efficiency
  in optimization} (Wiley-Interscience).

\bibitem[{Nesterov(2012)}]{nesterov2012make}
Nesterov Y (2012) How to make the gradients small. \emph{Optima. Mathematical
  Optimization Society Newsletter} (88):10--11.

\bibitem[{O’Neill \protect\BIBand{} Wright(2021)}]{o2021log}
O’Neill M, Wright SJ (2021) A log-barrier {Newton-CG} method for bound
  constrained optimization with complexity guarantees. \emph{IMA Journal of
  Numerical Analysis} 41(1):84--121.

\bibitem[{Sahin et~al.(2019)Sahin, Eftekhari, Alacaoglu, Latorre,
  \protect\BIBand{} Cevher}]{sahin2019inexact}
Sahin MF, Eftekhari A, Alacaoglu A, Latorre F, Cevher V (2019) An inexact
  augmented {Lagrangian} framework for nonconvex optimization with nonlinear
  constraints. \emph{Advances in Neural Information Processing Systems},
  13965--13977.

\bibitem[{Sun \protect\BIBand{} Hong(2019)}]{sun2019distributed}
Sun H, Hong M (2019) Distributed non-convex first-order optimization and
  information processing: lower complexity bounds and rate optimal algorithms.
  \emph{IEEE Transactions on Signal processing} 67(22):5912--5928.

\bibitem[{Thekumparampil et~al.(2019)Thekumparampil, Jain, Netrapalli,
  \protect\BIBand{} Oh}]{thekumparampil2019efficient}
Thekumparampil K, Jain P, Netrapalli P, Oh S (2019) Efficient algorithms for
  smooth minimax optimization. \emph{Advances in Neural Information Processing
  Systems} 32.

\bibitem[{Xiao et~al.(2023)Xiao, Liu, \protect\BIBand{}
  Toh}]{xiao2023dissolving}
Xiao N, Liu X, Toh KC (2023) Dissolving constraints for {Riemannian}
  optimization. \emph{Mathematics of Operations Research.}
  \urlprefix\url{https://doi.org/10.1287/moor.2023.1360}.

\bibitem[{Xie \protect\BIBand{} Wright(2021)}]{xie2021complexity}
Xie Y, Wright SJ (2021) Complexity of proximal augmented {Lagrangian} for
  nonconvex optimization with nonlinear equality constraints. \emph{Journal of
  Scientific Computing} 86:1--30.

\bibitem[{Xu(2017)}]{xu2017accelerated}
Xu Y (2017) Accelerated first-order primal-dual proximal methods for linearly
  constrained composite convex programming. \emph{SIAM Journal on Optimization}
  27(3):1459--1484.

\bibitem[{Xu(2021)}]{xu2021first}
Xu Y (2021) First-order methods for constrained convex programming based on
  linearized augmented {Lagrangian} function. \emph{INFORMS Journal on
  Optimization} 3(1):89--117.

\bibitem[{Xu(2024)}]{xu2024d-GDMax}
Xu Y (2024) Decentralized gradient descent maximization method for composite
  nonconvex strongly-concave minimax problems. \emph{SIAM Journal on
  Optimization} 34(1):1006--1044.

\bibitem[{Xu et~al.(2023)Xu, Zhang, Xu, \protect\BIBand{} Lan}]{xu2023unified}
Xu Z, Zhang H, Xu Y, Lan G (2023) A unified single-loop alternating gradient
  projection algorithm for nonconvex--concave and convex--nonconcave minimax
  problems. \emph{Mathematical Programming} 201(1):635--706.

\bibitem[{Yashtini(2022)}]{yashtini2022convergence}
Yashtini M (2022) Convergence and rate analysis of a proximal linearized {ADMM}
  for nonconvex nonsmooth optimization. \emph{Journal of Global Optimization}
  84(4):913--939.

\bibitem[{Zhang et~al.(2020)Zhang, Lin, Jegelka, Sra, \protect\BIBand{}
  Jadbabaie}]{zhang2020complexity}
Zhang J, Lin H, Jegelka S, Sra S, Jadbabaie S (2020) Complexity of finding
  stationary points of nonconvex nonsmooth functions. \emph{International
  Conference on Machine Learning}, 11173--11182 (PMLR).

\bibitem[{Zhang \protect\BIBand{} Luo(2020)}]{zhang2020proximal}
Zhang J, Luo Z (2020) A proximal alternating direction method of multiplier for
  linearly constrained nonconvex minimization. \emph{SIAM Journal on
  Optimization} 30(3):2272--2302.

\bibitem[{Zhang \protect\BIBand{} Luo(2022)}]{zhang2022global}
Zhang J, Luo Z (2022) A global dual error bound and its application to the
  analysis of linearly constrained nonconvex optimization. \emph{SIAM Journal
  on Optimization} 32(3):2319--2346.

\bibitem[{Zhang et~al.(2022)Zhang, Pu, \protect\BIBand{}
  Luo}]{zhang2022iteration}
Zhang J, Pu W, Luo Z (2022) On the iteration complexity of smoothed proximal
  {ALM} for nonconvex optimization problem with convex constraints,
  \urlprefix\url{https://arxiv.org/pdf/2207.06304.pdf}.

\bibitem[{Zhang et~al.(2021)Zhang, Yang, Guzm{\'a}n, Kiyavash,
  \protect\BIBand{} He}]{zhang2021complexity}
Zhang S, Yang J, Guzm{\'a}n C, Kiyavash N, He N (2021) The complexity of
  nonconvex-strongly-concave minimax optimization. \emph{Uncertainty in
  Artificial Intelligence}, 482--492 (PMLR).

\bibitem[{Zhu et~al.(2023)Zhu, Chen, Zhang, \protect\BIBand{}
  Wen}]{zhu2023optimal}
Zhu Z, Chen F, Zhang J, Wen Z (2023) On the optimal lower and upper complexity
  bounds for a class of composite optimization problems.
  \urlprefix\url{https://arxiv.org/pdf/2308.06470.pdf}.

\end{thebibliography}
